\documentclass[11pt]{article}

\usepackage{fullpage}

\usepackage{tipa}
\usepackage{dsfont}
\usepackage{tikz-cd}
\usepackage{adjustbox}

\usepackage{graphicx}

\usepackage{xr}
 
\usepackage{amsxtra}
\usepackage{amsmath}
\usepackage{amssymb}
\usepackage{amsfonts}
\usepackage[all,2cell]{xy}
\usepackage{mathrsfs}
\usepackage{amsthm}
\usepackage{enumitem}
\usepackage{hyperref}
\usepackage{subfigure}
\usepackage{mathabx}
\usepackage{euscript}
\usepackage[bbgreekl]{mathbbol}

\usepackage{tikz, tikz-3dplot, pgfplots}
\usetikzlibrary{decorations.markings, patterns, decorations.pathmorphing}

\tikzset{->-/.style={decoration={
markings,
mark=at position #1 with {\arrow{>}}},postaction={decorate}}}

\newtheoremstyle{example}{\topsep}{\topsep}%
     {}
     {}
     {\bfseries}
     {.}
     {2pt}
     {\thmname{#1}\thmnumber{ #2}\thmnote{ #3}}

\theoremstyle{example}
\newtheorem{exa}[equation]{Example}
\newtheorem{exas}[equation]{Examples}
\newtheorem{ex}[equation]{Example}
\newtheorem{rem}[equation]{Remark}
\newtheorem{rems}[equation]{Remarks}
\newtheorem{defi}[equation]{Definition}

\newtheorem{thm}[equation]{Theorem}
 \newtheorem{cor}[equation]{Corollary}
\newtheorem{lem}[equation]{Lemma}
\newtheorem{prop}[equation]{Proposition}
\newtheorem{conj}[equation]{Conjecture}

\setlist[enumerate,1]{label=(\arabic{*})}
\setlist[enumerate,2]{label=(\roman{*})}
\setlist[enumerate,3]{label=(\alph{*})}

\def\Ac{\mathcal{A}}
\def\Bc{\mathcal{B}}
\def\Cc{\mathcal{C}}
\def\Kc{\mathcal{K}}
\def\Dc{\mathcal{D}}

\def\Fc{\mathcal{F}}

\def\Lc{\mathcal{L}}

\def\Nc{\mathcal{N}}

\def\Oc{\mathcal{O}}

\def\Sc{\mathcal{S}}
\def\Tc{\mathcal{T}}

\def\Vc{\mathcal{V}}

\newcommand{\co}[1]{ {\langle #1 \rangle}}
\newcommand{\arr}[1]{ {\{ #1 \} }}

\newcommand*\circled[1]{\tikz[baseline=(char.base)]{  
            \node[shape=circle,draw,inner sep=2pt] (char) {#1};}}   
            
\def\<<{\langle {}\hskip -.1cm {}\langle}
\def\>>{\rangle \hskip -.1cm \rangle}
\def\={{\, \simeq\, }}

\def\A{ {\EuScript A}}

\def\Ac{\mathcal {A}}

\def\Adj{{\bf{Adj}}}
\def\adm{{\on{adm}}}

\def\B{ {\EuScript B}}
\def\be{\begin{equation}}

\def\C{{\EuScript C}}

\def\CC{\mathbb{C}}

\def\Cat{{\mathcal Cat}}
\def\Cat{{\EuScript Cat}}

\def\Cb{{\mathbf{C}}}

\def\Ch{\operatorname{Ch}}

 \def\Cof{{\on{Cof}}}
 
 \def\Coh{{\on{Coh}}}

\def\Cone{\operatorname{Cone}}

\def\Conv{{\on{Conv}}}

\def\cm{\langle m \rangle}
\def\cn{\langle n \rangle}

\def\Di{{\EuScript D}}
\def\D{\Di}
\def\Dc{{\mathcal{D}}}

\def\del{{\partial}}

 \def\dg{{\on{dg}}}

 \def\dgVect{{\on{dgVect}}}
 \def\drarr{\draw  [decoration={markings,mark=at position 0.9 with
{\arrow[scale=1.5,>=stealth]{>}}},postaction={decorate},
line width=.15mm]}

\def\E{{\EuScript E}}

\def\ee{\end{equation}}

\def\eps{{\varepsilon}}
\def\ev{\on{ev}}

\def\EXT{\operatorname{EXT}}
\def\Ext{\operatorname{{Ext}}}

\def\F{{ \EuScript F}}

\def\Fc{\mathcal{F}}

\def\Fib{ \operatorname{Fib}}
\def\fib{\operatorname{fib}}
\def\Filt{{\on{Filt}}}

\def\Fun{\operatorname{Fun}}

\def\gr{{\on{ gr}}}
\def\Grad{{\on{Grad}}}

\def\h{\operatorname{h}\!}

\def\HH{\mathbb{H}}

\def\Hom{{\operatorname{Hom}}}

\def\hra{\hookrightarrow}

\def\I{{\EuScript {I}}}

\def\Id{\on{Id}}

\def\J{ {\EuScript J}}

\def\K{ {\EuScript K}}
\def\Kc{{\mathcal K}}
\def\k{\mathbf k}

\def\Ker{{\on{Ker}}}

\def\Lc{{\mathcal{L}}}

\def\L{{\EuScript L}}

\def\lra{\longrightarrow}
\def\lla{\longleftarrow}
\def\llra{\longleftrightarrow}

\def\M{{\EuScript M}}
\def\Map{\operatorname{Map}}

\def\Mod{\on{Mod}}

\def\Mor{\on{Mor}}

\def\Nc{{\EuScript N}}
\def\N{\operatorname{N}}

\def\Ob{ \operatorname{Ob}}

\def\Oc{{\mathcal O}}
\def\ol{\overline}
\def\on{\operatorname}

\def\oo{{\infty}}
\def\op{{\operatorname{op}}}

\def\Perv{{\on{Perv}}}

\def\phi{{\varphi}}

\def\pt{ {\on{pt}}}

\def\Q{ {\EuScript Q}}

\def\Se{{\EuScript S}}

\def\Set{ {\operatorname{Set}}}
\def\Sh{{\on{Sh}}}

\def\Sph{\on{Sph}}
\def\sph{{\on{sph}}}

\def\St{{\EuScript S}t}

\def\U{ {\EuScript U}}

\def\V{ {\EuScript V}}
\def\Vect{ {\on{Vect}}}

\def\wt{\widetilde}

\def\X{{\EuScript X}}

\def\Y{{\EuScript Y}}

\def\ZZ{\mathbb{Z}}

\def\Sh{\on{Sh}}

\DeclareMathSymbol\DDelta\mathord{bbold}{"01}
\DeclareMathSymbol\GGamma\mathord{bbold}{"00}
\DeclareMathSymbol\SSigma\mathord{bbold}{'117}
\DeclareMathSymbol\LLambda\mathord{bbold}{'003}

\title{Spherical adjunctions of stable $\infty$-categories and the relative
S-construction}
\author{Tobias Dyckerhoff, Mikhail Kapranov, Vadim Schechtman, Yan Soibelman}

\begin{document}

\maketitle

\begin{abstract}
We develop the theory of semi-orthogonal decompositions and spherical functors in the framework of
stable $\oo$-categories. Building on this, we study the relative Waldhausen S-construction
$S_\bullet(F)$ of a spherical functor $F$ and equip it with a natural paracyclic structure
(``rotational symmetry''). 

This fulfills a part of the general program to provide a rigorous account of perverse schobers which
are (thus far conjectural) categorifications of perverse sheaves. Namely, in terms of our previous
identification of perverse sheaves on Riemann surfaces with Milnor sheaves, the relative $S$-construction
with its paracyclic symmetry amounts to a categorification of the stalks of a Milnor sheaf at a
singularity of the corresponding perverse sheaf. The action of the paracyclic rotation is a
categorical analog of the monodromy on the vanishing cycles of a perverse sheaf. 

Having this local categorification in mind, we may view the S-construction of a spherical functor as
defining a schober locally at a singularity. Each component $S_n(F)$ can be interpreted as a
partially wrapped Fukaya category of the disk with coefficients in the schober and with $n+1$ stops
at the boundary.
\end{abstract}

\tableofcontents

\addcontentsline{toc}{section}{Introduction}


\numberwithin{equation}{subsection}

\section*{Introduction}

This paper is the second, after \cite{dkss:1}, part of our project devoted to a systematic
development of the theory of perverse schobers which are (enhanced, triangulated) categorical
analogs of perverse sheaves, see \cite  {kapranov-schechtman:schobers} for a general discussion. It
can be read independently and provides results which may be of interest in their own way. Its
position in the general program can be explained as follows. 

Let $(X,N)$ be a {\em stratified surface}, i.e., a topological surface $X$, possibly with boundary,
together with a finite set $N$ of interior points. It gives rise to the abelian category
$\Perv(X,N)$ of perverse sheaves  on $X$ with singularities in $N$.  In \cite{dkss:1}, we
established a ``purely abelian'' description of  $\Perv(X,N)$ in terms of so-called {\em Milnor
sheaves}.  The latter are certain functors $\F$ which associate a vector space $\F(A,A')$ to any
pair (Milnor disk) consisting of a closed disk $A\subset X$ meeting $N$ in at most one point  and a
nonempty union of closed arcs $A'\subset \del A$.  Explicitly, for a perverse sheaf $F\in\Perv(X,N)$
the corresponding Milnor sheaf has the form $\F(A,A') = \HH^0(A,A'; F)$, the $0$th relative
hypercohomology (the only nontrivial hypercohomology!) with coefficients in $F$. 

This description, being purely abelian, is suitable for categorification. In this paper we study the
features of such categorification for the case when a disk $A$ is fixed, and only $A'\subset \del A$
varies.  In this case the categorical analog of $\F$ reduces to data of (enhanced triangulated)
categories $\F_n$, one for each $n\geq 0$ (corresponding to the case when $A'$ consists of $(n+1)$
arcs) which are connected by a {\em paracyclic structure} , i.e., by  an action of the category
$\Lambda_\oo$ covering  the  cyclic category $\Lambda$ of  Connes,  see  \cite{dkss:1}.  The
categorical analog of $F$ itself  is a perverse schober on the disk $A$, and $\F_n$ can be thought
as the Fukaya category of $A$ with coefficients in the schober and with supports (``stops'') at
$A'$. 
  
As proposed in \cite {kapranov-schechtman:schobers}, a  schober on a disk should be described  by
a spherical adjunction, i.e., an adjoint pair of functors 
\[
  \xymatrix{
  \C \ar@<.4ex>[r]^F&\D\ar@<.4ex>[l]^G
  }
\]
such that the cones of the unit and the counit of the adjunction are equivalences of categories,
see  \cite{AL:spherical} and compare with the description of perverse sheaves on the disk
\cite{beil:gluing, GGM}.  Our main discovery is that from the point of view of a spherical functor
$F$, the collection of categories $\F_n$ turns out to be nothing but  $S_\bullet(F)$, the {\em
relative Waldhausen $S$-construction} of $F$, cf.  \cite{waldhausen}. More precisely,
$S_\bullet(F)$ can be defined for any exact functor $F$  of stable $\oo$-categories and is a
simplicial object (in the category of stable $\oo$-categories). We prove that for a spherical
functor $F$, this simplicial structure is naturally extended to a paracyclic one, matching the
intuition related to a disk on a surface with a bunch of arcs in the boundary.  The action of the
center of $\Lambda_\oo$ corresponds, geometrically, to a categorical analog of the monodromy of a
perverse sheaf. 
    
Thus, from the point of view of perverse schobers on stratified surfaces, the structures studied in
this paper constitute the local data from which a general schober is glued. From the standpoint of
spherical functors however, these structures are absolutely  remarkable features which have no
conceivable explanation outside of the analogy with perverse sheaves.  In a paper to follow, we plan
to glue the structures defined here (associated naturally to disks), into global data on a
stratified surface. 

In this paper we adopt the language of stable $\oo$-categories as our preferred way of enhancing
triangulated categories. The main reason for this is the following. 

From the very beginning, the theory of spherical functors has been burdened by  technical
difficulties. Namely, as pointed out  in  \cite{AL:spherical}, the context of classical
triangulated categories is not sufficient since we need to take the cones of unit and counit natural
transformations. We need therefore some enhancement of the triangulated structure.   The most
immediate such enhancement (and the most frequently used in this context) is that of
pre-triangulated dg-categories \cite{BK:enhanced}. 

However, the concept of adjoint functors in the dg-enhanced triangulated context was always
considered rather subtle as it appears to require a tower of explicit higher coherence data, see
\cite{AL:Fourier, AL:spherical,  AL:bar}. For example, for adjoint dg-functors $F$ and $G$ the
``identification'' of the complexes $\Hom^\bullet(F(x), y)$ and $\Hom^\bullet(x, G(y))$ should have
the form of a {\em  quasi}-isomorphism and one immediately hits the questions like ``which
direction?" or ``what are the relations among such quasi-isomorphisms?" etc.  

The approach with $\oo$-categories allows one to avoid answering such questions explicitly while
retaining a sufficient pool of possible answers.  This is done via encoding $\oo$-functors between
$\oo$-categories in terms of (co)Cartesian fibrations over the interval $\Delta^1$, see
\cite[\S5.2]{lurie:htt}  and a discussion below in \S
\ref{subset:stab-triang}\ref{par:cocart-implicit} More precisely, an $\oo$-functor $F: \C \to\D$ of
$\oo$-categories gives rise to a coCartesian fibration $\Gamma(F)\to\Delta^1$ (covariant
Grothendieck construction).  The functor $F$ {\em has a right adjoint}, if this fibration is also
Cartesian and so {\em can be identified} with the contravariant Grothendieck construction of some
$\oo$-functor $G: \D\to\C$. Such $G$ ({\em an adjoint}, one among many) is defined uniquely up to a
contractible space of choices but  we can delay making this choice  until really necessary (or avoid
it altogether) working with $\Gamma(F)$ instead. This ``print-on-demand'' approach saves a lot of
effort in the end. 

For this reason we spend some time explaining, once for all future occasions, the details of the
concept of $\oo$-categorical adjunctions and its compatibility with other approaches to adjunctions,
esp. in the dg-setting. The outcome is that our approach (i.e., that of \cite{lurie:htt}) covers all
the other ones used previously, and so various earlier results retain their value in our context. 
However, while it is known that (pre-triangulated) dg-categories give rise to (stable)
$\oo$-categories \cite{lurie:ha, {faonte:nerve}}, the simplifications mentioned above do not seem to
simplify further by specializing to the dg-case. 
 
Another advantage of $\oo$-categories is that various auxiliary totalities (e.g. the collection of
all spherical adjunctions) readily make sense as $\oo$-categories but in the dg-setting they do not
form dg-categories, so some form of $\oo$-categories is needed anyway. 

In fact, the Grothendieck construction of a functor plays quite a prominent role in our study. In
addition to the above use of this concept, the procedure of {\em gluing a semiorthogonal
decomposition (SOD) from a gluing functor}, introduced in the dg-context in \cite{kuznetsov-lunts},
amounts to nothing but (the category of sections of) the Grothendieck construction. Further, the
first level $S_1(F)$ of the relative Waldhausen S-construction of a functor $F$, is again nothing
but the category of sections of the Grothendieck construction (i.e., the category with an SOD glued
via $F$).

\vskip .2cm

The paper consists of three chapters.

In Chapter 1 we develop the theory of spherical functors and spherical adjunctions in the framework
of stable $\oo$-categories (the  existing treatments in the literature use dg-categories).  \S
\ref{subset:stab-triang} is a general reminder on (stable) $\oo$-categories.  We pay particular
attention to the ``implicit'' way of defining $\oo$-functors in terms of (co)Cartesian fibrations
over $\Delta^1$.  In  \S \ref{subsec:adj-oo} we discuss, following \cite[\S 5.2]{lurie:htt}, the
$\oo$-categorical approach to adjunctions.  In \S \ref{subsec:comparison}, we give a comparison of
this approach with other approaches that are used in the literature: the one involving Fourier-Mukai
kernels in various dg-contexts \cite{AL:Fourier, AL:bar} and the one involving higher coherence data
\cite{{riehl-verity:2-cat}, {riehl-verity:yoneda},  {riehl-verity}}.   In \ref{subsec:sphad} we
define spherical adjunctions as those for which the cones (fibers) of the unit and the counit are
equivalences.  We give a few elementary examples related, in one way or another, to the concept of
sheaves on spheres.  They can be used to motivate a more sophisticated example in \S
\ref{sec:massey} which provides a family of spherical functors, one for each $n\geq 0$. The
importance of this family  is that it can be seen as the representing object of the spherical
S-construction  discussed later in \S \ref {subsec:rel-S-sph}. 

\vskip .2cm

Chapter \ref{sec:SOD} is devoted to the theory of semi-orthogonal decompositions (SOD) in the
framework of stable $\oo$-categories. While the concept of an SOD itself makes sense in the context
of classical triangulated categories, and we recall it in \S \ref {subsec:SOD-triang}, deeper
aspects such as gluing SODs via gluing functors, require an enhancement.  In \S \ref{sec:semi} we
develop these aspects when the enhancement is $\oo$-categorical. In line with our emphasis on
Cartesian and coCartesian fibrations as encoding data of $\oo$-functors, we introduce the concepts
of Cartesian and coCartesian SOD's which appear very naturally in the $\oo$-categorical context.
They lead to not one but two concepts of gluing functors which, when they exist, go in opposite
directions and it is convenient to distinguish them.  In \S \ref{subsec:admissibility}, we relate
the concept of admissibilty of an SOD (traditionally defined as the condition that various
orthogonals form SODs as well) with the existence of ($\oo$-categorical) adjoints of the gluing
functors. In \S \ref{subsec:mutation}, we study the mutation functor (identification of the left and
right orthogonal of an admissible subcategory) and upgrade it to a ``coordinate change functor'', a
categorical analog of a $2\times 2$ matrix, see \S \ref{subsec:mutation} \ref{par:mut-change} In \S
\ref{subsec:4-period} we adapt to the stable $\oo$-categorical context the interpretation, due
originally to Halpern-Leistner and Shipman \cite{halpern-shipman}, of spherical functors as gluing
functors of $4$-periodic SODs. 

\vskip .2cm

In Chapter \ref{sec:spher-S} we  study $S_\bullet(F)$, the Waldhausen S-construction of a spherical
functor $F$ and prove  that it has a natural structure a paracyclic $2$-Segal object in $\Cat_\oo$
with a certain prescribed pattern of adjunctions between face and degeneracy functors.  In \S
\ref{subsec:rel-S} we recall the paracyclic category $\Lambda_\oo$ and the original construction of
$S_\bullet(F)$ by Waldhausen \cite{waldhausen}.  In \S \ref{subsec:rel-S-sph} we interpret
$S_\bullet(F)$ in the particular case of  spherical $F$ in terms of the corresponding $4$-periodic
SOD.  This definition (we call it the {\em spherical S-construction}) makes the paracyclic symmetry
manifest.  Its equivalence with the original $S_\bullet(F)$, proved in Theorem \ref{thm:sphericalS},
is one of the forms of the main result of this paper. A less technical formulation, Theorem
\ref{thm:main}, says that $S_\bullet(F)$ admits a natural paracyclic structure.  In the final \S
\ref{subsec:sph-admis} we complement this result with additional properties holding for any
$S_\bullet(F)$ with $F$ spherical or not: the {\em $2$-Segal property} \cite{DK:HSS} and {\em
admissibility} (a certain pattern of adjunctions between faces and degeneracies)
\cite{dyckerhoff:DK}. We conjecture that all these properties taken together, characterize spherical
functors. More precisely, Conjecture \ref{conj:sphad=adm2Segpar} says that the $\oo$-category of
spherical adjunctions is equivalent to that of paracyclic $2$-Segal admissible stable
$\oo$-categories. 

\vskip .2cm

We are grateful to A. Bodzenta and A. Bondal for useful discussions. 
T.D. acknowledges the support
of the VolkswagenStiftung through the Lichtenberg Professorship Programme. The research of T.D. is
further supported by the Deutsche Forschungsgemeinschaft under Germany's Excellence Strategy -- EXC
2121 ``Quantum Universe'' -- 390833306. 
The research of M.K. was supported by World Premier
International Research Center Initiative (WPI Initiative), MEXT, Japan.  The research of Y.S. was
partially supported by Munson-Simu award of KSU.


\section{Spherical adjunctions}

\subsection{(Stable) $\oo$-categories and triangulated 
categories}\label{subset:stab-triang}

As in  \cite{dkss:1}, we use the language of $\oo$-categories \cite{lurie:htt}. In this section, we
review some basic notions for convenience of the reader and  to fix the notation. 

\paragraph{Basic reminder.} 
We denote by $\Set_\Delta$ the category of simplicial sets. A simplicial set $X$ will be sometimes denoted
$X_\bullet$ or $(X_n)_{n\geq 0}$ with $X_n$ being the set of $n$-simplices.  We denote by $X^\op$
the opposite simplicial set to $X$.

If $\k$ is a field, we denote by $\dgVect_\k$
  the category of cochain complexes (dg-vector spaces) over $\k$.  A cochain complex $V$ will be sometimes
  denoted $V^\bullet$ to indicate the grading. 
  
  \vskip .2cm
  
  By an $\oo$-{\em category} we mean a simplicial set  $\C$ satisfying the weak Kan condition
  (i.e., an  $(\oo,1)$-category, in a more precise terminology). 
  An $\oo$-functor between $\oo$-categories $\C,\D$  is simply  a morphism of simplicial sets. 
  Similarly, the {\em functor $\oo$-category} $\Fun(\C,\D)$ is simply the simplicial mapping space between $\C$ and $\D$
  as simplicial sets. 
  
    \vskip .2cm

  By an  {\em object} of $\C$  we mean a vertex, i.e.,
  an element of $\C_0$.  By a {\em morphism}
   of $\C$ we mean an edge, i.e., an element of $\C_1$, and we use the notation
  $\Hom_\C(x,y)$ for the set of edges of $\C$ beginning at $x$ and ending at $y$. 
  An ordinary small category $C$ can be considered as an $\oo$-category by forming the nerve of $C$, denoted $\N C$.
  In the opposite direction, any $\oo$-category $\C$ gives an ordinary category $\h\C$ with
  the same objects known as the {\em homotopy category} of $\C$. Its morphisms are certain
  equivalence classes of morphisms of $\C$. A morphism in $\C$ is called an {\em equivalence},
  if its image in $\h\C$ is an isomorphism. 
  
  \vskip .2cm
  
  For any two objects $x,y$ of an $\oo$-category $\C$ we have a Kan complex $\Map_\C(x,y)\in \Set_\Delta$
  (the {\em mapping space})
  whose set of vertices is $\Hom_\C(x,y)$. In fact, there are several different versions of this spaces,
  all homotopy equivalent to each other. They are functorial with respect to morphisms, i.e., we have 
  morphisms of simplicial sets (with $\Hom_\C$ being a discrete simplicial set). 
  \be\label{eq:map-functor}
  \Hom_\C(y,z) \times \Map_\C (x,y) \lra \Map_\C(x,z), \quad \Map_\C(y,z) \times\Hom_\C(x,y) \lra
  \Map_\C(x,z).
  \ee
  However, there is no preferred  composition law
  \[
   \Map_\C(y,z) \times \Map_\C (x,y) \lra \Map_\C(x,z)
  \]
   for
  such mapping spaces although one can define it up to contractible space of choices,  see \cite[1.2.2]{lurie:htt}. 
  The mapping spaces can be used to define $\oo$-categorical initial and final objects  (e.g., $x$ is initial,  if each $\Map_\C(x,y)$ is contractible) and therefore other types of universal constructions such as
  limits $\varprojlim$ and colimits $\varinjlim$. 
  
  \vskip .2cm

 In the opposite direction, let $\Cc$ be a category enriched in Kan complexes.  Then we can
 associate to $\Cc$ an $\oo$-category $N^\Delta \Cc$ by means of the {\em simplicial nerve
 construction}, see \cite[1.1.5]{lurie:htt} and references therein. 
  
        \vskip .2cm
  
  By a {\em dg-category} we will mean a category enriched in $\dgVect_\k$. Given a dg-category $\Ac$,
  one can associate to it an $\oo$-category $N^\dg \Ac$ by means of the {\em dg-nerve construction} 
  (the simplicial set formed by ``Sugawara simplices'' \cite[App.2]{HS}) \cite[1.3.1]{lurie:ha} \cite{faonte:nerve}. The homotopy category $\h N^\dg A$ is identified with the $0$th
  cohomology category $H^0 \Ac$. 
  More generally, any $A_\oo$-category gives an $\oo$-category by means of the {\em $A_\oo$-nerve construction}
   \cite{faonte:nerve}.    In this way various Fukaya $A_\oo$-categories appearing in symplectic geometry
 \cite{seidel}   can be included in the $\oo$-categorical formalism.  

  
   \paragraph{The $\oo$-category of $\oo$-categories} \label{par:oo-2-level}
   
The simplicial nerve construction provides the fundamental object: $\Cat_\oo$, the {\em  $\oo$-category of
all (small) $\oo$-categories}. It is defined as the simplicial nerve of the simplicially enriched category
$\Cat_\oo^\Delta$ whose objects are $\oo$-categories as above and $\Map(\C,\D)$ is
the largest Kan subcomplex in $\Fun(\C,\D)$, see \cite[Def. 2.0.0.2]{lurie:htt}. 
    
       
           \vskip .2cm

This approach is very convenient and  will be used later in this paper  to speak about
$\oo$-categories formed by related data such as, $\oo$-functors,  adjunctions of $\oo$-categories
etc. It will be eventually used to speak about $\oo$-categories of perverse schobers. 
 
 \vskip .2cm
 
Sometimes, a more general approach is necessary. Indeed, the term ``$\oo$-category'' adopted in
\cite{lurie:htt} and here, stands for ``$(\oo,1)$-category'', with the ``$1$'' signifying that all
$p$-morphisms with $p>1$ are weakly invertible. So this concept does not include all classical, even
strict, $2$-categories. In particular, it  cannot describe  $\Cat$,  the classical {\em $2$-category
of all usual categories} (objects: categories, morphisms: functors, $2$-morphisms: natural
transformations). One can form the corresponding part of $\Cat_\oo$ as defined above,  but it
captures only {\em invertible} natural transformations. 
 
  \vskip .2cm
 
 The more general structure possessed by  ``the category of all $\oo$-categories''
 (and capturing that on $\Cat$) is that of    $(\oo,2)$-category. Let us give
 a simplified (partially strict) definition, following \cite{riehl-verity}.
 See \cite{lurie:ha} for a ``fully lax'' definition. 
 
 \begin{defi}\label{def:oo-2}
 An {\em $(\oo,2)$-category} is a simplicially enriched category $\Cc$ such that
 for any objects $c,d\in \Cc$ the simplicial set $\Hom_\Cc(c,d)$ is an $\oo$-category. 
 \end{defi}
  
  This is analogous to defining a strict classical $2$-category as a category
  enriched in the category of categories.      In fact, any strict $2$-category $C$
  can be considered as an $(\oo,2)$-category $\N^\Hom C$ by taking the nerves of all the
  categories $\Hom_C(c,d)$. 
  
    \vskip .2cm
    
   Further,  we have the $(\oo,2)$-category $\Cat_\oo^{(2)}$ of all
   $\oo$-categories. Its objects are all $\oo$-categories, and the
   simplicial Hom-sets are the functor $\oo$-categories $\Fun(\C,\D)$.  
   This $(\oo,2)$-category refines the $\oo$-category $\Cat_\oo$ defined above.

   
   \paragraph{Stable $\oo$-categories.}\label{par:stable-oo}
   We refer to \cite{lurie:stable, lurie:ha} for a general treatment. It is convenient to take,
   as the definition, the following reformulation
     \cite[Prop. 1.1.3.4]{lurie:ha}: 
   \begin{defi}
   An $\oo$-category $\C$ is called {\em stable}, if:
   \begin{itemize}
   \item[(S1)] $\C$ has finite limits and colimits, in particular, an initial and  and a final object.
   
   \item[(S2)] $\C$ is pointed, i.e., initial and final objects coincide. Such objects are called 
   {\em zero objects} and denoted $0$.
   
   \item [(S3)] A commutative square (i.e., an $\oo$-functor $\Delta^1\times\Delta^1\to\C$)
   \[
   \xymatrix{
   X\ar[d] \ar[r]&Y\ar[d]
   \\
   T\ar[r]& Z
   }
   \]
   is Cartesian (i.e., $X$ is equivalent to $\varprojlim\, \{ T\to Z \leftarrow Y\}$) if an only if it is coCartesian
   (i.e., $\varinjlim\, \{ T\leftarrow X\to Y\}$  is equivalent to $Z$). Such squares are called
   {\em biCartesian}. 
   \end{itemize} 
   \end{defi}
   
   If $\C$ is stable, then the category $\h\C$ is triangulated in the classical sense of
   Verdier, see  \cite[Thm. 3.11]{lurie:stable}  \cite[Thm. 1.1.2.14]{lurie:ha}. 
   Distinguished triangles in $\h\C$ come
   from biCartesian squares as in (S3) with $T=0$. Such squares  are called {\em triangles} in $\C$. 
   Let $f: \C\to\D$ be an $\oo$-functor of stable $\oo$-categories. It is known \cite[Prop. 1.1.4.1]{lurie:ha}
   that $f$ being {\em left exact} (i.e., commuting with finite limits) is equivalent to $f$ being
   {\em right exact} (i.e., commuting with finite colimits). If this is the case, we call $f$ {\em exact}. 
   If $f$ is  exact, then the induced functor $\h f: \h\C\to\h\D$ is  triangulated,
   i.e., preserves the shift and the class of distinguished triangles. 
   In this way stable $\oo$-categories provide ``enhancements" of classical triangulated categories. 
   
   \vskip .2cm

   Unlike classical triangulated categories, in a stable $\oo$-category $\C$ we can speak about
  (homotopy)  canonical cones of morphisms by fitting them into triangles in $\C$. 
     In fact, it is convenient to separate this concept into
   two: the {\em fiber} of a morphism $f: X\to Y$, denoted $\Fib(f)$ and the {\em cofiber}, or
   {\em cone}, denoted $\Cof(f)$ as well as $\Cone(f)$. They are defined as
   \be\label{eq:(co)fiber}
   \begin{gathered}
   \Fib(f) \,=\, \varprojlim \, \{ X\buildrel f\over\to Y \leftarrow 0\},
    \quad \Cone(f)=\Cof(f) \,=\,   \varinjlim \, \{ 0\leftarrow X \buildrel f\over\to Y\}. 
   \end{gathered}
   \ee
   As any limits and colimits in an $\oo$-category, the fiber and cofiber are defined uniquely up
   to a contractible choice.
   The suspension and desuspension functors in $\C$ (corresponding to the shift 
   by $\pm 1$ in the triangulated category $\h\C$) are defined as
   \[
   X[1] \,=\,\Cof\{X\to 0\}, \quad X[-1] \,=\,\Fib\{ 0\to X\}. 
   \]
   We further have an equivalence $\Fib(f) \simeq \Cof (f)[-1]$.
     Another advantage of stable $\oo$-categories is the following:
   
   \begin{prop}
   (a) If $\A, \C$ are $\oo$-categories and $\C$ is stable, then $\Fun(\A, \C)$ is stable.
   
   \vskip .2cm
   
   (b) Let $\A$ and $\C$ be stable $\oo$-categories. The fiber and cofiber of any morphism
   of exact $\oo$-functors $\A\to \C$ is exact. 
   \end{prop}
   
   \noindent{\sl Proof:} (a) is a particular case of \cite[Prop. 4.1]{lurie:stable} (in fact, one can take
   for $\A$ any simplicial set, not necessarily an $\oo$-category). To show (b), note that
   projective limit is a left exact operation, so $\Fib(f)$, being a projective limit of left
   exact functors, is left exact and hence is exact since $\A$ and $\C$ are stable.
   Similarly, $\Cof(f)$, being an inductive limit of right exact functors, is right exact hence exact. \qed

   The {\em $\oo$-category of stable $\oo$-categories}  $\St$ is defined as the simplicial subset
   in $\Cat_\oo$ formed by simplices whose vertices correspond to stable $\oo$-categories and
   edges to exact $\oo$-functors.
   
   \vskip .2cm
   
   Another type of enhancements used in the literature is that of {\em dg-enhancements},
   given by  pre-triangulated dg-categories \cite{BK:enhanced}. It can be included
   into the $\oo$-categorical picture: as shown in \cite{faonte:nerve}, for a pre-triangulated
   dg-category $C$, its dg-nerve $N^\dg C$ is a stable $\oo$-category, which gives the same
   triangulated category $\h N^\dg C = H^0 C$.  Note, however, an interesting difference
    with
   the $\oo$-categorical case: 
 for  any dg-functor $f: C\to D$ between pre-triangulated dg-categories, the induced functor
 $H^0 f: H^0 C\to H^0 D$ is triangulated,  there is no need to impose any further conditions on $f$. This is
 because the concept of an distinguished triangle, at the dg-level, can be encoded by
 a dg-functor from some particular dg-category with three objects, see \cite[\S I.3]{seidel}


\paragraph{Cartesian and coCartesian fibrations.}
Let $C$ be an ordinary category. According to the classical theory of Grothendieck \cite[Ch.3]{fantechi},
 fibered categories $X\to C$ over $C$ provide an alternative, implicit way 
 of encoding contravariant pseudo-functors $F: C\to\Cat$
from $C$ to the $2$-category of categories.  Covariant 
pseudo-functors are similarly encoded by cofibered categories. 

An  $\oo$-categorical generalization of this point of view  is given by Lurie's theory of
Cartesian and  coCartesian fibrations of $\oo$-categories
\cite[\S 2.4]{lurie:htt}.  Let us recall
the necessary concepts at the level of generality we need.

\begin{defi}\label{def:cart-edge}
(a) Let  $p: \X\to\C$ be an $\oo$-functor of $\oo$-categories.
An edge (morphism) $f: x\to y$ in $\X$ is called {\em $p$-Cartesian}, if the natural map
\[
\X_{/f} \lra \X_{/y} \times_{\C_{/p(y)}} \C_{/p(f)}
\]
is a trivial Kan fibration of simplicial sets. 

\vskip .2cm

(b) An edge $f$ is called $p$-{\em coCartesian}, if it corresponds to a Cartesian edge for
the opposite $\oo$-functor $p^\op: \X^\op\to\C^\op$. 

\end{defi}

Here $\X_{/y}$ is the over($\oo$-)category whose objects are objects in $\X$ over $y$, and similarly
for $\C_{/p(y)}$. The notation $\X_{/f}$ means the over($\oo$-)category whose objects are morphisms
in $\X$ over $f$, and similarly for $\C_{/p(f)}$, see \cite[Def. 1.2.9.2]{lurie:htt} for the general concept. 

\vskip .2cm

Explicitly,  $f$ being $p$-coCartesian can be  expressed  by a condition involving
undercategories instead of overcategories
and the source rather than the target of $f$. 

Notationally, for  an edge $f: a \to b$ in $\X$, we will write $a \overset{!}{\to} b$ if $f$ is coCartesian and $a \overset{*}{\to} b$ if $f$ is
Cartesian.

\begin{defi}(\cite{lurie:htt}, Def. 2.4.2.1 + Prop. 2.3.1.5) \label{def:cart-fib}
Let $\C$ be (the nerve of) an ordinary category and  $p: \X\to\C$ be an $\oo$-functor of $\oo$-categories.
\begin{itemize}
\item[(a)]
 We say that $p$ is a {\em Cartesian fibration}, if for every edge $g: c\to d$ of $\C$ and any vertex (object)
 $\wt d$ of $\X$ with $p(\wt d)=d$, there is a $p$-Cartesian edge $\wt g: \wt c\to\wt d$ with
 $p(\wt g)=g$. 
 
 \item[(b)] We say that $p$ is a {\em coCartesian fibration}, if $p^\op: \X^\op\to \C^\op$ is a Cartesian fibration. 
 \end{itemize}
\end{defi}

\begin{rem}\label{rem:marked}
An important technical tool for study of Cartesian (say) fibrations over $\C$ is Lurie's model category
$(\Set^+_\Delta)_{/\C}$ whose objects are {\em marked simplicial sets over $\C$}, i.e.,
simplicial sets $S\to\C$ over $\C$ with a distinguished set of edges called {\em marked edges}
\cite[\S 3.1]{lurie:htt}.
Cartesian fibrations (with Cartesian edges as marked edges) are characterized as fibrant objects
in $(\Set^+_\Delta)_{/\C}$, so one can use fibrant resolutions etc. We do not go into further details
in this reminder but will use this technique occasionally in the proofs. 
\end{rem}

  Examples of (co)Cartesian fibrations are provided by the $\oo$-categorical version of the
 Grothendieck construction (relative nerve)  \cite[Def. 3.2.5.2]{lurie:htt}. 
Let $C$ be an ordinary small category with nerve $\C=\N\C$ and $q: C \to \Set_{\Delta}$ be a functor such that, for every
object $c$ of $C$, the simplicial set $q(c)$ is an $\infty$-category. Let also $[n]$ be the poset
${0<1\cdots <n}$ considered as a category,  with the nerve $\Delta^n$, the $n$-simplex.

\begin{defi}\label{def:groth}
 (a) The {\em covariant Grothendieck construction} of $q$ is the simplicial set $\Gamma(q)$
	defined as follows. An $n$-simplex in $\Gamma(q)$ consists of
	\begin{enumerate}
		\item a functor $\sigma: [n] \to C$,
		\item for every $I \subset [n]$, a map $\Delta^I \to q(\sigma(\max(I)))$,
	\end{enumerate}
	such that, for every $J \subset I \subset [n]$, the diagram 
	\[
		\begin{tikzcd}
			\Delta^J \ar{r}\ar{d} & q(\sigma(\max(J))) \ar{d}\\
			\Delta^I \ar{r} & q(\sigma(\max(I))) 
		\end{tikzcd}
	\]
	commutes. There is a natural forgetful map ($\oo$-functor) $\pi: \Gamma(q) \to \C=\N(C)$.
	
	\vskip .2cm
	
	(b) The {\em contravariant Grothendieck construction} $\chi(q)$ of $q$ is the simplicial set
$\Gamma((-)^\op \circ q)^{\op}$.   It comes equipped with a natural forgetful map $\pi: \chi(q) \to
\C^\op$.
\end{defi}

As shown in \cite[\S 3.2.5]{lurie:htt}, the functor $\Gamma(q) \to \C$ is a coCartesian fibration
while  the functor $\chi(q) \to \C^{\op}$ is a Cartesian fibration.

Conversely, see   \cite[Th.3.2.0.1]{lurie:htt} any coCartesian fibration $p: \X\to\C$ gives rise to a {\em classifying $\oo$-functor}
$Q: C \to\Cat_\oo$ defined uniquely up to a contractible space of
choices. For   an object $c$ of $C$,  the value $\Q(c)$ 
is the preimage $p^{-1}(c)\subset\X$ in the sense of simplicial sets.

 If $\X=\Gamma(q) \to \C^\op$
 is the covariant Grothendieck construction of a functor $q$ as in Definition \ref{def:groth},
  then $\N(q)$ is a classifying $\oo$-functor.  
 
 Dually, any Cartesian fibration $\X$ over $\C$ gives rise to a classifying functor $Q: \C^\op\to \Cat_\oo$.
 Again, the value $Q(c)$ on an object $c\in C$ is the preimage $p^{-1}(c)\subset \X$. 
 If $\X=\chi(q)\to\C^\op$ is the contravariant Grothendieck construction as above, then $\N(q)$
 can be taken as a classifying $\oo$-functor. 
 
 \begin{rem}\label{rem:subtlety}
 If $\C=\N C$ is the nerve of an ordinary category (as we assume),
  then an $\oo$-functor $Q: \C\to\Cat_\oo$
 is a more general datum than an input of Definition  \ref{def:groth}, i.e., an ordinary functor $q: C\to\Set_\Delta$ with each $q(c)$ an $\oo$-category. More precisely, $Q$  does give for each object $c\in C$ an $\oo$-category
 $Q(c)$ and for each morphsim $u: c\to c'$ an $\oo$-functor $Q(u): Q(c)\to Q(c')$,
 but the compatibility of these data with composition holds not strictly but up to higher homotopies
 which are hidden in the definition of $\Cat_\oo$ as a simplicial nerve. For example,
 for a composable pair of morphisms $u,v$ in $C$ we are given a $2$-simplex with edges $Q(uv), Q(u), Q(v)$ and so on. 
  So $Q$ is an analog of a pseudo-functor $C\to\Cat$ in the classical sense. 
 \end{rem}

 \paragraph{(Co)Cartesian fibrations over $\Delta^1$: implicit $\oo$-functors.}
 \label{par:cocart-implicit}
 
 In general, the construction of a classifying $\oo$-functor of a (co)Cartesian fibration is quite involved
 and uses the model category  of marked simplicial sets over $\C$ (see Remark \ref{rem:marked}). 
 Let us work out explicitly the simplest example, which will be of particular importance for
 the discussion of adjoint functors later in the paper. 
 
 \begin{ex}[(Grothendieck construction for a single $\oo$-functor)]\label{ex:groth-single}
 Let $C=[1]=\{0\to 1\}$ and $\C=\N C=\Delta^1$ be the  ($\oo$-) category  with two objects $0,1$ and one
non-trivial morphism. 
A covariant $\oo$-functor $q: \C\to\Cat_\oo$ reduces to a single $\oo$-functor
$F: \A=q(0) \to \B=q(1)$ of $\oo$-categories,   the subtlety
 from Remark \ref{rem:subtlety} not appearing here.  In this case we denote the covariant Grothendieck construction of
$q$ by $\Gamma(F)$.

 Similarly, a contravariant functor $q: \C\to\Cat_\oo$
 is the same as a single $\oo$-functor $G: \B \to\A $.  In this case we denote the contravariant Grothendieck
 construction of $q$ by $\chi(G)$.

 Let us describe the covariant Grothendieck construction $\Gamma(F)$ associated
 to  an $\oo$-functor  $F:\A\to\B$ explicitly.   
 For any integers $0\leq m\leq n$ we have the standard embedding $[m]\hra [n]$. 
 For any simplicial set $X$ we will denote $\del_{[m]}:  X_n \to X_m$ 
 the iterated face map associated to this embedding.
 Specializing Definition \ref{def:groth} to our case, we see that an $n$-simplex $s\in \Gamma(F)_n$
  is a datum of:
 \begin{itemize}
 \item[(1)] An integer $-1\leq m\leq n$.
 
 \item[(2)] An  $n$-simplex $s_\B\in \B_n$ and, in case $m\geq 0$, 
 an $m$-simplex $s_\A\in \A_m$  such that
 $\del_{[m]} s_\B = F(s_\A)$.
  \end{itemize}
  The projection $p: \Gamma(F)\to\Delta^1$ sends an $n$-simplex $(m, s_\A, s_\B)$
  to (the degeneration of) $0$, if $m=n$, to (the degeneration of) $1$, if $m=-1$ and to
  (the degeneration of) the $1$-dimensional simplex of $\Delta^1$ otherwise.
  
   In particular, a {\em section} of $p$, i.e., an edge (morphism)  projecting to the edge of $\Delta^1$,
  is a datum of an object $a=s_\A \in\A$, $b\in \B$ and a morphism $s_\B: F(a)\to b$ in $\B$. 
  Among these,   coCartesian edges of $\Gamma(F)$ 
  are distinguished by the condition that $s_\B$ is an equivalence. In particular, we can take
  for $a$ any object of $\A$ and put $s_\B=\Id: F(a)\to F(a)$.

  We see that any object in $\A= p^{-1}(0)$ is a source of a coCartesian edge, so $p$ is a coCartesian fibration,
  by Definition \ref{def:cart-fib}(b). 
 \end{ex}

 \begin{ex}[(Implicit $\oo$-functors)]\label{ex:implicit}
 Conversely, 
 given a Cartesian (resp. coCartesian) fibration $p: \X\to\Delta^1$, we have two $\oo$-categories
 $\A=p^{-1}(0)$, $\B=p^{-1}(1)$ (the preimages as simplicial sets).  A classifying $\oo$-functor for $p$
 reduces to a single $\oo$-functor $G: \B\to \A$ (resp. $F: \A\to\B$).
   We can say that $G$ (resp. $F$)
 is an {\em implicit $\oo$-functor} defined by a Cartesian (resp. coCartesian) fibration $p$. 
 Let us explain how to construct $G$ for a Cartesian fibration $p$ (the case of a coCartesian fibration beingn dual). 
 We define $G$ (a morphism of simplicial sets)
 inductively on simplices of dimensions $0,1,2,\cdots$. 
 \begin{itemize}
 \item[(0)] Let $y$ be an object (vertex) of $\B$, so $p(y)=1$. 
  Let $g: 0\to 1$ be the nontrivial edge of $\Delta^1$.
  By Definition \ref{def:cart-edge}(a), 
 there is a Cartesian edge $f_y: x\to y$ in $\X$ covering $g$, so that $x$ is a vertex of $\A$.
 Choosing such an edge, we put $G(y)=x$, thus defining $G$ on $0$-simplices. 
 
 \item[(1)]  Let $v: y_0\to y_1$ be a morphism (edge) in $\B$. Let $x_i=G(y_i)$ be as constructed before. 
 Because $\X$ is an $\oo$-category
 (weak Kan complex) we can form a pseudo-composition of $v$ and $f_{y_0}$, filling a horn of type
 $\Lambda^2_1\subset \Delta^2$ , see \cite[Def. 1.1.2.1]{lurie:htt} for notation. This gives a
 $2$-simplex (triangle) denoted \circled{1}  in Fig. \ref{fig:implicit1}. The third side of this
 triangle is a morphism which we denote $c: x_0= G(y_0) \to y_1$. 
 
 \begin{figure}[h]
 \centering
 \begin{tikzpicture}[scale=0.3]
 
 \node at (-6,0){\small$\bullet$}; 
  \node at ( 6,0){\small$\bullet$}; 
  
  \draw [decoration={markings,mark=at position 0.7 with
{\arrow[scale=1.5,>=stealth]{>}}},postaction={decorate},
line width = .3mm]  (-6,0) -- (6,0) ;

\node at (-7,0){$0$}; 
\node at (7,0){$1$}; 
\node at (0,-1) {$g$}; 

\node at (6,4) {\small$\bullet$}; 
\node at (6,9) {\small$\bullet$}; 
\node at (-6,4) {\small$\bullet$}; 
\node at (-6,9) {\small$\bullet$}; 

\drarr (6,9) -- (6,4); 
\drarr (-6,9) -- (-6,4); 
\drarr (-6,9) -- (6,9); 
\drarr (-6,4) -- (6,4); 
 \drarr (-6,9) -- (6,4); 
 
 \node at (2,7.5){\circled{1}}; 
 \node at (-2, 5.5){\circled{2}}; 
 
 \node at (7,9){$y_0$}; 
 \node at (7,4){$y_1$};  
 \node at (-7,9){$x_0$}; 
 \node at (-7,4){$x_1$}; 
 
 \node at (7, 6.5){$v$}; 
  \node at (-7, 6.5){$u$}; 
  
  \node at (0,10){$f_{y_0}$}; 
 \node at (0,3){$f_{y_1}$}; 
 \node at (0.5, 7){$c$}; 
 
 \draw (6, 6.5) ellipse (2.5cm and 5cm); 
  \draw (-6, 6.5) ellipse (2.5cm and 5cm); 
 
 \node at (-8.5, 11){$\A$}; 
  \node at (8.5, 11){$\B$}; 
 
 \end{tikzpicture}
 
 \caption{An implicit $\oo$-functor from a Cartesian fibration: values on morphisms.}\label{fig:implicit1}
 \end{figure}
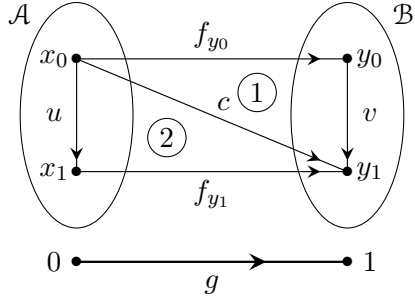

 Next, the condition for $f_{y_1}$ to be a Cartesian edge reduces in our particular case to saying that
 $\X_{/f_{y_1}} \to \X_{/y_1}$ is a trivial Kan fibration, in particular, it is surjective on $0$-simplices.
 Now, $c$ is an object ($0$-simplex)  of $\X_{y_1}$, while a $0$-simplex of $\X_{f_{y_1}}$ covering $c$
 is, by definition, a triangle of the form \circled{2}  in Fig. \ref{fig:implicit1}. 
 Denote $u: x_0\to x_1$ its third side. We put $G(v)=u$, thus defining $G$ on $1$-simplices. 
 
 \item[(2)] Next, let $\tau$ be a $2$-simplex in $\B$ with vertices $y_0,y_1,y_2$
 and edges $v_{ij}: y_i\to y_j$, $i<j$. 
  Let $x_i= G(y_i)$ and $u_{ij}= G(v_{ij})$, already constructed. 

 The previous choices give us triangles filling a part of the boundary of the $3$-dimensional prism 
 $\Delta^1\times\Delta^2$, see Fig. \ref{fig:implicit2}.  More precisely, $\tau$ fills one (the right) foundation of the prism and the
 triangles of the form \circled{1} and \circled{2} for the $v_{ij}$ fill the walls of the prism (which are
 therefore triangulated), but the left foundation is not filled. 
 Now, we have the standard triangulation of the prism, compatible with the triangulations of the walls.
 It has 3 tetrahedra which we denote by the sets of their vertices:
 $[x_0y_0y_1y_2]$, $[x_0 x_1 y_1 y_2]$ and $[x_0 x_1 x_2 y_2]$. We fill these tetrahedra (by $3$-simplices
 of $\X$) as follows.
 
 \begin{figure}[h]
 \centering
 \begin{tikzpicture}[scale=0.4]
 
 \node at (10,0){\small$\bullet$}; 
  \node at (10,9){\small$\bullet$}; 
 \node at (6,4){\small$\bullet$};  
  \node at (-10,0){\small$\bullet$}; 
 \node at (-10,9){\small$\bullet$};  
  \node at (-14,4){\small$\bullet$}; 
  
  \draw [decoration={markings,mark=at position 0.7 with
{\arrow[scale=1.5,>=stealth]{>}}},postaction={decorate},
line width = .3mm]  (-10,0) -- (10,0) ;  
 
  \draw [decoration={markings,mark=at position 0.7 with
{\arrow[scale=1.5,>=stealth]{>}}},postaction={decorate},
line width = .3mm]  (-10,9) -- (10,9) ;   

  \draw [decoration={markings,mark=at position 0.7 with
{\arrow[scale=1.5,>=stealth]{>}}},postaction={decorate},
line width = .3mm]  (10,9) -- (6,4) ;  

  \draw [decoration={markings,mark=at position 0.7 with
{\arrow[scale=1.5,>=stealth]{>}}},postaction={decorate},
line width = .3mm]   (6,4)-- (10,0) ;  

  \draw [decoration={markings,mark=at position 0.7 with
{\arrow[scale=1.5,>=stealth]{>}}},postaction={decorate},
line width = .3mm]  (10,9) -- (10,0) ;  

  \draw [decoration={markings,mark=at position 0.7 with
{\arrow[scale=1.5,>=stealth]{>}}},postaction={decorate},
line width = .3mm]  (-14,4) -- (6,4) ;  

  \draw [decoration={markings,mark=at position 0.7 with
{\arrow[scale=1.5,>=stealth]{>}}},postaction={decorate},
line width = .3mm]  (-10,9) -- (-14,4) ;  

  \draw [decoration={markings,mark=at position 0.7 with
{\arrow[scale=1.5,>=stealth]{>}}},postaction={decorate},
line width = .3mm]  (-14,4) -- (-10,0) ;  

  \draw [dashed, decoration={markings,mark=at position 0.7 with
{\arrow[scale=1.5,>=stealth]{>}}},postaction={decorate},
line width = .35mm]  (-10,9) -- (-10,0) ;  

\drarr (-14,4) -- (10,0); 
\drarr (-10,9) -- (6,4); 

  \draw [dashed, decoration={markings,mark=at position 0.7 with
{\arrow[scale=1.5,>=stealth]{>}}},postaction={decorate},
line width = .15mm]  (-10,9) -- (10,0) ;

\node at (11,9){$y_0$}; 
\node at (7,4){$y_1$}; 
\node at (11,0){$y_2$}; 
\node at (-11,9){$x_0$}; 
\node at (-15,4){$x_1$}; 
\node at (-11,0){$x_2$}; 

\node at (11, 4.5){$v_{02}$}; 
\node at (7,7){$v_{01}$}; 
\node at (8,3){$v_{12}$}; 

\node at (0, 10){$f_{y_0}$}; 
\node at (0,-1){$f_{y_2}$}; 
\node at (-3,4.8){$f_{y_1}$}; 

\node at (-13,7){$u_{01}$}; 
\node at (-13,2){$u_{12}$}; 
\node at (-9,5){$u_{02}$}; 

\node at (8.5,5){\circled{$\tau$}}; 
 
\end{tikzpicture}
\caption{An implicit $\oo$-functor from a Cartesian fibration: values on $2$-simplices.}\label{fig:implicit2}
 \end{figure}
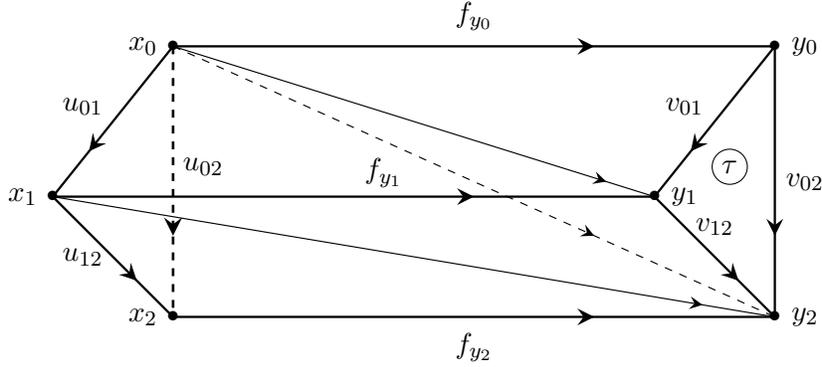

 First, out of the desired $[x_0y_0y_1y_2]$ we already have the horn $\Lambda^3_1$, so we can fill it by
 the weak Kan property of $\M$. After this, we have a horn $\Lambda^3_2$ out of the desired $[x_0 x_1 y_1 y_2]$,
 so we fill it as well. After this, the triangle $[x_0 x_1 y_2]$ represents a $1$-simplex in $\X_{/y_2}$
 whose vertices $[x_0y_1]$ and $[x_1y_2]$ are lifted to $\X_{/f_{y_2}}$. Now, as before, the fact
 that $f_{y_2}$ is a Cartesian edge, means that $\X_{/f_{y_2}} \to \X_{/y_2}$ is a trivial Kan fibration. 
 Therefore any map of $\Delta^1$ to the $ \X_{/y_2}$ with a lift over $\del \Delta^1$ to $\X_{/f_{y_2}}$
 can be lifted to a map of  $\Delta^1$  to $\X_{/f_{y_2}}$, i.e., to a $1$-simplex of  $\X_{/f_{y_2}}$.
 Such a $1$-simplex is, by definition \cite[Prop. 1.2.9.2]{lurie:htt} nothing but a tetrahedron ($3$-simplex of $\X$)
 of the desired shape $[x_0 x_1 x_2 y_2]$.
 The $4$th side of this tetrahedron is a $2$-simplex $\sigma$ in $\A$ with sides $u_{01}, u_{12}, u_{02}$. 
 We put $G(\tau)=\sigma$, thus defining $G$ on $2$-simplices. 
 
 \item[(...)] 
 
 \item[($n$)] Suppose that $G$ is defined on simplices of dimension $<n$ and $\tau$ is an $n$-simplex
 of $\B$, with vertices $y_0,\cdots, y_n$.
 Then we have a filling of the part of the boundary of the prism $\Delta^1\times\Delta^n$
 obtained by omitting one of the two foundations. This prism  has a standard triangulation into $n+1$ simplices
 of dimension $n+1$. We fill them inductively.  The first $n$ fillings are possible by the weak Kan property
 of $\X$, filling  horns of shapes $\Lambda^{n+1}_i$, $i=1,\cdots, n$. After this, we obtain an $(n-1)$-simplex
 in $\X_{/y_n}$ whose boundary is lifted into $\X_{/f_{y_n}}$, so the condition that $f_{y_n}$ is a Cartesian edge
 allows us to extend the lifting from the boundary to an entire simplex, thus obtaining an $(n-1)$-simplex
 in $\X_{/f_{y_n}}$. This $(n-1)$-simplex is interpreted as an $(n+1)$-simplex in $\X$ which completes
 the triangulation of the prism. The last face of this $(n+1)$-simplex is an $n$-simplex $\sigma$ 
 filling the other foundation of the prism.  
 We define  $G(\tau)=\sigma$, this defining $G$ on $n$-simplices. 
   \end{itemize}
   
   \noindent  Comparing with (the contravariant modification of) Example \ref{ex:groth-single}, 
   we note that the all the simplices constructed in the above inductive process, form  in fact a copy of the
   contravariant Grothendieck construction of $G$.  Therefore we have an $\oo$-functor $\chi(G)\to \X$
   compatible with the projection to $\Delta^1$. 
   One can further show that it is an equivalence of $\oo$-categories. 
 \end{ex}
 
\paragraph{Kan extensions.} Our frequent technical tool will be the theory of $\oo$-categorical Kan
extensions see \cite[\S 4.3] {lurie:htt} or a short summary in \cite[\S A.5]{dkss:1}. Since it will
be used only in the proofs, we do not give a summary here.  We will frequently be using the
following fundamental result about Kan extensions (see \cite[Prop. 4.3.2.15]{lurie:htt}):

\begin{prop}\label{prop:kan-ext-equiv}
(a) Let $\J $ be an $\oo$-category and $\I\subset\J$ a full subcategory.
Let $\D$ be a category with colimits and $\E\subset\Fun(\I,\D)$
be a full subcategory. Let $\E^!\subset \Fun(\J,\D)$
be the full subcategory formed by $\oo$-functors which are
left Kan extensions of functors from $\E$. Then the
restriction functor $\E^!\to\E$ is  a trivial Kan fibration, in particular, it is
an equivalence of $\oo$-categories.

(b) A similar statement for right Kan extensions, if $\D$ has limits. \qed
\end{prop}

In fact, we only need the existence of those (co)limits that enter in the standard pointwise
formulas for Kan extensions \cite[Eq. (A.5.2)]{dkss:1}. 
We will also be using a variant of Proposition \ref{prop:kan-ext-equiv} for relative Kan extension
which is also included in \cite[Prop. 4.3.2.15]{lurie:htt}.


\subsection {Adjunctions for $\oo$-categories}\label{subsec:adj-oo}

\paragraph{Grothendieck constructions and adjunctions.}
The Grothendieck constructions provide a very effective means to study
adjunctions of $\infty$-categories  \cite[\S 5.2.2]{lurie:htt}.

\begin{defi}\label{def:bicart}
Let $\C = \N C$ be the nerve of an ordinary category. An $\oo$-functor $p: \X\to\C$ is called
a {\em biCartesian fibration}, if it is both a Cartesian and coCartesian fibration. 
\end{defi}

A biCartesian fibration gives rise to two $\oo$-functors $\C\to\Cat_\oo$: a covariant
$\oo$-functor $Q^{\on{coCart}}$ classifying $p$ as a coCartesian fibration and a contravariant
$\oo$-functor $Q^{\on{Cart}}$ classifying $p$ as a Cartesian fibration. The values of $Q^{\on{coCart}}$ 
and  $Q^{\on{Cart}}$ on an object $c\in\C$ are the same:  the preimage $p^{-1}(c)\subset \X$, so for a morphism
$g: c\to d$ in $C$ we have a pair of  $\oo$-functors in the opposite directions
 \[
\xymatrix{
p^{-1}(c)  \ar@<.4ex>[rr]^{Q^{\on{coCart}}(g) } &&p^{-1}(d) \ar@<.4ex>[ll]^{Q^{\on{Cart}}(g) }. 
}
\]
By analogy with the classical category theory, it is natural to regard such $\oo$-functors as adjoint.  
By focusing on  pairs of single functors as in Example \ref{ex:groth-single}, this leads to the following definition.

\begin{defi} \begin{enumerate}
	\item An {\em adjunction} of $\infty$-categories is a biCartesian
	fibration 
	\[
		p: \X \lra \Delta^1.
	\]

	\item Let $F: \A \to \B$ be an $\oo$-functor of $\infty$-categories. We say
	that
	\begin{enumerate}
		\item $F$ admits a right adjoint if the covariant Grothendieck construction
			$\Gamma(F) \to \Delta^1$ (always a coCartesian fibration) is  also a Cartesian fibration.
		\item $F$ admits a left adjoint if the contravariant Grothendieck construction (always a Cartesian fibration)
			$\chi(F) \to \Delta^1$ is  also a coCartesian fibration.
	\end{enumerate}
	\end{enumerate}
\end{defi}

Given an adjunction $p: \X \to \Delta^1$, we define $\infty$-categories $\A = p^{-1}(0)$ and 
$\B= p^{-1}(1)$ as in  Example \ref{ex:groth-single}. Further, as in Example \ref{ex:implicit},
we construct  implicit $\oo$-functors
\be\label{eq:adj}
\xymatrix{
\A \ar@<.4ex>[r]^{F } &\B  \ar@<.4ex>[l]^{G }, 
}
\ee
defined by $p$ as a Cartesian ($G$) and coCartesian ($F$) fibration.

\begin{defi}
Let $F,G$ be $\oo$-functors in opposite directions as in \eqref{eq:adj}. We say that $G$ 
{\em is right adjoint to} $F$ and $F$ {\em is left adjoint to} $G$
(or equivalently,  $(F,G)$ {\em form an adjoint pair}), if there exists a biCartesian
fibration $p: \X\to\Delta^1$ such that $\A=p^{-1}(0)$, $\B=p^{-1}(1)$ and moreover, $F$ is an
implicit $\oo$-functor defined by $p$ as a coCartesian fibration while $G$ is an implicit
$\oo$-functor defined by $p$ as a Cartesian fibration. We will  use the notation $F\dashv G$ to
signify that $(F,G)$ form an adjoint pair, as well as to (implicitly) signify the adjunction itself. 
\end{defi}

Thus being an adjoint pair is a property and not an additional structure.


\paragraph{Units and counits.}  
A natural expectation for a workable concept of adjoint $\oo$-functors $F,G$ would be data of weak equivalences
 of
mapping spaces
\be
\Map_\B(F(a), b) \,\simeq \, \Map_\C(a,G(b))
\ee
or, the data of unit and counit morphisms of $\oo$-functors
\be
u: \Id_\A \lra GF, \quad c: FG\to\Id_\B
\ee
inducing such weak equivalences via functorialities \eqref{eq:map-functor}. 
This is achieved by the following result.

\begin{prop}\label{prop:oo-unit}
For a pair of $\oo$-functors
 $\xymatrix{
\A \ar@<.4ex>[r]^{F } &\B  \ar@<.4ex>[l]^{G } 
}$
 the following are equivalent:
\begin{itemize}
\item[(i)] $(F,G)$ is an adjoint pair.

\item [(ii)] There exists a morphism $u: \Id_\A\to GF$ (a {\em unit}) such that for any objects $a\in \A$, $b\in\B$ the composition
\[
\Map_\B(F(a), b) \buildrel G_*\over\lra \Map_\A(GF(a), G(b)) \buildrel u^*\over\lra \Map_\A(a, G(b))
\]
is a weak equivalence.

\item[(iii)] There exists a morphism $c: FG\to\Id_\B$ (a {\em counit}) such that for any objects
$a\in \A$, $b\in\B$ the
composition
\[
\Map_\A(a, G(b)) \buildrel F_* \over\lra \Map_\B(F(a), FG(B)) \buildrel c_*\over\lra \Map_\B(F(a), b)
\]
is a weak equivalence. 
\end{itemize}
\end{prop}

\noindent{\sl Proof:} The equivalence of (i) and (ii) is \cite[Prop. 5.2.2.8]{lurie:htt}. The equivalence of
(i) and (iii) is dual. \qed

\begin{cor}
If $(F,G)$ is an adjoint pair of $\oo$-functors as above, then  at the level of homotopy categories
we get an adjoint paiir of functors $\xymatrix{
\h\A \ar@<.4ex>[r]^{\h F } &\h\B  \ar@<.4ex>[l]^{\h G } 
}$ in the usual sense. \qed
\end{cor}

\vskip .2cm

We will now provide a description of the unit and counit, based on the universal properties of
Cartesian and coCartesian edges as effectively captured within Lurie's model structures on marked
simplicial sets. The diagramatic manipulations that we will use to construct and relate functors of
$\infty$-categories are very efficiently controlled by Lurie's theory of relative Kan extensions
\cite{lurie:htt}. Later in the text, we will omit details about analogous arguments.

\begin{lem}[(Units and counits implicitly)] \label{lem:unit}
Let $p: \X \to \Delta^1$ be an adjunction of  $\infty$-categories,
with $\A=p^{-1}(0)$, $\B=p^{-1}(1)$. 

\vskip .2cm

(a) Let 
	$
		\U \subset \Fun(\Delta^2, \X)
	$
	be the full subcategory spanned by $2$-simplices of the form
	\begin{equation}\label{eq:unitdiag}
		\begin{tikzcd}
			a \ar{rd}{!}\ar{d} & \\
			a' \ar{r}{*} & b
		\end{tikzcd}
	\end{equation}
	with $a,a'$ objects of $\A$, $b$ an object of $\B$, the edge $a \to b$ coCartesian, and the
	edge $a' \to b$ Cartesian. Then the functor 
	$
		\U \to\A
	$
	given by evaluation on the vertex $\{0\}$ (so a triangle  \eqref{eq:unitdiag}
	is sent into $a$), is a trivial Kan fibration (in particular, an equivalence of
	$\oo$-categories). 
	
	\vskip .2cm
	
	(b) Let $\C\subset \Fun(\Delta^2, \X)$ be the full subcategory spanned by $2$-simplices of the form 
		\begin{equation}\label{eq:counitdiag}
		\begin{tikzcd}
			b' \ar{rd} & \\
			a \ar{r}{*}  \ar{u}{!}& b
		\end{tikzcd}
	\end{equation}
	with $a$ an object of $\A$, $b,b'$ objects of $\B$, the edge $a\to b$
	Cartesian and the egde $a\to b'$ coCartesian. 
	Then the functor $\C\to\B$ given by evaluation on the vertex $\{2\}$
	(so a triangle \eqref{eq:counitdiag} is set into $b$) is
	a trivial Kan fibration  (in particular, an equivalence of
	$\oo$-categories). 

\end{lem}
\begin{proof} We prove (a), since (b) is similar.
First consider the $\infty$-category $\V \subset \Fun(\Delta^1, \X)$ spanned by the
	coCartesian edges $a \to b$ that cover the edge $0 \to 1$. By \cite[Ex. 4.3.1.4]{lurie:htt},
	the edges in $\V$ are precisely the left $p$-Kan extensions of their restriction to the
	vertex $a$ along the inclusion $\{0\} \subset \Delta^1$. By 
	\cite[Prop. 4.3.2.15]{lurie:htt} we therefore
	obtain that the functor $\V \to \A$ given by evaluation at $\{0\}$ is a trivial Kan
	fibration. Similarly, let $\V' \subset \Fun(\Delta^1, \X)$ denote the full subcategory
	spanned by the Cartesian edges $a' \to b$ in $\X$ covering $0 \to 1$. The projection functor
	$\V' \to \B$ given by evaluation at $\{1\}$ is a trivial Kan fibration, since the edges in
	$\V'$ are precisely the right $p$-Kan extensions of their restriction to $\{1\}$. Since
	trivial Kan fibrations are stable under pullback and composition, we obtain that the
	composite of
	\[
			\V \times_{\B} \V' \to \V \to \A
	\]
	is a trivial Kan fibration. Further, the horn inclusion of marked simplicial sets
	\[
		\Lambda^2_2 \subset \Delta^2
	\]
	with degenerate edges and the edge $\{1,2\}$ marked. This inclusion is marked anodyne
	\cite[Def. 3.1.1.1]{lurie:htt}. Further, the simplicial set $\X$, being a Cartesian
	fibration, defines (by marking Cartesian edges)
	 a
	fibrant object in the Cartesian model structure on 
	$(\Set_{\Delta}^+)_{/\Delta^1}$ (in $\Delta^1$ every edge is marked). 
	Therefore, the induced pullback map on mapping $\infty$-categories
	\[
		\Map^{\flat}_{\Delta^1}(\Delta^2, X) \lra \Map^{\flat}_{\Delta^1}(\Lambda^2_2, X),
	\]
	 is a trivial Kan
	fibration. Here the superscript $\flat$ means that we consider mapping $\oo$-categories with respect to  the above mentioned markings, see \cite[\S 3.1.3]{lurie:htt}. 
	We conclude the proof by noting that the map $\U \to \V \times_{\B} \V'$ is
	obtained as a pullback of the latter trivial Kan fibration so that it is a trivial Kan
	fibration as well. 
\end{proof}

\begin{rem}
Any morphism $a \to a'$ appearing as a part of a diagram \eqref{eq:unitdiag} is a potential choice of a
unit morphism $a \to G(F(a))$. Lemma \ref{lem:unit}(a) expresses the fact that, for every object $a$ of
$\A$, there is an essentially unique unit morphism $a \to a'$.  
Similarly, any morphism $b' \to b$ appearing as a part of a diagram
 \eqref{eq:counitdiag} is a potential choice of a
counit morphism $F(G(b))\to b$. Lemma \ref{lem:unit}(b) expresses the fact that, for every object $b$ of
$\B$, there is an essentially unique counit morphism $b'\to b$.  
\end{rem}


\paragraph{Exactness; adjoints of the composition.} 
For future use we record the analogs of some familiar properties of adjoints known in the classical
category theory. 

\begin{prop}\label{prop:adj-exact} \cite[Prop. 5.2.4.3]{lurie:htt}
If an $\oo$-functor $F:\A\to\B$  has a right adjoint, then it is left exact (i.e., preserves all colimits which exist in $\A$).
If $F$ has a left adjoint, then it is right exact (i.e., preserves all limits that exist in $\A$). 
\qed
\end{prop}

\begin{prop}\cite[Prop. 5.2.2.6]{lurie:htt}
Let $\xymatrix{
\A \ar@<.4ex>[r]^{F } &\B  \ar@<.4ex>[l]^{G } 
 \ar@<.4ex>[r]^{F' } &\C  \ar@<.4ex>[l]^{G'}
}$ 
be a diagram of $\oo$-categories and 
$\oo$-functors such that $(F,G)$ and $(F', G')$ are adjoint pairs.
Then $(F'\circ F, G\circ G')$ is also an adjoint pair. \qed
\end{prop}


\subsection{Comparison with other approaches to adjunctions}\label{subsec:comparison}

\paragraph{Fourier-Mukai kernels.}
As mentioned in \S\ref{subset:stab-triang}, any dg-category $\Ac$ gives an $\oo$-category
$\N^\dg\Ac$; any dg-functor $f: \Ac\to\Bc$ gives an $\oo$-functor
$\N^\dg f: \N^\dg \Ac \to \N^\dg \Bc$. Therefore the $\oo$-categorical theory of adjoints
provides an approach to defining adjointness for dg-functors. In this and the next
paragraph we compare it with  other approaches to this question that have been
used in the literature.  

\vskip .2cm

Let $X$ be a smooth projective algebraic variety over a field $\k$.
We denote by $\Coh_X$ the abelian category of coherent sheaves on $X$,
and by $D^b\Coh_X$ its bounded derived category. We can think of 
objects of $D^b\Coh_X$ as complexes of quasicoherent sheaves with
bounded coherent cohomology. 
We have eqiivalences of
(usual) categories
\[
D^b\Coh_X \,\simeq\, H^0 D^b_\dg \Coh_X\,\simeq\, \h D^b_{(\oo)} \Coh_X,
\]
where: 
\begin{itemize}
\item $D^b_\dg \Coh_X$ is the dg-enhancement of $D^b\Coh_X $ obtained
by taking an injective  resolution of each  object and
forming the Hom-complexes between such resolutions \cite{BK:enhanced}.

\item $D^b_{(\oo)} \Coh_X = \N^\dg D^b_\dg \Coh_X$ is the dg-nerve
of the dg-enhancement. 
\end{itemize}

Let  $X,Y$ be two smooth projective varieties.
Denote $X\buildrel p_X\over\leftarrow X\times Y \buildrel p_Y\over\to Y$
the projections. For an object    $\Kc\in D^b\Coh_{X\times Y}$
 we have the {\em Fourier-Mukai functor}
\[
F_\Kc   : D^b\Coh_X \lra D^b\Coh_Y, \quad \Fc\mapsto Rp_{Y*} (p_X^*\Fc \otimes \Kc). 
\]
The complex $\Kc$ is referred to as the {\em Fourier-Mukai kernel} (of $F_\Kc$), see 
\cite{huybrechts}. 
By passing to appropriate resolutions, $F_\Kc$ can be lifted to a dg-functor and
then to 
to an $\oo$-functor
\[
F_\Kc^\dg:  D^b_\dg \Coh_X\lra  D^b_\dg \Coh_Y,\quad 
F_\Kc^{(\oo)} =\N^\dg F_\Kc^\dg: D^b_{(\oo)} \Coh_X\lra D^b_{(\oo)} \Coh_Y
\]
respectively. 

If $W$ is another smooth projective variety and $\Lc\in D^b\Coh_{W\times X}$ is
another Fourier-Mukai kernel, then the composition
$F_\Kc\circ F_\Lc: D^b\Coh_W\to D^b\Coh_X$ is identified with $F_{\Kc*\Lc}$,
where the complex $\Kc*\Lc$, known as the {\em convolution} of $\Kc$ and $\Lc$,  is defined as
\[
\Kc*\Lc\,=\, Rp_{WX*} (p_{XY}^*\Kc \otimes p_{WY}^*\Lc),
\]
with $p_{WX}, p_{XY}, p_{WY}$ are the projections of $W\times X\times Y$ to the pairwise
products. The classical Serre-Grothendieck duality implies 
\cite[Lem. 1.2]{bondal-orlov}
that a left adjoint 
of $F_\Kc$ in the sense of usual categories is provided by $F_{^\bigstar\hskip -0.05cm  \Kc}$,
where 
\[
^\bigstar \hskip -0.05cm \Kc \,=\, (\sigma^*\Kc)\otimes p_X^* \Omega^n_X[n] \,\in\, D^b\Coh_{Y\times X},
\quad n=\dim(X), \,\, \sigma:Y\times X\to X\times Y,
\]
with $\sigma$ being the permutation. Further
\cite[Th.  3.1]{AL:Fourier} this theory provides a natural trace morphism
\[
{^\bigstar\hskip -0.05cm  \Kc}*\Kc \lra \Oc_{\Delta_X},
\] 
which induces the counit
for the adjoint pair $(F_{{^\bigstar\hskip -0.05cm  \Kc}}, F_\Kc)$. 
Proposition \ref{prop:oo-unit} implies therefore the following
corollary which can   be seen as an instance of \cite[Prop. 5.2.2.12]{lurie:htt}.

\begin{cor}
 $(F_{{^\bigstar\hskip -0.05cm  \Kc}}^{(\oo)}, F_\Kc^{(\oo)})$
 form an adjoint pair of $\oo$-functors. \qed
\end{cor}

Therefore  the   theory of dg-enhanced adjoints  for Fourier-Mukai functors
as developed in  \cite{AL:Fourier}
and used in \cite{AL:spherical},  can be embedded into the
$\oo$-categorical approach adopted in this paper. 


\paragraph{ Bimodule dg-functors.} In a similar spirit,
let $\Ac$ be a small dg-category, and $\Mod_\Ac$ be the dg-category of right 
$\Ac$-modules, i.e., of contravariant dg-functors $\Ac\to\dgVect_\k$. 
The {\em derived category} of $\Ac$, denoted $D(\Ac)$,  is the localization of
$H^0 \Mod_\Ac$ (the homotopy category) by quasi-isomorphisms. 
It is a triangulated category which can be represented as
\[
D(\Ac) \, \simeq \, H^0 D_\dg(\Ac) \,\=\, \h D_{(\oo)}(\Ac), 
\]
 where:
 \begin{itemize}
 \item $D_\dg(\Ac)$ is the dg-enhancement of $D(\Ac)$ constructed in
 \cite{AL:bar}  and denoted $\ol{\bf Mod}\text{-} \Ac$. It is defined
 using the standard  bar-resolution $\on{Bar}(E)$ for
 any $E\in\Mod_\A$, that is
 \[
 \Hom^\bullet_{D_\dg(\Ac)}(E,F) \, := \, \Hom^\bullet_{\Mod_\Ac)}(\on{Bar}(E), F).
 \]
 Thus closed degree $0$ morphisms in $D_\dg(\Ac)$ are $A_\oo$-morphisms
 of dg-modules. 
 
 \item $ D_{(\oo)}(\Ac) = N^\dg D_\dg(\Ac)$ is the $\oo$-category obtained
 as the dg-nerve of the dg-enhancement. 
 \end{itemize}
 
 Let $\Ac, \Bc$ be two dg-categories. An $(\Ac, \Bc)$-{\em bimodule}  
 is, by definition, a covariant dg-functor $\Ac\times \Bc^\op \to\dgVect_\k$. Any  
 $(\Ac, \Bc)$-bimodule $M$
 defines a dg-functor $F_M^\dg:  D_\dg(\Ac)\to D_\dg(\Bc)$ by the derived
 tensor product over $\Ac$,  see \cite[Def. 3.9]{AL:bar}. It induces an
 $\oo$-functor $F_M^{(\oo)}: D_{(\oo)}(\Ac)\to D_{(\oo)}(\Bc)$ by passing
 to dg--nerves and an exact functor of triangulated categories
 $F_M: D(\Ac)\to D(\Bc)$ by passing to $H^0$-categories. 
We define the {\em left and right transposes} of $M$ to be the 
$(\Bc, \Ac)$-bimodules whose values on the objects are
given by
\[
\begin{gathered}
^\bigstar \hskip -0.05cm M(b,a)  \,=\,\Hom^\bullet_{D_{\dg}(\Ac^\op)}
\bigl( M(-,b),\, \Hom^\bullet_\Ac(a,-)\bigr),\\
 M^\bigstar(b,a)\,=\, \Hom^\bullet_{D_{\dg}(\Bc)}
 \bigl(M(a,-),\, \Hom^\bullet_\Bc(-,b)\bigr), 
\end{gathered} 
\]
cf.  \cite[Def. 3.31]{AL:bar}. Theorem 1.1 of \cite{AL:bar} together with 
Proposition \ref{prop:oo-unit} above, give the following fact.
\begin{prop}
(a) If $M$ is $\Bc$-perfect, then $F^{(\oo)}_{M^\bigstar}$ is a right
$\oo$-categorical adjoint of $F^{(\oo)}_M$, and $F_{M^\bigstar}$
is a right adjoint of $F_M$ in the sense of usual categories.

\vskip .2cm

(b) If $M$ is $\Ac$-perfect, then $F^{(\oo)}_{^\bigstar \hskip -0.05cm M}$
is a left $\oo$-categorical adjoint of $F^{(\oo)}_M$, and
$F_{^\bigstar \hskip -0.05cm M}$ is a left adjoint of $F_M$ in the sense of
usual categories. \qed
\end{prop}

Therefore  the   theory of dg-enhanced adjoints  for bimodule functors
 developed in  \cite{AL:bar}
  can be embedded into the
$\oo$-categorical approach adopted in this paper. 


\paragraph{Homotopy coherent adjunctions. } 
The concept of adjunction can give rise to a certain ambiguity
which is present  already in the classical
case of functors 
 $\xymatrix{
A \ar@<.4ex>[r]^{F } &B  \ar@<.4ex>[l]^{G } 
}$
 between usual categories. In this case the strongest and most natural adjunction
 data for $(F,G)$ (a {\em coherent adjunction}) consists 
  \cite{maclane} of two transformations $u: \Id_A\to GF$ (unit) and $c: FG\to\Id_B$ (counit)
satisfying the following compatibility condition:
\begin{itemize}
\item[(AC)] The compositions
\[
F=F\circ \Id_A \buildrel F\circ u \over\lra F\circ G\circ F \buildrel c\circ F\over\lra \Id_B\circ F =F,
\quad
G=\Id_A\circ G \buildrel u\circ G\over\lra G\circ F\circ G \buildrel G\circ c\over\lra G\circ
\Id_B=G
\]
are equal to the identity transformations of $F$ and $G$ respectively. 
\end{itemize}
This condition implies that for any objects $a\in A$, $b\in B$ the maps
\[
\xymatrix{
\Hom_B(F(a),b)  \ar@<.4ex>[r]^{u^*_{a,b} } &\Hom_A(a,G(b)) 
 \ar@<.4ex>[l]^{c_*^{a,b} } 
}
\]
 defined as the compositions (cf. 
   Proposition \ref{prop:oo-unit})
 \[
 \begin{gathered}
\Hom_B(F(a), b) \buildrel G_*\over\lra \Hom_A (GF(a), G(b)) \buildrel u^*\over\lra 
\Hom_A(a, G(b)), 
\\
 \Hom_A(a, G(b)) \buildrel F_* \over\lra \Hom_B(F(a), FG(b)) \buildrel c_*\over\lra \Hom_B(F(a), b),
 \end{gathered}
\]
 are inverse to each other. 
 
 Now, a seemingly weaker datum would be that of a
 single transformation, say $u$,
   such that all the $u^*_{a,b}$ are bijections (or $c$ such that all the
   $c_*^{a,b}$ are bijections). However, it is a classical
   fact that any such partial datum can be always completed to
   a  pair $(u,c)$ satisfying (AC). So the possible ambiguity in the definition
   of an adjunction is, in the case of usual categories, an inessential one. 
   
   \vskip .2cm
   
  In the $\oo$-categorical case we have a similar ambiguity. To address it,  one needs to
   do two things:
   \begin{itemize}
   \item[(1)]  Formulate an $\oo$-categorical analog of the condition (AC) 
   (i.e., define the stronger concept of a {\em homotopy coherent adjunction}). 
   
   \item[(2)]   Compare this concept with the approach via biCartesian fibrations
   over $\Delta^1$ discussed above. 
   \end{itemize} 
   
  \noindent  This has been done in \cite{riehl-verity}.  For completeness, let us
  give some details of their analysis. 
  
  \vskip .2cm
  
\begin{figure}[h]
  \centering
  \begin{tikzpicture}[scale=0.4]
 
 \node at (0,0){$\bullet$}; 
  \node at (-1,0){$a$}; 
  
  \node at (12,0){$\bullet$}; 
 \node at (13,0){$b$}; 
  
  \draw  [dashed, decoration={markings,mark=at position 0.6 with
{\arrow[scale=1.5,>=stealth]{>}}},postaction={decorate},
line width=.2mm]  (0,0) arc (180:0:6cm and 2cm); 

  \draw  [ decoration={markings,mark=at position 0.6 with
{\arrow[scale=1.5,>=stealth]{>}}},postaction={decorate},
line width=.2mm]  (12,0) arc (360:180:6cm and 2cm); 

\draw (6,0) circle (6cm); 

\node at (4.5,-2.7){$g$}; 
\node at (7.6, 2.8){$f$}; 

\node[rotate=45] at (1.4,1.7) {\huge $\Rightarrow$}; 
\node[rotate=45] at (1,2.3) {$u$}; 

\node[rotate=45] at (11, -2.2) {\huge $\Rightarrow$}; 
\node[rotate=45] at (10.4,-1.8) {$c$}; 
  
  \end{tikzpicture}
  \caption{Generating morphisms of $\Adj^{\leq 2}$ and  cells of the $2$-sphere.}
  \label{fig:2-sphere}
  \end{figure}
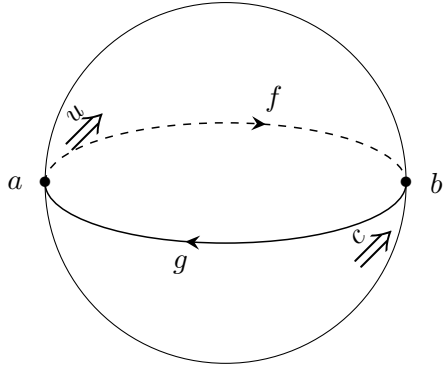

  Following \cite{auderset, schanuel-street}, let $\Adj^{\leq 2}$ be the {\em strict (classical)
  $2$-category generated by the universal adjunction}. That is, $\Adj^{\leq 2}$ has
  two objects $a$ and $b$ corresponding to indeterminate categories $A$ and $B$ above.
  As a $2$-category, it is generated by
  two $1$-morphisms $f: a\to b$ and $g: b\to a$ corresponding to indeterminate
  functors $F,G$ above and two $2$-morphisms $u: \Id_a\to gf$ and $c: fg\to\Id_b$. 
  which are subject only to the two relations mimicking  (AC). Geometrically, the
  objects, $1$- and $2$-morphisms  of $\Adj^{\leq 2}$ can be depicted as $0$- $1$-
  and $2$-dimensional cells of the standard ``equatorial'' cell
  decomposition of the $2$-sphere $S^2$, and the two relations in (AC) as two
  $3$-cells which extend $S^2$ to a $3$-sphere $S^3$ for which $S^2$ is
  the equator, see Fig. \ref{fig:2-sphere}. 
  
  By construction, a (coherent) adjunction of usual categories is the same
  as a strict $2$-functor  of strict $2$-categories $\Adj^{\leq 2} \to\Cat$. 
  
  \begin{ex}\label{ex:univ-monad}
  The above definition of $\Adj^{\leq 2}$ is implicit, but one can describe it completely
  explicitly   \cite{auderset, schanuel-street}. For example, 
  the category
  $\Hom_{\Adj^{\leq 2}}(a,a)$ is isomorphic to $\Delta_+$, the augmented simplicial
  category formed by all finite ordinals, including the empty one, so that the
  ordinal $\{1,\cdots, n\}$ (taken to be $\emptyset$ for $n=0$)
  corresponds to $(gf)^n$. Similarly, 
  $\Hom_{\Adj^{\leq 2}}(b,b)$ is isomorphic to $\Delta^\op_+$, the opposite category,
  with  $\{1,\cdots, n\}$  corresponding to $(fg)^n$. This reflects the fact that
  for an adjoint pair $(F,G)$ the composition $GF$ is a comonad, while $FG$
  is a monad. 
  \end{ex}
  
  As with any strict $2$-category, $\Adj^{\leq 2}$ gives rise to an $(\oo,2)$-category
  $\N^\Hom \Adj^{\leq 2}$ by taking the nerves of all the $\Hom$-categories,
  see Definition \ref{def:oo-2}.
  
  \begin{defi}\cite{riehl-verity}
  A {\em homotopy coherent adjunction of $\oo$-categories} is a functor of
  simplicial categories $ N^\Hom \Adj^{\leq 2} \to\Cat_\oo^{(2)}$, 
   where $\Cat_\oo^{(2)}$ is the $(\oo,2)$-category of all $\oo$-categories,
   see \S \ref{subset:stab-triang}\ref{par:oo-2-level}
  \end{defi}
  
  \begin{thm}\cite{riehl-verity}
   Let  $\xymatrix{
\A \ar@<.4ex>[r]^{F } &\B  \ar@<.4ex>[l]^{G } 
}$ be a pair of $\oo$-functors between $\oo$-categories. The following are equivalent:
\begin{itemize}
\item[(i)] $(F,G)$ is an adjoint pair.

\item[(ii)] There exists a homotopy coherent adjunction $N^\Hom \Adj^{\leq 2} \to\Cat_\oo^{(2)}$ sending $a,b$ to $\A,\B$ and $f,g$ to $F,G$ respectively
\end{itemize}
\end{thm}

 \noindent{\sl Proof:} (i)$\Rightarrow$(ii) is the content of \cite[Th. 4.3.9]{riehl-verity} taken 
 together with  equivalence between the two definition of adjoints for $\oo$-functors:
 that of \cite{lurie:htt} which we use and that of \cite{riehl-verity}. That equivalence
 is established in \cite {riehl-verity:2-cat, riehl-verity:yoneda}, see the introduction
 to \cite     {riehl-verity}. The implication (ii)$\Rightarrow$(i) follows from Proposition
  \ref{prop:oo-unit}. 
 
 \qed
 
 \vskip .2cm
 
 Therefore the possible ambiguity as to the notion of adjointness for $\oo$-categories
 is also inessential.  
 \begin{rem}\label{rem:adj-sphere}
 The geometric datum (Fig. \ref{fig:2-sphere}) used to define $\Adj^{\leq 2}$
 has an analog for any $n$. This is the sphere $S^{n+1}$ with its equatorial
 cell decomposition (two $p$-cells for each $0\leq p\leq n+1$). These cells
 are oriented so that for any $p$ the two $p$-cells comprising $S^p\subset S^{n+1}$
 induce the same orientation of $S^p$. This means that their ``directions'' 
 are, intuitively, the opposite of each other, unlike the standard globe
 diagram  \cite{street:parity}, where these directions are the same.
 
 It would be interesting to construct an $n$-category (in some sense)
 $\Adj^{\leq n}$ ``freely generated'' by the above $S^{n+1}$ so that the $p$-cells
 for $p\leq n$ give generating $p$-morphisms, and the two $(n+1)$-cells
 serve as relations.  Intuitively, $\Adj^{\leq n}$ should describe adjunctions
 for $(n-1)$-categories. 
 
 It is not clear whether one can construct $\Adj^{\leq n}$
 as a classical strict $n$-category in the spirit of Street's orientals
 \cite{street:orientals}, as this $S^{n+1}$ is not a loop-free pasting diagram
 \cite{street:parity}. 
 
 A study of adjunctions for $2$-categories undertaken in
 \cite{macdonald-stone}, leads to a construction of $\Adj^{\leq 3}$
 as a Gray $3$-category (i.e., with $2$-dimensional associativity holding
 up to a canonical ``braiding'' $3$-isomorphism). From the point of view of
 Example \ref{ex:univ-monad}, the $2$-category $\Hom(a,a)$ in that $\Adj^{\leq 3}$
 (the universal $2$-comonad) is a certain $2$-categorification of the category
 $\Delta_+$ in which faces and degenerations become adjoint
 $1$-morphisms \cite{macdonald-scull}.
  A strict $2$-categorification of $\Delta_+$ with these
 properties can be obtained \cite[\S 3.1.3]{dyckerhoff:DK} by viewing
 each finite ordinal (a poset) as a category. It will be used below
 in \S \ref {subsec:sph-admis}\ref{par:admis}
 \end{rem}


 \subsection{Spherical adjunctions of stable $\oo$-categories}\label{subsec:sphad}

\paragraph{Adjunctions of stable $\oo$-categories.} 
By an {\em adjunction of stable $\infty$-categories} 
we will mean an adjunction of $\oo$-categories, i.e., a biCartesian fibration 
$p: X \to \Delta^1$,
such that both $\oo$-categories  $\A= p^{-1}(0)$ and $\B= p^{-1}(1)$ are stable.

\begin{prop}
Let $\xymatrix{
\A \ar@<.4ex>[r]^{F } &\B  \ar@<.4ex>[l]^{G }, 
}$ be an adjoint pair of $\oo$-functors with $\A$ and $\B$ being stable $\oo$-categories.
Then both $F$ and $G$ are exact functors. 
\end{prop}

\noindent{\sl Proof:} By Proposition \ref{prop:adj-exact}, $F$ is left exact and $G$ is
right exact. But between stable $\oo$-categories, any left or right exact functor
is automatically exact, see  \S\ref{subset:stab-triang}\ref{par:stable-oo} \qed


\paragraph{Twist and cotwist.} 
Let $p: \X\to\Delta^1$ be an adjunction of stable $\oo$-categories and 
$\xymatrix{
\A \ar@<.4ex>[r]^{F } &\B  \ar@<.4ex>[l]^{G } 
}$
be an adjoint pair of $\oo$-functors associated to $p$. By Proposition \ref {prop:oo-unit}
we can choose unit and counit transformations
\[
u: \Id_\A \lra GF, \quad c: FG\lra \Id_\B. 
\]
By analogy with the pre-triangulated dg-setting \cite{AL:spherical} 
we  define  the {\em twist}  and {\em cotwist functors} associated
to the adjunction are  as the cofiber of the unit and  the fiber of the counit,  see \eqref{eq:(co)fiber}:
\be\label{eq:twist-as-cofiber}
T_\A\,=\,\Cof(u): \A \lra \A, \quad T_\B \,=\,\Fib(c): \B\lra\B. 
\ee

We will now provide a more implicit description of twist and cotwist based on universal properties
which will be used extensively lateron.

\begin{lem}\label{lem:twist} Let $p:X \to \Delta^1$ be an adjunction of stable $\infty$-categories. Let 
	\[
		\M \subset \Fun(\Delta^1 \times \Delta^1 \amalg_{\Delta^1} \Delta^2, X)
	\]
	denote the full subcategory spanned by diagrams of the form
	\[
		\begin{tikzcd}
			0 \ar{d} & \ar{l} a \ar{rd}{!}\ar{d} & \\
			a'' & \ar{l} a' \ar{r}{*} & b
		\end{tikzcd}
	\]
	where the right-hand triangle forms a diagram as in \eqref{eq:unitdiag} and the left-hand
	square is a biCartesian square in $\A$. Then the functor 
	\[
		\M \lra \A
	\]
	given by restriction to $a$ is a trivial Kan fibration.
\end{lem}
\begin{proof} The diagrams in $\A$ of the form 
	\begin{equation}\label{eq:rightkan}
			\begin{tikzcd} 
				a \ar{d}\ar{r} & a'\\
				0
			\end{tikzcd}
	\end{equation}
	are precisely the right Kan extension of their restriction to $a \to a'$. Further, the
	biCartesian squares in $\A$ of the form
	\begin{equation}\label{eq:bicart}
			\begin{tikzcd} 
				a \ar{d}\ar{r} & a' \ar{d}\\
				0 \ar{r} & a''
			\end{tikzcd}
	\end{equation}
	are precisely the left Kan extensions of their restriction to \eqref{eq:rightkan}. Let $\Nc
	\subset \Fun(\Delta^1 \times \Delta^1, \A)$ denote the full subcategory spanned by the
	biCartesian squares of the form \eqref{eq:bicart}. Applying \cite[4.3.2.15]{lurie:htt}
	twice, we obtain that the restriction functor
	\[
		\Nc \lra \Fun(\Delta^1, \A)
	\]
	given by restriction to the edge $a \to a'$ is a trivial Kan fibration. We have a pullback
	square
	\[
		\begin{tikzcd}
			\M \ar{r}\ar{d} & \Nc \ar{d} \\
			\V \ar{r} & \Fun(\Delta^1,\A)
		\end{tikzcd}
	\]
	so that the map $\M \to \V$ is a trivial Kan fibration. We conclude the proof by Lemma
	\ref{lem:unit}.
\end{proof}

Given an adjunction $p: \X \to \Delta^1$ of stable $\infty$-categories, we consider the
correspondence
\[
	\A \overset{r}{\lla} \M \overset{s}{\lra} \A
\]
where $r$ is the trivial Kan fibration from Lemma \ref{lem:twist} and $s$ is the projection
to $a''$. We obtain a functor
\begin{equation}\label{eq:twist}
		T_{\A}: \A \lra \A
\end{equation}
as the composite of a section of $r$ with $s$. Comparing with \eqref{eq:twist-as-cofiber},
we see that 
 $T_{\A}$ is  the  twist associated with the
adjunction (which justifies the same notation).
 Indeed, it  assigns to an object $a$ in $\A$ a cofiber of a unit
morphism $a \to G(F(a))$ as constructed in Lemma \ref{lem:unit}.

Applying the same construction to $p^{\op}$, we obtain a functor
\begin{equation}\label{eq:cotwist}
		T_{\B}: \B \lra \B,
\end{equation}
which is identified with  the  cotwist  \eqref{eq:twist-as-cofiber},
as it associates to an object $b$ in $\B$ a fiber of a counit morphism 
$F(G(b)) \to b$.


\paragraph{Spherical adjunctions.} 

\begin{defi}\label{def:sphad}
 An adjunction of stable $\infty$-categories $p: X \to \Delta^1$ is called {\em
	spherical} if its associated twist and cotwist are equivalences of $\infty$-categories. An adjoint pair  $\xymatrix{
\A \ar@<.4ex>[r]^{F } &\B  \ar@<.4ex>[l]^{G }
}$
associated to a spherical adjunction, will be called an {\em adjoint pair of
spherical $\oo$-functors}. An $\oo$-functor will be called {\em left spherical},
resp {\em right spherical}, if it is the respective part of an adjoint
pair $(F,G)$ of spherical $\oo$-functors. 
\end{defi}

We will show later (Corollary \ref{cor:left-sph-right}) that the concepts
of a left and a right spherical functors are in fact equivalent, so the term
``spherical functor'' can be used unambiguously. 
Further, the usual definition of 
 of a spherical functor, cf.  \cite[Th. 1.1]{AL:spherical} includes
 a requirement of the existence of one more adjoint, for example,
 the right adjoint of $G$. We will show later (Corollary \ref{cor:sph-oo-adjoints})
 that this requirement is also unnecessary. 
 

\paragraph{Examples of spherical adjunctions.}

\begin{ex}[(``Classical'' dg-examples)]
Because of the comparison done in \S  \ref{subsec:comparison}, all the
examples of spherical dg-functors between pre-triangulated dg-categories
that have been considered in the literature, translate into our framework.
Let us consider just two examples. 

\vskip .2cm

(a) (Calabi-Yau fibrations) Let $X\buildrel p\over\to B$ be a smooth projective
morphism of algebric varieties over $\k$. Assume that all the fibers $X_b=p^{-1}(b)$,
$b\in B$
are $n$-dimensional Calabi-Yau varieties, i.e., $H^i(X_b, \Oc)$ is
isomorphic to $\k_b$ (the residue field of the point $b$) for $i=0,n$ and is zero otherwise.
Then 
\[
\xymatrix{
D^b_{(\oo)}\Coh_B  \ar@<.4ex>[r]^{p^* } &D^b_{(\oo)}\Coh_X  \ar@<.4ex>[l]^{Rp_* }, 
}
\]
form an adjoint pair of spherical $\oo$-functors. 
\vskip .2cm

(b) (Sphere fibrations) For a CW-complex $Y$ we denote by $\Sh_Y$ the
category of sheaves of $\k$-vector spaces on $Y$ and by $D^b_{(\oo)}\Sh_Y$
the corresponding derived $\oo$-category. Let $S\buildrel p\over\to B$
be a locally trivial fibration of CW-complexes with fibers $S_b=p^{-1}(b)$
homeomorphic to the $n$-sphere $S^n$.  Then 
\[
\xymatrix{
D^b_{(\oo)}\Sh_B  \ar@<.4ex>[r]^{p^* } &D^b_{(\oo)}\Sh_S  \ar@<.4ex>[l]^{Rp_* }, 
}
\]
form an adjoint pair of spherical $\oo$-functors. 

One can obtain a smaller example by restricting to full $\oo$-subcategories consisting of complexes
with locally constant cohomology sheaves, i.e., to local systems with values in the stable
$\oo$-category $D^b_{(\oo)}\Vect_\k$. This class of examples was generalized to local systems with
values in an arbitrary stable $\oo$-category in \cite{christ}. 

\end{ex}

\begin{ex}[(Versions of the sphere example)]\label{ex:versions-sphere}

(a) (Version  of the $S^0$-example) Let $\A$ be any stable $\oo$-category.
Let $\delta: \A\to \A\times A$ be the diagonal embedding. Then
\[
\xymatrix{
\A  \ar@<.4ex>[r]^{\hskip -0.5cm\delta } & \A \times\A   \ar@<.4ex>[l]^{\hskip -0.5cm \oplus }, 
}
\]
 form an adjoint pair of spherical $\oo$-functors. The twist $T_\A$ is the identity,
 while $T_{\A\times\A}$ is the permutation functor. 

\vskip .2cm

(b) (Version of    the $S^1$-example) A local system of $\k$-vector spaces on $S^1$
is the same as a $\k[t,t^{-1}]$-module, and pullback of (trivial) local systems 
under $S^1\to\pt$ corresponds to the pullback, denote it $\eps$, 
 of $\k$-modules under the
evaluation map $\k[t, t^{-1}]\to\k$, $t\mapsto 1$.   
At the purely algebraic level, one can see directly that
\[
\xymatrix{
D^b_{(\oo)} \Vect_\k   \ar@<.4ex>[r]^{ \hskip -0.5cm \eps } &D^b_{(\oo)}\Mod_{\k[t,t^{-1}]}  \ar@<.4ex>[l]^{\hskip -0.5cm \eta }, 
}
\quad \eta (M) =  R\Hom_{\k[t,t^{-1}]}(\k, M), 
\]
 form an adjoint pair of spherical $\oo$-functors. The functor $\eta$
 can be calculated using the $2$-term  free resolution   $\{\k[t,t^{-1}]\buildrel t-1\over\to 
 \k[t,t^{-1}]\}$ of $\k$.
 This gives 
the twist on $D^b_{(\oo)}\Vect_k$ to be   identified with the shift 
functor $V\mapsto V[-1]$, while the cotwist on $D^b_{(\oo)}\Mod_{\k[t,t^{-1}]}$ is  identified with the identity.

We also have an essentially equivalent
 version with the left adjoint $\delta$ of $\eps$ (which corresponds to taking
 homology instead of cohomology of local systems on $S^1$):
\[
\xymatrix{  D^b_{(\oo)}\Mod_{\k[t,t^{-1}]} 
  \ar@<.4ex>[r]^{ \hskip 0.5cm \delta } & D^b_{(\oo)} \Vect_\k 
 \ar@<.4ex>[l]^{\hskip 0.5cm \eps }, 
}
\quad \delta (M) =   M\otimes^L_{\k[t,t^{-1}]}\k. 
\]
Here the twist on $ D^b_{(\oo)}\Mod_{\k[t,t^{-1}]} $ is the identity while
the cotwist on $D^b_{(\oo)} \Vect_\k$ is the functor $V\mapsto V[1]$. 

\vskip .2cm

(c) (Global algebro-geometric version of (b)) Let $X$ be a smooth algebraic
variety over $\k$ and $i: Y\hra X$ is an irreducible hypersurface, so that we have
the line bundle (sheaf of ideals)
 $\Oc_X(-Y)
\subset \Oc_X$  with $\Oc_X(-Y)|_Y = N^*_{Y/X}$, the conormal bundle. 
Then
\[
\xymatrix{
 D^b_{(\oo)}\Coh_X  \ar@<.4ex>[r]^{\delta } &  D^b_{(\oo)}\Coh_Y
\ar@<.4ex>[l]^{\eps=i_*  }
},
\quad \delta = (-)\otimes^L_{\Oc_X} \Oc_Y. 
\]
 form an adjoint pair of spherical $\oo$-functors. 
 The twist on $D^b_{(\oo)}\Coh_X$ is given by
 tensoring with $\Oc_X(-Y)$, 
 while the cotwist on
 $D^b_{(\oo)}\Coh_Y$ is given by tensoring with
 $N^*_{Y/X}[1]$,  see \cite{segal}. 

\vskip .2cm

(d) (Graded $\k[t]$-modules) We will need yet another version of the above.
Let $\Vect_\k^\ZZ$ be the category of $\ZZ$-graded $\k$-vector spaces and
$\Mod_{\k[t]}^\ZZ$ be the category of $\ZZ$-graded $\k[t]$-modules, with $\deg(t)=1$.
Thus an object of $\Mod_{\k[t]}^\ZZ$ is a diagram of vector spaces and linear operators
\be\label{eq:graded-k[t]}
M\,=\, \bigl\{ \cdots \buildrel t\over\lra M^{-1}  \buildrel t\over\lra M^0  \buildrel t\over\lra M^1
 \buildrel t\over\lra\cdots\bigr\}
\ee
The embedding  $\eps: \Vect_\k^\ZZ\hra \Mod_{\k[t]}^\ZZ$ is the pullback
under the evaluation map $\k[t]\to \k$ at $t=0$.  As in (b),
\[
\xymatrix{
 \A=D^b_{(\oo)}\Mod_{\k[t]}^\ZZ 
 \ar@<.4ex>[r]^{   \delta } &  D^b_{(\oo)} \Vect^\ZZ _\k =\B  
\ar@<.4ex>[l]^{  \eps }, 
}
\quad \delta  = (-) \otimes^L_{\k[t]}\k, 
\]
  form an adjoint pair of spherical $\oo$-functors. 
 Note that  a  graded $\k[t]$-module $M$ is $\delta$-flat, if $t$ is injective, i.e.,
 \eqref{eq:graded-k[t]}  represents a filtration. In this case $\delta (M)= \gr (M)$
 is the associated graded space of this filtration. 
 
 As before, calculating $\delta$ using the $2$-term graded
 resolution $\{\k[t]\buildrel t\over\to \k[t]\}$ of $\k$ gives that
 the twist $T_\A$   is given by
 the shift of grading of the graded $\k[t]$-modules $(T_\A M)_i = M_{i-1}$,
 The cotwist $T_\B$
 is given by $(T_\B(V))_i=V_{i-1} [1]$.    
 \end{ex}


\subsection{Spherical adjunctions from Massey products in stable $\infty$-categories}
\label{sec:massey}

In this section we work out some particular examples of spherical adjunctions which
will be used later in the study of the relative Waldhausen $S$-construction.
Recall \cite{waldhausen:free} that the main idea behind the ``classical''
$S$-construction (for exact categories in the sense of Quillen) is to
rigidify the concept of a  filtered object by including
all the quotients of the filtration into the set of data. Adapting this to
the context of stable $\oo$-categories and to filtrations of infinite length,
we arrive at the following concepts. 

\paragraph{Filtered and graded objects.} Let $\D$ be a stable $\oo$-category. 

\begin{defi}
(a) A {\em filtered object} in $\D$ is an $\oo$-functor $Z: (\ZZ,\leq)\to\D$, i.e.,
a diagram of objects and morphsims in $\D$
\[
Z\,=\,\bigl\{ \cdots\to  Z_i \lra Z_{i+1} \lra Z_{i+2} \to \cdots\bigr\}. 
\]

(b) A {\em graded object} in $\D$ is an $\oo$-functor from $\ZZ$, considered as
a discrete category, into $\D$, i.e., simply a family of objects $A=(A_i)_{i\in\ZZ}$
in $\D$. 
\end{defi}

We have $\oo$-categories $\Filt(\D)$ and $\Grad(\D)$ formed by filtered and
graded objects in $\D$. They are connected by a pair of $\oo$-functors
\[
\xymatrix{
\Filt(\D)   \ar@<.4ex>[r]^{  \gr} &  \Grad(\D)  \ar@<.4ex>[l]^{  \eps }, 
}
\]
where $\eps$ considers any graded object as a filtered object with zero
morphisms, and $\gr$ is the  functor of
``taking the associated graded":
\[
\gr_i(Z) \,=\, \Cof\, \{Z_{i-1} \lra Z_i\}. 
\]
The following is an analog of Example \ref{ex:versions-sphere}(d).

\begin{prop}\label{prop:filt-and-gr}
The pair $(\gr, \eps)$ is a spherical adjoint pair. The associated twist
$T_{\Filt(\D)}$ is given by the shift of the filtration: $(T_{\Filt(\D)}(Z))_i = Z_{i-1}$. 
The cotwist
$T_{\Grad(\D)}$ is given  by $(T_{\Grad(\D)}(A))_i = A_{i-1}[1]$. 
\end{prop}

\noindent {\sl Proof:} The adjunction for  $(\gr, \eps)$ 
is obtained by defining the unit and counit maps
\[
u_Z: Z \lra \eps (\gr (Z)), \quad c_A: \gr(\eps(A))\lra A
\]
explicitly. That is, 
  $u_Z$ consists of the canonical maps $Z_i \to\Cof \{Z_{i-1}\to Z_i\}$,
while $c_A$ consists of the canonical maps $\Cof\{ A_{i-1}\buildrel 0\over\to A_i\} \to A_i$.  The fact that the adjunction is spherical, follows by 
straightforward calculation of the twist and cotwist and their identifications with what is
claimed in the proposition. 
  \qed


\paragraph{Filtered objects and chain complexes.} 
 In classical homological algebra, a filtration gives rise to
  a kind of chain complex formed by the  quotients of the filtration
 and the $\Ext^1$-classes between successive quotients. 
 We adapt this point of view to our case as follows.
 
 \vskip .2cm

Let $\ZZ^{\lhd}$ denote the poset obtained from $\ZZ$ by adjoining a minimal element which
we denote by $(-\infty)$. Let $P \subset \ZZ^\lhd \times \ZZ$ be the
subposet 
consisting of those pairs $(i,j)$ which satisfy $i \le j$. 
To avoid cluttering the
notation, in this section we will identify a poset with its nerve.

\begin{defi}\label{def:chain-cx} 
Let $\D$ be a stable $\infty$-category. 
A {\em chain complex in $\D$} is a functor
	$X: P \to \D$
	satisfying the following conditions:
	\begin{enumerate}
		\item For every $i \in \ZZ^\lhd$, the object $X(i,i)$ is a zero object.
		\item For every triple $(i,j,k)$ of elements of $\ZZ^\lhd$ with $i \leq k\leq j\leq l$, the square
			\[
				\begin{tikzcd}
					X(i,j) \ar{r}\ar{d} & X(i,l) \ar{d}\\
					X(k,j) \ar{r} & X(k,l)
				\end{tikzcd}
			\]
			is biCartesian in $\D$. 
	\end{enumerate}
	We denote by $\Ch(\D)$ the full subcategory of $\Fun(P,\D)$ spanned by the chain
	complexes. 
\end{defi}

\begin{prop}\label{prop:filtered}
The $\oo$-functor 
$\Ch(\D)
 \to \Filt (\D)= \Fun(\ZZ, \D)$
 given by restriction along the inclusion of posets
\[
	\alpha: \ZZ \hra P, \; i \mapsto (-\infty, i),
\]
is an equivalence. 
\end{prop}

\noindent{\sl Proof:} 
Let us factor the embedding
$\alpha: \ZZ\hra P$ into $\ZZ\buildrel\beta\over\hra Q \buildrel
\gamma\over\hra P$, where $Q$ is obtained by adding to
$\ZZ=\{(-\oo, i)\}$ elements of the form $(i,i)$, $i\in\ZZ$. 
We can identify $\Fun(\ZZ,\D)=\Fun'(Q,\D)$ where $\Fun'(Q,\D)\subset
\Fun(Q,\D)$ consists of functors equal to $0$ on each $(i,i)$.
One can seen this extension by $0$'s as an instance of $\beta_*$,
the right Kan
extension from $\ZZ$ to $Q$ (along $\beta$), since the undercategories appearing in the pointwise formulas, are empty.

Next consider $\gamma_!$, the left Kan extension  from $Q$ to $P$
and its effect on functors from $\Fun'(Q,D)$. If $Z$ is such a functor,
then the pointwise formula gives the value of $(\gamma_!Z)(i,j)$ 
as the colimit over the overcategory
 of $(i,j)$ in $Q$, and this colimit reduces
to
\[
\varinjlim \,\bigl\{ Z(i,i) \leftarrow Z(-\oo,i) \to Z(-\oo,j)\bigr\} \,=\,
\Cof\bigl\{ Z(-\oo,i) \to Z(-\oo,j)\bigr\}, 
\]
which is, naively, the ``quotient'' of the filtered object represented
by $Z$. It is clear that $\gamma_!Z$, i.e., the system of these ``quotients'', satisfies the
conditions of Definition \ref{def:chain-cx} so that $\Ch(\D)$
consists precisely of the functors which are left Kan extensions
of those from $\Fun'(Q,\D)$. Applying Proposition 
\ref{prop:kan-ext-equiv} finishes the proof. \qed
 

\paragraph{Chain complexes as coherence data.} 
To relate the data comprising a chain complex $X$ in a stable $\infty$-category to more traditional
terminology, we define the {\em $i$th term } of $X$ to be the object
$Y_i = X(i,i+1)[i]$. 
The connecting morphisms of the exact triangles
\[
	\begin{tikzcd}
		X(i,i+1) \ar{r}\ar{d} & X(i,i+2) \ar{d} \ar{r} & 0 \ar{d}\\
		X(i+1,i+1) \ar{r} & X(i+1,i+2) \ar{r} & X(i,i+1)[1]
	\end{tikzcd}
\]
induce morphisms $d: Y_{i+1} \to Y_{i}$, called the {\em differentials} of $X$. To interpret the
remaining data, we introduce the subposet $P(n) \subset P$ consisting of those pairs $(i,j)$ of
integers satisfying $j - i \le n$. With this notation we have
\begin{enumerate}
	\item The datum $X|P(1)$ is equivalent to the collection of terms $\{Y_i\}$ of the complex.
	\item The datum $X|P(2)$ is equivalent to the sequence 
		\[
			\dots \overset{d}{\to} Y_{i+1} \overset{d}{\to} Y_i \overset{d}{\to} Y_{i-1}
			\overset{d}{\to} \dots
		\]
		of differentials. 
	\item The datum $X|P(3)$ is equivalent to the sequence of differentials together with a
		homotopy $d^2 \simeq 0$ for every consecutive pair of differentials.
	\item[($n$)] The datum $X|P(n)$ is equivalent to the sequence of differentials together with a
		coherent system of homotopies $d^k \simeq 0$ for $2 \le k < n$. These 
		homotopies can be interpreted as trivializations of higher Massey products
		associated to the sequence of differentials.
		The triviality of these products is,
		 classically, the condition for the next level
		Massey products being defined. 
\end{enumerate}


\paragraph{Pseudo-complexes.}
We use the above analysis of chain complexes to
interpolate between filtered and graded objects.

\begin{defi}
Let $n \ge 0$. By a {\em pseudo-complex of order $n$} in
a stable $\oo$-category $\D$ we mean a functor
$X: P(n) \to \D$
satisfying:
\begin{enumerate}
	\item For every $i \in \ZZ$, the object $X(i,i)$ is a zero object.
	\item For every triple $(i,j,k)$ of integers with $i < j < k$ and $k-j \le n$, the square
		\[
			\begin{tikzcd}
				X(i,j) \ar{r}\ar{d} & X(i,k) \ar{d}\\
				X(j,j) \ar{r} & X(j,k)
			\end{tikzcd}
		\]
		is biCartesian in $\D$. 
\end{enumerate}
We denote by 
\[
	\Ch^{(n)}(\D) \subset \Fun(P(n), \D)
\]
the full $\oo$-subcategory consisting of pseudo-complexes of order $n$. 
\end{defi}

\begin{ex}
A pseudo-complex of order $1$ is the same as a graded object in $\D$,
i.e., $\Ch^{(1)}(\D)=\Grad(\D)$. 
A chain complex can be formally considered as a pseudo-complex of order 
$\oo$. 

\end{ex}

The following is a generalization of Proposition \ref{prop:filt-and-gr}.

\begin{prop}\label{prop:spherical}
 Let $\D$ be a stable $\infty$-category and $n \ge 0$. 
Let $F:   \Ch(\D) \to \Ch^{(n)}(\D)$ be  the functor
	given by restriction along $P(n) \subset P$. Then: 
	\begin{enumerate}
		\item[(a)]  After the identification of $\Ch(\D)$ with $\Filt(\D)$ given by Proposition 
		\ref{prop:filtered}, the functor $F$ associates to a filtered object 
			\[
			\dots \to Z_{i-1} \to Z_i \to Z_{i+1} \to \dots
			\]
			the pseudo--complex $X: P(n) \to \D$ with $X(i,j) = \Cof(Z_i \to Z_j)$.
			
		\item[(b)]  The functor $F$ has a right adjoint $G$ which associates to a pseudo--complex 
		$X : P(n) \to \D$,
			the  (complex corresponding to the) filtered object 
			\[
				\dots \to X(i-1,i-1+n) \to X(i,i+n) \to X(i+1,i+1+n) \to \dots.
			\]
		\item[(c)]  The adjoint pair $(F,G)$ is spherical. The corresponding spherical twists
			are given on objects by the formulas (in the first one we identify $\Ch(\D)$ with $\Filt(\D)$):
			\[
			\begin{gathered}
				T_{\Ch(\D)} = T_{\Filt(\D)}: (Z_i)_{i \in \ZZ} \mapsto (Z_{i-n}[1])_{i \in \ZZ}
		 \\
				T_{\Ch^{(n)}(\D)}: (X_{(i,j)})_{(i,j) \in P(n)} \mapsto
				(X_{(i-n,j-n)})_{(i,j) \in P(n)}. 
				\end{gathered} 
			\]
	 
	\end{enumerate}
\end{prop}

\noindent{\sl Proof:} Part (a) is immediate from the way the identification of Proposition 
		\ref{prop:filtered} is constructed. 
		
	\vskip .2cm
		
		To prove (b), 
let $p: \X=\Gamma(F) \to \Delta^1$ be the covariant Grothendieck construction of the functor $\Delta^1 \to
	\Cat_{\infty}$ classifying $F$, see Example \ref{ex:groth-single}.  	
	Let  also  $\L(\X) = \Fun_{\Delta^1}(\Delta^1, \X)$
	be the $\infty$-category of sections of $p$. Thus an object of $\L(\X)$
	is an edge $s$ of $\X$ covering the non-degenerate edge of $\Delta^1$. 
	
		\vskip .2cm
	
	We use two lemmas. First, 
	we describe $\L(\X)$ directly, 
	 using the description of  $F$ in the form given in (a). More precisely, 
denote by $Q(n) \subset P$ the subposet spanned by $P(n)$ together with the pairs $(-\oo,i)$. 

\begin{lem}\label{lem:groth-pseudocx}
  $\L(\X)$  is equivalent to the full subcategory of
	$\Fun(Q(n),\D)$ consisting of those diagrams 
	$X: Q(n) \to \D$
	whose restriction to $P(n)$ is a pseudo--complex. 
\end{lem}

\noindent{\sl Proof of Lemma \ref{lem:groth-pseudocx}:}
By definition, see Example \ref{ex:groth-single}, an object of $\L(\X)$
	is a datum of:
	\begin{itemize}
\item[($1$)]  A filtered object, denote it $Z$.  
\item[($2$)]  A pseudo-complex of order $n$, denote it $X'$.

\item[($3$)]  A morphism  $\phi: F(Z)\to X'$  in $\Ch^{(n)}(\D)$. 
\end{itemize} 
Let us compare this with the data contained in a functor $X:Q(n)\to\D$ as in the lemma. Such an 
$X$ gives
us, first of all:
 
\begin{itemize}
\item[($1'$)]  A filtered object, denote it $Z$, represented by the restriction of $X$ to $\ZZ=\{(-\oo,,i)\}$

\item[($2'$)]  A pseudo-complex, denote it $X'$, represented by the restriction of $X$ to $P(n)$. 
\end{itemize} 
These data match (1) and (2) above. The remaning data contained in $X$, are the
morphisms $\xi_{ijk}: X(-\oo,i) \to X(j,k)$ obtained as the values of $X$ on the morphisms
(inequalities) $(-\oo,i) \leq (j,k)$ between $\ZZ$ and $Q(n)$. The compatibillity of the $\xi_{ijk}$
with the values of $X$ on $\Mor(Q(n))$ implies that they factor through the
\[
\xi_{ij}: X(\-\oo,j) \lra X(i,j), \quad (i,j)\in Q(n).
\]
Since $X(i,i)=0$, we see that $\xi_{ij}$ factors through a morphism
\[
\phi_{ij}: F(Z)_{ij} \,=\,\Cof\bigl\{ X(-\oo,i) \to X(-\oo,j)\bigr\} \lra X(i,j).
\]
These morphisms account for the data (3) above. In this way we identify $\L(\X)$ with
the category of functors $X$ as in the lemma, at the level of objects. The
  remaining details are left to the reader. Lemma \ref{lem:groth-pseudocx} is proved. 
  
  \vskip .2cm
  
  Next, we use the lemma just proved to describe (co)Cartesian sections (edges) of $p$. 
  
  \begin{lem}\label{lem:cart-edges-X}
  Let $s$ be a section of $p$    and  $X: Q(n)\to\D$ be a functor corresponding to $s$
  by  Lemma \ref{lem:groth-pseudocx}. Then: 
  	\begin{enumerate}
		\item[(a)]  $s$ is a  coCartesian edge  if and only if, for every $i < j \le i+n$,
			the corresponding square
			\[
				\begin{tikzcd}
					X(-\oo,i ) \ar{r}\ar{d} & X(i, i) \ar{d} \\
					X(-\oo,j) \ar{r} & X(i,j)
				\end{tikzcd}
			\]
			is biCartesian in $\D$.
		\item[(b)]   $s$ is a Cartesian edge if and only if, for every $j \in \ZZ$, the
			morphism $X(-\oo ,j) \to X(j-n,j)$
			 is an equivalence. 
	\end{enumerate}
  \end{lem} 
  
  \noindent{\sl Proof of Lemma \ref{lem:cart-edges-X}:} (a) Recall the description of coCartesian
  edges of $\Gamma(F)$ in Example
  \ref{ex:groth-single}. Applying it to our case, we see that $s$ is coCartesian
   if and only if the morphism $\alpha: F(Z)\to X'$
  in the datum (3) in the proof of Lemma \ref{lem:groth-pseudocx}, is an equivalence.
  This is just an alternative formulation of the condition in (a).
  
  \vskip .2cm
  
  (b) Let $s$ be represented by a filtered object $Z$, a pseudo-complex $X'$ and a
  morphism $\alpha: F(Z)\to X'$, as above. Thus $s: Z\to X'$ as an edge of $\Gamma(F)$.
  The condition that $s$ is Cartesian means that $\X_{/s} \to \X_{/X'}$ is a trivial Kan fibration, cf. 
Example \ref{ex:implicit}(1). Let us consider this condition
 at the level of objects. Objects of $\X_{/s}$ are triangles ($2$-simplices) in $\X$
\[
\xymatrix{
W\ar[d]_\phi \ar[dr]^c&
\\
Z
\ar[r]_s&X'
}
\]
with $W\in \Filt(\D)$, while objects of 
$\X_{/X'}$ are arrows $c$ which are a priori parts of such triangles. So the condition means,
in particular, that any $c$ can be completed to such a triangle. Now $c$, by definition of $\Gamma(F)$,
consists of a morphism, call it $\beta$,  from $F(W)=(\Cof\{W_i\to W_j\})_{(i,j)\in P(n)}$ to
$X'=(X_{ij})_{(i,j)\in P(n)}$. So the condition at the level of objects implies that for
any such $W,\beta$ we can find a morphism of filtered objects $\phi: W\to Z$ fitting
into a triangle as above. Notice that we can take $W$ to consist of $W_j = X(j-n,j)$,
and this implies that the condition  stated in part (b) of our lemma, is necessary.
The sufficiency follows easily. Lemma \ref{lem:cart-edges-X} is proved.

\vskip .2cm

The two above lemmas imply  that $p$ is a biCartesian fibration and that the adjoint
functor $G$ of $F$ is as claimed in Proposition \ref{prop:spherical}(b), thus proving this part.
Finally, part (c) of Proposition \ref{prop:spherical}
now follows in a straightforward manner
by explicitly computing the spherical twists via the descriptions  of the
	Cartesian and coCartesian sections.

This finishes the proof of 
Proposition \ref{prop:spherical}.

	
\section{Semiorthogonal decompositions}\label{sec:SOD}

In this chapter, we  reformulate the notion of a spherical adjunction in terms of 
semiorthogonal
decompositions and deduce various corollaries that arise from this beautiful point of
 view, which was pioneered
in \cite{halpern-shipman}.

Our first goal will be to describe an approach to semiorthogonal decompositions 
via stable
$\infty$-categories. We start by recalling the standard setup \cite{BK:SOD, bondal-orlov}.

\subsection{Semiorthogonal decompositions of triangulated categories}\label{subsec:SOD-triang}

Let $\Cc$ be a  category.  
 A full subcategory $\Sc\subset\Cc$ will be called {\em strictly full}, if any object of $\Cc$
isomorphic to an object in $\Sc$, is in $\Sc$. Thus, every full subcategory $\Ac\subset\Cc$
has its {\em strictly full closure} $\Ac^{\on{sf}}\supset \Ac$  which is the full subcategory consisting of all objects isomorphic to
objects of $\Ac$. Clearly, the embedding $\Ac\subset \Ac^{\on{sf}}$ is an equivalence. 
For this reason will  sometimes not distinguish between full and stricly full subcategories. 

\vskip .2cm

Let $\Cc$ be a triangulated category. Let $\Bc, \Ac \subset\Cc$ be  full triangulated subcategories.
We say that $\Bc$ and $\Ac$ are {\em orthogonal}, if for any objects  $b\in\Bc$
we have $\Hom_\Cc(b,a)=0$. 

Suppose given one full triangulated subcategory $\Ac\subset \Cc$.   
The {\em left} and {\em right orthogonals} of $\Ac$ in $\Cc$ 
are the full subcategories $^\perp \Ac, \Ac^\perp\subset\Cc$ whose objects
are defined as follows:
\be\label{eq:orth-triang}
^\perp\Ac\,=\, \bigl\{ c\in\Cc\bigl| \, \Hom(c,a)=0,\,\,\forall \, a\in\Ac\bigr\}, 
\quad \Ac^\perp \,=\, \bigl\{ c\in\Cc\bigl| \, \Hom(a,c)=0,\,\,\forall \, a\in\Ac\bigr\}. 
\ee
They are strictly full  
triangulated subcategories in $\Ac$.

 \begin{defi}\label{def:SOD-triang} 

An ordered pair $(\Ac,\Bc)$ of full triangulated subcategories of
$\Cc$ is called a {\em semiorthogonal decomposition} (or {\em SOD}, for short)
 if 
\begin{enumerate}
	\item $\Bc$ and $\Ac$ are orthogonal.
	
	\item for every object $x \in \Cc$, there exists a distinguished triangle
		\begin{equation}\label{eq:extab}
				b \lra x \lra a \lra b[1]
		\end{equation}
		with $a \in \Ac$ and $b \in \Bc$.
\end{enumerate}

\end{defi}
\noindent
 Let $(\Ac, \Bc)$ be an SOD of $\Tc$.
Denote by $\EXT(\Ac, \Bc)$  
 the category of distinguished triangles in $\Cc$ of the form
\eqref{eq:extab}. We have the obvious
 forgetful functor $\Xi: \EXT(\Ac,\Bc) \to \Cc$. Condition (2)
implies that this functor is essentially surjective. Condition (1) further implies, that $\Xi$ induces
a bijection on isomorphism classes of objects. Therefore, $\Xi$ is very close to being an equivalence
of categories, but we cannot quite deduce that:
 the attempt to construct an inverse to $\Xi$ fails due
to the non-functoriality of cones. 

\vskip .2cm

It is straightforward to show that, if $(\Ac,\Bc)$ forms an SOD of $\Cc$, then
$\Bc$ is, essentially, uniquely determined by $\A$: 
the embedding $\Bc\subset {{}^\perp\Ac}$ is an
equivalence of categories, so $^\perp\Ac$ is the strictly full closure of $\Bc$. 
Similarly, $\Ac\subset\Bc^\perp$ is an equivalence, so $\Bc^\perp$ is the strictly full closure
of $\Ac$. 

 \vskip .2cm
 
  Vice versa, we can start with any full
triangulated subcategory $\Ac \subset \Cc$ and ask when the ordered pair $(\A,{^\perp}\A)$
 is  n SOD of $\Cc$. The condition for that is that 
 the inclusion functor $i: \Ac \to \Cc$ has a left
adjoint. The adjoint in that case associates to any $x\in\Cc$ the $a$-term in the triangle \eqref{eq:extab}, this term being (in that case) functorial in $x$. 
We call such subcategories $\Ac \subset \Cc$ {\em left admissible}.
Dually, the condition that $i$ has a right adjoint, means that $(\Ac^\perp, \Ac)$
is an SOD of $\Cc$. In this case we call $\Ac$
{\em right admissible}.

 If $i$ has both a left and a
right adjoint, in which case $\Ac \subset \Cc$ is called {\em admissible}, then the restriction of a right
adjoint to ${^\perp}\A$ defines a functor 
\[
	g: {^\perp}\A \to \A.
\]
A calculation similar to the above analysis of $p$ shows that, morally, we ought to be able to
recover the triangulated category $\Cc$ from the triple $(\A,{^\perp}\A, g)$. However, in trying to
prove this statement, we run into obstacles of the same kind as above.

\vskip .2cm 

In conclusion, the axiomatic framework of triangulated categories does not allow for a fully
satisfying framework of semiorthogonal decompositions. It is well-known that the above issues can be resolved by passing to a suitable context of enhanced triangulated categories. 
See \cite{kuznetsov-lunts} for a treatment using dg-enhancements. 
Our choice here is to
use stable $\infty$-categories. 


\subsection{Semiorthogonal decompositions of stable $\oo$-categories and gluing functors}
\label{sec:semi}

\paragraph{Orthogonal subcategories.} 
Let $\C$ be an $\oo$-category. A full $\oo$-subcategory $\Se\subset\C$ will be called
{\em strictly full}, if any object of $\C$ equivalent to an object of $\Se$, is in $\Se$.
As in the case of the usual categories in \S \ref{subsec:SOD-triang}, any full
$\oo$-subcategory $\A\subset \C$ has its strictly full closure $\A^{\on{sf}}\subset \C$ which is the
full $\oo$-subcategory
 spanned by all objects equivalent to objects of $\A$. Thus $\A \subset \A^{\on{sf}}$
 is an equivalence of $\oo$-categories. For this reason we will sometimes not distinguish
 between full and strictly full $\oo$-subcategories. 

\vskip .2cm

A full subcategory $\A$ of a stable $\infty$-category $\C$ is called {\em stable} if it contains a zero
object and is closed under forming fibers and cofibers. Such an $\A$ is, in particular, a stable
$\oo$-category in the intrinsic sense. In particular, $\h\A$ is a full triangulated subcategory
in the triangulated category $\h\C$. 

\vskip .2cm

Let $\C$ be a stable $\oo$-category.  Let $\B, \A\subset\C$ be two stable subcategories.
We say that $\B$ and $\A$ are {\em orthogonal} if  for any objects  $b\in\B$, $a\in\A$,
the space $\Map_\C(b,a)$ is contractible. Note that this is equivalent to the condition
that the subcategories $\h\B$ and $\h\A$ in the triangulated category $\h\C$
are orthogonal in the sense of \S \ref{subsec:SOD-triang}. 

\vskip .2cm
For a stable subcategory $\A\subset \C$ we define its {\em left} and {\em right orthogonals} to be
the full $\oo$-subcategories $^\perp\A, \A^\perp\subset\C$ whose objects are specified by:
\[
\begin{gathered}
^\perp\A\,=\, \bigl\{ c\in \C\bigl| \, \Map_\C(c,a) \text{ is contractible, }\,\,\forall \, a\in\Ac\bigr\}, 
\\
\A^\perp \,=\, \bigl\{ c\in \C\bigl| \, \Map_\C(a,c) \text{ is contractible, }\,\,\forall \, a\in\Ac\bigr\}. 
\end{gathered}
\]
They are strictly full stable subcategories in $\C$. Their associated triangulated subcategories
in $\h\C$ can be identified as follows:
\[
\h\, ({^\perp\A})\,=\, ^\perp(\h\A),\quad \h(\A^\perp) \,=\, (\h\A)^\perp,
\]
where the orthogonals in the right hand sides are understood as in \eqref{eq:orth-triang}. 
 

\paragraph{The $\oo$-categories $\co{\A,\B}$ and $\arr{\A,\B}$.} 

Let $\C$ be a stable $\infty$-category and
let $\A,\B$ be stable  subcategories of $\C$. We denote by 
\[
	\co{\A,\B} \subset \Fun(\Delta^1 \times \Delta^1, \C)
\]
the full subcategory spanned by  biCartesian squares in $\C$ of the form
\begin{equation}\label{eq:semidiag}
		\begin{tikzcd}
			b\arrow{r} \arrow{d} & x \arrow{d}\\
			0 \arrow{r} &a
		\end{tikzcd}
\end{equation}
such that $a$ is an object of $\A$, $b$ is an object of $\B$, and $0$ denotes a zero object in $\C$. We have the $\oo$-functor $\xi: \co{\A,\B}\to \C$ given by the projection to the vertex $x$. 

\vskip .2cm

 Let further
	\[
		\arr{\A,\B} \subset \Fun(\Delta^1, \C)
	\]
	denote the full subcategory spanned by edges $a \to b$ in $\C$ with $a$ in $\A$ and $b$ in
	$\B$. We have the $\oo$-functor $\Fib: \arr{\A,\B}\to\C$ which associated to a morphism
	its fiber.

\begin{prop}\label{prop:<AB>}

(a) The $\oo$-categories $\co{A,\B}$ and $\arr{\A,\B}$ are equivalent.

\vskip .2cm

(b) The equivalence can be chosen to as to send the $\oo$-functor $\xi$ to the $\oo$-functor 
$\Fib$. 

\vskip .2cm

(c) Assume that $\B, \A$ are orthogonal. Then the $\oo$-functor in (b) is fully faithful, i.e., it defines an equivalence of $\co{\A,\B}\simeq
\arr{\A,\B}$ with a full $\oo$-subcategory in $\C$. 
\end{prop}

\noindent{\sl Proof:} 
	Consider  the full subcategory
	\[
		\M \subset \Fun(\Delta^1 \times \Delta^2, \C)
	\]
	spanned by diagrams
	\[
		\begin{tikzcd}
			b\arrow{r} \arrow{d} & x \arrow{d} \arrow{r} & \arrow{d} 0\\
			0 \arrow{r} & a \arrow{r} & b' 
		\end{tikzcd}
	\]
	with both squares biCartesian, $a$ in $\A$, and $b, b'$ in $\B$. The functors 
	\[
		\co{\A,\B} \lla \M  \lra \arr{\A,\B}
	\]
	given by projecting to $x$ and $a \to b'$, respectively, are equivalences. The functor
	$\Fib$ can be identified with the composite
	\[
		\arr{\A,\B} \overset{\simeq}{\lra} \co{\A,\B}\buildrel\xi\over \lra \C
	\]
	which  proves (a) and (b). 	
	
	\vskip .2cm
	
	Let us prove (c) in the following form: the functor $\Fib: \arr{\A,\B}\to\C$
	induces weak equivalences on the $\Map$-spaces, i.e., for any two objects
	$a\buildrel f\over\to b$ and $a'\buildrel f'\over\to b'$ of $ \arr{\A,\B}$ the map
	 \[
	\Map_{\arr{\A,\B}}(f, f') \lra \Map_\C(\Fib(f), \Fib(f'))
	\]
	is a weak equivalence.  Note that we have a weak equivalence
	\[
	\Fib\{a\buildrel f\over\lra b\}\, \simeq\,  \Cof\{a[-1]\buildrel f[-1]\over\lra b[-1]\}.
	\]
	Therefore we have
	weak equivalences
	\[
	\begin{gathered}
	\Map_\C(\Fib(f), \Fib(f')) \,\simeq\, \Map_\C(\Cof(f[-1]), \Fib(f'))  \,\simeq\,
	\\
	\Fib^T\biggl\{ \Fib^T\bigl\{\Map_\C(b[-1],a') \to \Map_\C(b[-1], b')\bigr\} \lra
	\Fib^T\bigl\{ \Map_\C(a[-1], a') \to\Map_\C (a[-1], b')\bigr\}\biggr\} 
	\,\simeq \, 
	\\
	\Fib^T\biggl\{ \Fib^T\bigl\{\Map_\C(b,a') \to \Map_\C(b, b')\bigr\} \lra
	\Fib^T\bigl\{ \Map_\C(a, a') \to\Map_\C (a, b')\bigr\}\biggr\}, 
	 \end{gathered} 
	\]
where $\Fib^T$ means the homotopy fiber in the category of spaces. 	
	Now, since the $\Map_\C(b,a')$ is contractible,	the last iterated homotopy fiber above
	can be identified with the 
	homotopy fiber product
	\[
	 \Map_\C(b,b') \times^h_{\Map_\C(a,b')} \Map_\C(a,a')   \,\simeq \,
	 \Map_{\arr{\A,\B}}(f, f'). 
	\]
	This chain of identifications is clearly compatible with the one induced by the functor $\Fib$.
	This proves the statement (c). \qed

	
	\paragraph{Semiorthogonal decompositions.} 

\begin{defi}\label{def:oo-SOD}
 Let $\C$ be a stable $\infty$-category. We say that an ordered pair $(\A,\B)$ of full
	stable subcategories of $\C$ forms a {\em semiorthogonal decomposition} 
	(or an {\em SOD}, for short) of $\C$, if the
	  functor $\xi: \co{\A,\B} \to \C$ (or, equivalently, the functor
	  $\Fib: \arr{\A,\B}\to\C$) is an equivalence of $\infty$-categories.
\end{defi}

\begin{prop}\label{prop:perp}
	Suppose $\A$ and $\B$ are strictly full stable $\oo$-subcategories in $\C$.
	Suppose further that $(\A,\B)$ forms a semiorthogonal decomposition 
	 of a stable $\infty$-category $\C$.
	Then we have $\B = {^\perp}\A$ and $\A = \B^{\perp}$.
\end{prop}

\noindent{\sl Proof:}  Under the equivalence $\Fib: \arr{\A,\B} \to \C$, the subcategory
	$\A$ (resp. $\B$) of $\C$ is identified with the subcategory of $\arr{\A,\B}$ given by
	arrows of the form $a \to 0$ (resp. $0 \to b$). The mapping space $\Map_{\arr{\A,\B}}(0 \to
	b', a \to b)$ is equivalent to $\Map_{\B}(b',b)$ which is contractible for every $b'$ in
	$\B$ if and only if $b$ is a zero object. This implies $\A = \B^{\perp}$ and the first claim
	follows analogously.
\qed

\begin{cor}\label{SOD:oo.vs.triang}
A pair of stable subcategories
$(\A, \B)$ form an SOD of  a stable $\oo$-category $\C$ if and only if the pair of  full triangulated subcategories
$(\h\A, \h\B)$ form an SOD of the triangulated category $\h\C$ in the sense of
 \S \ref{subsec:SOD-triang}.
\end{cor}

The corollary means that we have a bijection between semiorthogonal decompositions of $\C$ 
as a stable $\oo$-category and of $\h\C$ as a triangulated category.

\vskip .2cm

\noindent{\sl Proof:} ``If'': Suppose that $(\h\A, \h\B)$ form an SOD of   $\h\C$. Then
$\h\B$ and $\h\A$ are orthogonal which implies that $\B$ and $\A$ are orthogonal.
Therefore, by Proposition \ref{prop:<AB>}(c)
 the $\oo$-functor $\xi$ is fully faithful and its essential image consists of
objects $x\in\C$ which include into a biCartesian square of the form
\eqref{eq:semidiag}. But by definition of the triangulated structure in $\h\C$,
existence of such a square is equivalent to existence of a distinguished triangle in $\h\C$
of the form \eqref{eq:extab}. Since $(\h\A, \h\B)$ form an SOD of   $\h\C$, such a triangle
does exist for any $x\in\Ob(\h\C) = \Ob(\C)$. Therefore a square \eqref{eq:semidiag}
exists as well. This means that the essential image of $\xi$ is all of $\C$, i.e.,
that $\xi$ is an equivalence. 

\vskip .2cm

``Only if'':  Suppose that $\xi$ is an equivalence. Then by Proposition \ref{prop:perp},
$\B$ and $\A$ are orthogonal, which means that $\h\B$ and $\h\A$ are orthogonal
in $\h\C$, thus fulfilling condition (1) of Definition \ref{def:SOD-triang}.
Further, essential surjectivity of $\xi$ together with the definition of
triangulated structure in $\h\C$, imply that condition (2) of the same definition
holds as well. \qed


\paragraph{Cartesian and coCartesian SODs.} 
Let $\C$ be a stable $\infty$-category, equipped with a semiorthogonal decomposition $(\A,\B)$. We
introduce a simplicial set $\chi(\A,\B)$ with a map $p:  \chi(\A,\B) \to \Delta^1$ as follows: An
$n$-simplex of $\chi(\A,\B)$ consists of
\begin{enumerate}
	\item[(1)] an $n$-simplex $j: [n] \to [1]$ of $\Delta^1$,
	\item[(2)] an $n$-simplex $\sigma: \Delta^n \to \C$ such that
	 $\sigma|\Delta^{j^{-1}(0)} \subset
		\A$ and $\sigma|\Delta^{j^{-1}(1)} \subset \B$.
\end{enumerate}
The map $p$ is the  forgetful map to part (1) of the data. 
 By construction, 
$p^{-1}(0)=\A$ and $p^{-1}(1)=\B$. 
It follows immediately from the inner horn filling property of
$\C$ that $\chi(\A,\B)$ is an $\oo$-category  (and this implies that $p$ it is an
inner fibration \cite[Prop.2.3.1.5]{lurie:htt}).

\vskip .2cm

The the $\infty$-category of sections of $p$ is identified as follows. 

\begin{prop}\label{prop:sections}
	There exists a canonical equivalence
	\[
		\Fun_{\Delta^1}(\Delta^1, \chi(\A,\B)) \,\simeq \, \C.
	\]
\end{prop}
 
 \noindent{\sl Proof:} 
Notice that $\Fun_{\Delta^1}(\Delta^1, \chi(\A,\B))$  is isomorphic as a simplicial
	set to the $\infty$-category $\arr{\A,\B}$. See Proposition \ref{prop:<AB>}(a)
	where it is proved that $\arr{\A,\B}$ is equivalent to  the $\oo$-category
	$\co{\A,\B}$  so  $\C$ is equivalent to them both by  the definition of an SOD. \qed

\begin{defi}\label{def:Cart-SOD}
A semiorthogonal decomposition $(\A,\B)$ of a stable $\oo$-category $\C$
is called {\em Cartesian} (resp. {\em coCartesian}), if the $\oo$-functor
$p: \chi(\A,\B)\to\Delta^1$ is a Cartesian (resp. coCartesian) fibration. 
\end{defi}

Recall (\S \ref{def:cart-fib}) that the condition of $p$ being a (co)Cartesian fibration is
formulated in terms of $p$-(co)Cartesian edges of $\chi(\A, \B)$. Let us reformulate
this condition using the specifics of our case.

\begin{defi}\label{def:(A,B)-cart}
	Let $\C$ be a stable $\infty$-category equipped with a semiorthogonal decomposition
	$(\A,\B)$. An edge $e: a \to b$ in $\C$ with $a$ in $\A$ and $b$ in $\B$ is called
	\begin{enumerate}
		\item[(1)]  $(\A,\B)$-{\em Cartesian} if $e$ defines a final object of the $\infty$-category
			$\C_{/b} \times_{\C} \A$.
		\item[(2)]  $(\A,\B)$-{\em coCartesian} if $e$ defines an initial object of the $\infty$-category
			$\C_{a/} \times_{\C} \B$.
	\end{enumerate}
\end{defi}

 Note that an edge  $e: a \to b$ in $\C$ with $a$ in $\A$ and $b$ in $\B$ is the same as
 as an edge of $\chi(\A,\B)$ covering the nontrivial edge of $\Delta^1$. 
 
 \begin{prop}\label{prop:cartesian}
  Let $\C$ be a stable $\infty$-category equipped with a semiorthogonal decomposition
	$(\A,\B)$. Let $e: a \to b$ in $\C$  be an edge with $a$ in $\A$ and $b$ in $\B$. Then the
	following are equivalent:
	\begin{itemize}
	\item[(i)] $e$ is  $(\A,\B)$-Cartesian (resp.  $(\A,\B)$-coCartesian).
	
	\item[(ii)] $e$, considered as an edge of $\chi(\A,\B)$, is $p$-Cartesian (resp.
	$p$-coCartesian). 	
	\end{itemize}
	\end{prop}
	
	\noindent{\sl Proof:} We prove the equivalence of the two versions of ``Cartesian''.
	The coCartesian case follows by duality. 
	We need to prove that $e$ is $p$-Cartesian if and only if the
	corresponding object in $\C_{/b} \times_{\C} \A$ is a final object. The edge $e$ being
	$p$-Cartesian is equivalent to the map 
	\[
		\chi(\A,\B)_{/e} \lra \chi(\A,\B)_{/b} \times_{\Delta^1_{/1}} \Delta^1_{/\{0 \to 1\}}
	\]
	being a trivial Kan fibration. Unravelling the definition, this map idenitifies with the map
	\[
		\C_{/e} \times_{\C} \A \lra \C_{/b} \times_{\C} \A
	\]
	which is a trivial Kan fibration if and only if $e$ defines a final object in $\C_{/b}
	\times_{\C} \A$.  \qed
	
	\begin{cor}
	A semiorthogonal decomposition $(\A,\B)$ is Cartesian 
	if and only if for every object $b$ of $\B$, there exists an
	$(\A,\B)$-Cartesian edge $a \to b$. Dually,  $(\A,\B)$ is coCartesian 
	if and only if for every object $a$ of $\A$, there exists an $(\A,\B)$-coCartesian edge $a \to b$. \qed
	\end{cor}


\paragraph{Gluing functors.} 

Let $\C$ be a stable $\infty$-category with Cartesian semiorthogonal decomposition
$(\A,\B)$. The Cartesian fibration $p: \chi(\A,\B) \to \Delta^1$ corresponds uniquely to an $\oo$-functor
 $G: \B \to \A$ so that $\chi(\A,\B)$ is identified with the contravariant Grothendieck
 construction $\chi(G)$. See Example \ref{ex:implicit}.  
 We call $G$ a {\em Cartesian gluing functor} of $(\A,\B)$.
 
 Dually,  given a coCartesian semiorthogonal
decomposition $(\A,\B)$, we similarly obtain an $\oo$-functor  $F: \A \to \B$ so that $\chi(\A,\B)$
is identified with the covariant Grothendieck construction $\Gamma(F)$. We
  will  refer to $F$ as a {\em coCartesian gluing functor} of $(\A,\B)$. 
  
  This analysis can be summarized as follows. 

\begin{prop}\label{prop:gluing} Let $\C$ be a stable $\infty$-category equipped with a 
	semiorthogonal decomposition $(\A,\B)$. 
	\begin{enumerate}
		\item[(1)] Assume that $(\A,\B)$ is Cartesian with Cartesian gluing functor $G: \B \to \A$. Then:
		\begin{enumerate}
			\item[(1a)] There is an equivalence of $\infty$-categories 
			\[ 
				\Fun_{\Delta^1}(\Delta^1, \chi(G)) \overset{\simeq}{\lra} \C.  
			\] 
			In particular, $\C$ can be uniquely recovered, up to equivalence, from the triple $(\A,\B,G)$.
			
			\item[(1b)] For two objects $a\in \A, b\in\B$, a morphism $a\to b$ in $\C$
			is $(\A, \B)$-Cartesian, if and only if the morphism $a\to G(b)$ 
			in $\A$, corresponding to it by (1a), is an equivalence.
			
			\end{enumerate}
		\item[(2)] Assume that $(\A,\B)$ is coCartesian with 
		coCartesian gluing functor $F: \A \to \B$. Then:
		\begin{enumerate}
		\item[(2a)] 	T here is an equivalence of $\infty$-categories 
			\[ 
				\Fun_{\Delta^1}(\Delta^1, \Gamma(F)) \overset{\simeq}{\lra} \C.  
			\] 
			In particular, $\C$ can be uniquely recovered, up to equivalence, from the triple $(\A,\B,F)$.
			
			\item[(2b)] For two objects $a\in \A, b\in\B$, a morphism $a\to b$ in $\C$
			is $(\A, \B)$-coCartesian, if and only if the morphism $F(a)\to b$ 
			in $\B$, corresponding to it by (2a), is an equivalence.
	\end{enumerate}
		\end{enumerate}
 
\end{prop}
\begin{proof} Parts (1a) and (2a) are immediate from Proposition \ref{prop:sections}  and  Example \ref{ex:implicit}
 (based on \cite[\S 3.3.2]{lurie:htt}). Parts (1b) and (2b) follow from the definition
 of the Grothendieck constructions. 
\end{proof}

\begin{rem}
Any $\oo$-functor $p: \X\to\Delta^1$, being an inner fibration, can be seen
as a {\em correspondence} between the  $\infty$-categories $p^{-1}(0)$ and $p^{-1}(1)$, see
 \cite[\S 2.3.1]{lurie:htt}. In these terms, 
 Proposition \ref{prop:sections} says that every semiorthogonal
decomposition $(\A,\B)$  of $\C$, (co)Cartesian or not,  arises from a correspondence between 
$\A$ and $\B$.  The Cartesian and coCartesian cases correspond to
the situation when this correspondence is given by a functor.
 
\end{rem}


\subsection{Admissiblity}\label{subsec:admissibility}

\paragraph{Admissibility and adjoints to the embedding.}
If $(\A,\B)$ forms an
SOD of a stable $\oo$-category $\C$ then, by 
  Proposition \ref{prop:perp},
the subcategories $\A$ and $\B$ determine each other
uniquely. It is therefore natural to ask which conditions a full stable subcategory $\A \subset
\C$ needs to satisfy so that $(\A,{^\perp}\A)$ forms an SOD. 

\begin{defi} Let $\C$ be a stable $\infty$-category, and let $\A \subset \C$ be a stable full
	subcategory. Let $i: \A \to \C$ denote the corresponding embedding functor. We call $\A
	\subset \C$
	\begin{itemize}
		\item {\em left admissible} if $i$ has a left adjoint,  denoted ${^*}i$; 
		\item {\em right admissible} if $i$ has a right adjoint, denoted $i^*$; 
		\item {\em admissible} if $i$ has both a left and a right adjoint.
	\end{itemize}
\end{defi}

\begin{prop}\label{prop:adm} Let $\C$ be a stable $\infty$-category and let $\A$ be a stable full
	subcategory. Then:
	\begin{enumerate}
		\item[(a)] The pair $(\A,{^\perp}\A)$ forms an SOD of $\C$ if and only if
			$\A \subset \C$ is left admissible. If these conditions hold, then we further
			have ${^\perp}\A = \Ker({^*}i)$ for any left adjoint  ${^*}i$ 
			of the inclusion $i: \A \to
			\C$.
			
		\item[(b)] The pair $(\A^\perp,\A)$ forms an SOD of $\C$ if and only if
			$\A \subset \C$ is right admissible. If these conditions hold,  then we further
			have $\A^\perp = \Ker(i^*)$ for any right adjoint $i^*$ of the inclusion $i: \A \to
			\C$.
	\end{enumerate}
\end{prop}

\noindent{\sl Proof:}  We show (a), the statement of (b) is dual.

\vskip .2cm

	Suppose that $i$ has a left adjoint ${^*i}: \C \to \A$. 
	Consider the triangulated category $\h\C$ and its full triangulated
	subcategories $\h\A$ and $\h({^\perp}\A)= {^\perp}\h\A$.
	 Note that
	$\Hom_{\h\C}(x,y) = \pi_0 \Map_\C(x,y)$. Therefore
	the equivalence 
	\[
		\Map_{\C}(x,i(a)) \simeq \Map_{\A}({^*i}(x),a)
	\]
	implies that the functor $\h\,  ({^*i}): \h\C\to\h\A$ is the left adjoint,
	in the usual categorical sense, of the
	embedding functor $\h i: \h\A\to\h\C$. This is a well known condition
	for $(\h\A, \h({^\perp}\A)= {^\perp} (\h\A))$ to form an SOD
	decomposition of $\h\C$ in the sense of triangulated categories.
	Therefore by Proposition \ref{SOD:oo.vs.triang}, $(\A, {^\perp \A})$
	form an SOD of $\C$ in the sense of stable $\oo$-categories.

	\vskip .2cm
	
	Conversely, suppose that $(\A, {^\perp\A})$ is an SOD of $\C$.
	By Definition \ref{def:oo-SOD}, cf. Proposition \ref{prop:<AB>}(a),
	this means that the functor $\Fib: \arr{\A, {^\perp\A}}\to\C$ is
	an equivalence. Under this equivalence, the embedding
	$i: \A\to\C$ is identified with the embedding 
	$i': \A \to  \arr{\A, {^\perp\A}}$ sending an object $a$ to the
	arrow $\{a\to 0\}$. But a left adjoint to $i'$ is immediately found to
	be the projection functor sending an arrow $\{a\to b\}$ to $a$. 
	\qed

\paragraph{ Admissibility and (co)Cartesian SODs.}

Admissibility is not only related to the existence of semiorthogonal decompositions but also to
their properties of being Cartesian and coCartesian:

\begin{prop}\label{prop:cartright} Let $\C$ be a stable $\infty$-category equipped with a
	semiorthogonal decomposition $(\A,\B)$, so that $\A$ is left
	admissible and $\B$ is right admissible. Then: 
	\begin{enumerate}
		\item[(a)] The semiorthogonal decomposition $(\A,\B)$ is Cartesian if and only if the
		subcategory $ \A \buildrel i_\A\over\hra \C$ is also right admissible. If these equivalent conditions
		hold, then the restriction to $\B$ of any right adjoint $i_\A^*: \C \to \A$ of $i_\A$ 
		is a  Cartesian gluing functor.
		
		\item[(b)]  The semiorthogonal decomposition $(\A,\B)$ is coCartesian if and only if the
		subcategory $\B \buildrel i_\B\over\hra \C$ is also left admissible. If these equivalent conditions
		hold, then the restriction to $\A$ of any left adjoint ${^*i}_\B: \C \to \B$ of $i_\B$ 
		is a coCartesian  gluing functor.
	\end{enumerate}
\end{prop}
\noindent{\sl Proof:}  We show (a), the statement of (b) is dual.

\vskip .2cm

	First, suppose $(\A,\B)$ is Cartesian. This means that
	the $\oo$-category
	$\chi(\A,\B)$ is identified, as an $\oo$-category over $\Delta^1$,
	with the contravariant Grothendieck construction
	$\chi(G)$, where $G: \B\to\A$ is a coCartesian gluing functor. 
Let $\L := \Fun_{\Delta^1}(\Delta^1, \chi(G))$ be the 
	$\infty$-category of sections of  $\chi(G)\to\Delta^1$.
	By definition, an object of $\L$ is a triple $(a,b,u)$,
	where $a\in\A$, $b\in\B$ are objects and $u: a\to G(b)$
	is a morphism in $\A$. Sometimes for simplicity we will
	write an object of $\L$ simply as $\{a\buildrel u\over\to G(b)\}$.

	Since, by definition, $\C$ is identified with the $\oo$-category
	of sections of $\chi(\A,\B)\to\Delta^1$,  we have an identification  
	$\C\=\L$. Note that this identification sends an object $a\in\A$
	to the triple  $(a,0,0)$ and an object $b\in\B$ to the triple
	$(0,b[1],0)$. Let  $\A'$ and $\B'$ be the full subcategories 
	of $\L$ spanned by triples of these types, so that $\A'\= \A$ and
	$\B'\=\B$.

	  \vskip .2cm
	  
	  Let $r: \X= \Gamma(i') \to \Delta^1$ be the covariant
	  Grothendieck construction of  the inclusion $i': \A' \subset \L$.
We need to show that $i'$ admits a
	right adjoint which is equivalent to the map $r$ being a Cartesian fibration. By definition, an edge $e$ of $\X$ covering the edge
$0 \to 1$ of $\Delta^1$ is a datum consisting of
objects $a,a'\in\A$, $b\in \B$ and a square  
\begin{equation}\label{eq:psicart}
		\begin{tikzcd}
			a' \ar{r} \ar{d} & a \ar{d}\\
			0 \ar{r} & G(b)
		\end{tikzcd}
	\end{equation}
	in $\A$. The following lemma is
	straightforward and its proof is left to the reader. 
	
	\begin{lem}
	An edge $e$ of $\X$ covering $0\to 1$ is Cartesian,
	if and only if the square \eqref{eq:psicart} is Cartesian
	(i.e., is biCartesian, since $\A$ is stable). \qed
	\end{lem}
	
	We conclude that $r$ is Cartesian, as any object $\{a\buildrel u\over
	\to G(b)\}$
	of $r^{-1}(1)= \L$ is the target of a Cartesian edge as above
	with $a'=\Fib(u)$. 
	
	\vskip .2cm
	
 Further, let us  determine the
	restriction to $\B'$ of any right adjoint, denote it $j$, of $i'$. 
	By construction, see Example \ref{ex:implicit},  the value 
	of $j$ on an object $L=\{ a\buildrel u\over
	\to G(b)\} $ of $\L$ is obtained as the source of a Cartesian edge
	terminating at $L$, i.e., as the fiber of $u$. If $L$ has $a=0$, then
	$\Fib(u)=G(b)[-1]$. By a similar analysis on all the simplices,
	we find that $j= G[-1]$. To bring this to a statement about the
	right adjoint of $i: \A\to\C$, we recall from the above
	 that the identification of
	$\B\subset\C$ with $\B'\subset \L$ is given by $b\mapsto
	(0, b[1], 0)$. 
	 Therefore, the restriction of the adjoint of the embedding
	 $\A \hra \L$ to $\B$ is
	simply given by $G$ so that we have verified
	the last claim of the proposition.
	
	\vskip .2cm

	Conversely, suppose that $i: \A \to \C$ admits a right adjoint,
	i.e., that 
	   the
	  covariant Grothendieck construction 
	$r:  \X=\Gamma(i) \to \Delta^1$ is Cartesian. 
 Let $\Y \subset \X$ be the full
	subcategory spanned by the objects of $\A = r^{-1}(0)$ and the vertices of 
	$\C =
	r^{-1}(1)$ that are contained in $\B$. Then it is immediate to verify that the
	restriction $r|_{\Y}: \Y \to \Delta^1$ is a Cartesian fibration and  
	the $r|_{\Y}$-Cartesian edges are
	precisely those whose image in $\X$ is $r$-Cartesian. Further, the simplicial set $\Y$ is
	isomorphic to $\chi(\A,\B)$ over $\Delta^1$. Thus $(\A,\B)$ is Cartesian.
 \qed

\begin{cor} Gluing functors are exact. \end{cor}
\begin{proof} By Proposition \ref{prop:cartright}, the gluing functor is given by the restriction of
	a right or a left adjoint.  So it is left  or left exact and therefore exact since it acts between
	stable $\oo$-categories. 
\end{proof}

\begin{rem}\label{rem:sections} Let $\C$ be a stable $\infty$-category equipped with a Cartesian
	semiorthogonal decomposition $(\A,\B)$ with right gluing functor $G: \B \to \A$. By Proposition
	\ref{prop:gluing}, we obtain a canonical equivalence 
	\[
	\phi: \C\buildrel \simeq\over\lra  \L =
	\Fun_{\Delta^1}(\Delta^1,\chi(G)).
	\]
	 Under this equivalence, we have the following
	identifications: 
	\begin{enumerate}
		\item $\phi(\A{^\perp}) \subset \L$ is the full subcategory
			spanned by the Cartesian sections. 
		\item $\phi(\A) \subset \L$ is the full subcategory
			spanned by sections of the form $a \to 0$ with $a \in \A$.
		\item $\phi(\B) \subset \L$ is the full subcategory
			spanned by sections of the form $0 \to G(b)$ with $b \in \B$.
	\end{enumerate}
	If we further assume that $(\A,\B)$ is coCartesian, then 
	\begin{enumerate}[resume]
		\item $\phi({^\perp}\B) \subset \L$ is the full subcategory
			spanned by the coCartesian sections. 
	\end{enumerate}
\end{rem}


\subsection{Mutation}\label{subsec:mutation}

\paragraph{Definition of the left and right mutation functors.}

Let $\C$ be a stable $\infty$-category and let $\A \subset \C$ be an admissible subcategory. 
This means that the embedding functor $i_\A: \A\hra\C$ has both a left adjoint 
${^*}i_\A: \C\to\A$ and a right adjoint $i_\A^*:\C\to\A$. 
By Proposition \ref{prop:adm}, we obtain that $(\A^\perp, \A)$ and $(\A, {^\perp}\A)$ form
semiorthogonal decompositions. This implies that:
\begin{itemize} 
\item $\A^\perp$ is left admissible, i.e., the embedding  $i_{\A^\perp}: \A^\perp\hra\C$
has a left adjoint ${^*}i_{\A^\perp}: \C\to\A^\perp$. 

\item ${^\perp}\A$ is right adissible, i.e., the embeddibg
$i_{{^\perp}\A}: {^\perp}\A\hra\C$ has a right adjoint $i_{{^\perp}\A}^*: \C\to {^\perp}\A$. 
\end{itemize}

\begin{defi}\label{defi:mutation}
We define the {\em left} and {\em right mutation $\oo$-functors} around $\A$ to be
the compositions
\[
\begin{gathered}
	\lambda_{\A}\, =\, ({^*}i_{A^\perp})\circ i_{{^\perp}\A}
	:  {^\perp}\A \lra \A^{\perp},
	\\
	\rho_{\A}\,=\, (i_{{^\perp}\A}^*) \circ i_{\A^\perp}
	: \A^{\perp} \lra {^\perp}\A.
	\end{gathered}  
\]
We will also 
  refer to the subcategory 
$\A^{\perp}$  as the {\em left mutation} of ${^\perp}\A$ around $\A$, and vice versa, 
refer to
${^\perp}\A$ as the {\em right mutation} of $\A^{\perp}$ around $\A$. 
\end{defi}

\paragraph{Mutation functors are equivalences.} 
By construction, $(\lambda_\A, \rho_\A)$ form an adjoint pair.  In fact, we have a stronger
statement. 

\begin{prop} \label{cor:mutation}  
$\lambda_\A$ and $\rho_\A$ are mutually quasi-inverse equivalences of $\oo$-categories. 
\end{prop}

 Before giving the proof, we introduce some terminology.  Let
 	\be\label{eq:M-sq}
		  \xymatrix{
			a'' \ar[r]^v \ar[d] _{ {  \textstyle  S\quad  =\quad  \quad\quad}} & a \ar[d]^u \\
			0 \ar[r] & a' 
		 }
		 \ee
 be a biCartesian square in $\C$. We will  say that $S$ is an {\em M-square}
 (short for {\em mutation square}), if
  $a''$ is an object of $\A^{\perp}$, while  $a$ is an object of $\A$, and $a'$ is an
			object of ${^\perp}\A$.

\begin{lem}\label{lem:flip}
	 Let $S$ be  a  biCartesian square in $\C$. Then the following are equivalent: 
	\begin{enumerate}
		\item[(i)]  $S$ is an M-square. 
		\item[(ii)]  The edge $a'' \to a$ is $(\A^{\perp},\A)$-coCartesian. 
		\item[(iii)]  The edge $a \to a'$ is $(\A,{^\perp}\A)$-Cartesian.
	\end{enumerate}
\end{lem}
 \noindent{\sl Proof of the lemma:} 
	We show (i) $\Leftrightarrow$ (iii), the equivalence  (i) $\Leftrightarrow$ (ii)
	follows by passing to opposites. 
	
	\vskip .2cm
	
	Note that 
 by Proposition \ref{prop:cartright},  
$(\A^\perp, \A)$ is a  coCartesian SOD  and $(\A, {^\perp}\A)$ is a Cartesian SOD. 
Set $\B =
	{^\perp}\A$ and write $i_\B=i_{{^\perp}\A}$. 
	
	\vskip .2cm
	
	Since  $(\A,\B)$ is Cartesian, we have a Cartesian gluing functor $G: \B \to \A$.
	We recall that $G$ is obtained by analyzing the Cartesian fibration
	 $p: \chi(\A,\B)\to\Delta^1$
	and identifying it with the contravariant Grothendieck construction
	$p: \chi(G)\to\Delta^1$.
	Thus a Cartesian edge $u: a\to a'$ with $a\in \A$,
	$a'\in\B$ can be, after composing with an equivalence in the source, identified
	with the canonical edge $G(a')\to a'$ in $\chi(G)$. 
	
	 At the same time $G$ is identified with
	$i_\A^* i_\B$ by Proposition \ref{prop:cartright}. The functor $i_\A^*$ 
	can be obtained from the diagram
	\[
	\C \buildrel \Cof\over  \lla \{\A^\perp, \A \} \buildrel \on{pr}_\A\over\lra \A. 
	\]
	Here $\Cof$ is an equivalence which we invert, and 
	 $\on{pr}_\A$ is the projection functor sending an arrow
	 $c\to a$ with $c\in \A^\perp$, $a\in\A$,  to $a$. This latter functor is obviously
	 right adjoint to the embedding $\A\to \{\A^\perp, \A \}$ sending
	 $a$ to $0\to a$. 
	 
	 Therefore $G(a')$, together with its canonical Cartesian edge to $a'$,  is found by lifting $a'$ to $ \{\A^\perp, \A \}$ via $\Cof$
	 and then projecting to $\A$ which means precisely forming an
	 M-square. Lemma is proved. \qed
	 
	 \vskip .2cm
	 
	 \noindent {\sl Proof of Proposition \ref{cor:mutation}:}
	 Let $\M \subset \Fun(\Delta^1 \times \Delta^1, \C)$ be the full subcategory formed
	 by M-squares. We have the projections
	 \[
	   \A^\perp \buildrel q''\over \lla \M \buildrel q'\over \lra {^\perp \A}. 
	 \]
	sending an M-square  as in \eqref{eq:M-sq} to $a'$ and $a''$ respectively.  
	We claim that both $q'$ and $q''$ are trivial Kan fibrations. Let us show this for
	$q'$, the case of $q''$ being similar. Indeed,  we can decompose $q''$ as
	the composition
	\[
	\M \buildrel q'_2 \over\lra \on{Cart} \buildrel q'_1\over\lra {^\perp\A},
	\]
where $\on{Cart}\subset \Fun(\Delta^1, \C)$ is the full $\oo$-subcategory
spanned by $(\A, {^\perp\A})$-Cartesian edges, $q'_2$ projects a square
\eqref{eq:M-sq} to $a\to a'$ and $q'_1$ projects $a\to a'$ to $a'$.
Now, $q'_1$ is a trivial Kan fibration since $(\A, {^\perp\A})$ is a Cartesian SOD,
and $q'_2$ is a trivial Kan fibration by the ``uniqueness of the triangle" property
\cite[Rem. 1.1.1.7]{lurie:ha}.  

Thus composing $q''$ with a section of $q'$, we get an equivalence $\A^\perp\to{^\perp\A}$
defined canonically up to a contractible choice. We finally note that
$\rho_{\A}= (i_{{^\perp}\A}^*) \circ i_{\A^\perp}$ can be identified with such
a composite equivalence. This can be seen by analyzing $i_{^\perp\A}^*$ in terms
of Cartesian edges of the covariant Grothendieck construction $\Gamma(i_{^\perp\A})$
in a way completely analogous to the proof of Proposition \ref{prop:cartright}.
We leave the details to the reader. Proposition  \ref{cor:mutation} is proved. \qed

 \paragraph{Mutation functors as  coordinate changes.}\label{par:mut-change}

	Let $\C$ be a stable $\infty$-category and $\A \subset \C$ an admissible subcategory. 
	By  Proposition \ref{prop:<AB>}(a), we have equivalences of $\infty$-categories
	\[
		\arr{\A^{\perp},\A} \overset{\fib}{\lra} \C \overset{\fib}{\lla} \arr{\A, {^{\perp}\A}}.
	\]
	Let  $\psi: \arr{\A^{\perp},\A} \buildrel \simeq\over\to  \arr{\A, {^{\perp}\A}}$
	be the resulting coordinate change
	equivalence. 
 Lemma \ref{lem:flip} allows for an explicit description of  $\psi$ as follows. Let 
	
	\begin{equation}\label{eq:m}
			\D \subset \Fun(\Delta^2 \times \Delta^2, \C)
	\end{equation}
	be the full $\oo$-subcategory
	consisting of diagrams of the form 
	\begin{equation}\label{eq:mdiag}
			\begin{tikzcd}
				x \ar{r}\ar{d} & b \ar{d}{!} \ar{r} & 0 \ar{d}\\
				a' \ar{r} \ar{d} & y \ar{r}{*}\ar{d} & b' \ar{d}\\
				0 \ar{r}  & a \ar{r} & x' 
			\end{tikzcd}
	\end{equation}
	such that
	\begin{itemize}
		\item $b$ in $\A^{\perp}$, $b'$ in ${^\perp}\A$, $a,a',y$ in $\A$, 
		\item all squares are biCartesian in $\C$,
		\item the edge $!$ is $(\A^{\perp},\A)$-coCartesian, and the edge $*$ is $(\A,
			{^{\perp}\A})$-Cartesian.
	\end{itemize}
	We have the projections
	\[
	\arr{\A^{\perp},\A} \buildrel p\over\lla \D \buildrel q\over\lra \arr{\A,
	{^{\perp}\A}},
	\]
  $p$ being the projection to the edge $b \to a$ and $q$ the projection to the edge $a' \to b'$. 
     
   	\begin{prop}\label{rem:basechange} 
(a) Both 	  $p$ and $q$ are trivial Kan fibrations. 

\vskip .2cm

(b) 	The  equivalence $\psi$ is  given by
	composing a section of $p$ with $q$.
	
	\vskip .2cm
	
	(c) The equivalence with $\C$ is given by
	projecting to the vertex $x$. 
	
\end{prop}
   
  \noindent{\sl Proof:} Parts (a) and (b)   follow immediately from Lemma \ref{lem:flip}.   To
  see (c), notice that
 the diagram \eqref{eq:mdiag} exhibits $x$ simultaneously
	as the fiber of the morphism $b \to a$ and the fiber of $a' \to b'$. \qed
	
	\paragraph{Relation between Cartesian and coCartesian  gluing functors.} 
Let	$(\A,\B)$ be a semiorthogonal decomposition (SOD)
of $\C$.
 Definition   \ref{def:Cart-SOD}  and Proposition \ref{prop:cartright} can be
  summarized as follows.
  
  \vskip .2cm
 
If $(\A,\B)$ is Cartesian, then we have a Cartesian  gluing functor $G: \B\to\A$, so that 
$\chi(\A,\B)$ is identified with the contravariant Grothendieck construction $\chi(G)$. 
In this case the pair $(\A^\perp, \A)$ also forms an SOD, which is coCartesian
and therefore has a coCartesian  gluing functor $\Phi: \A^\perp\to\A$. 

\vskip .2cm

If $(\A,\B)$ is coCartesian, then we have a coCartesian  gluing functor 
$F: \A\to\B$ so that
$\chi(\A,\B)$ is identified with the covariant Grothendieck construction $\Gamma(F)$. 
In this case the pair $(\B, {^\perp\B})$ also forms an SOD, which is Cartesian
and therefore has a Cartesian  gluing functor  $\Psi:  {^\perp\B}\to\B$.

 \begin{prop}\label{prop:L-R-glue}
 (a) If  $(\A,\B)$ is Cartesian, then we have an identification $\Phi \= G\circ\rho_\A$. 
 
 (b)  If  $(\A,\B)$ is coCartesian, then we have an identification  $\Psi\= F\circ\lambda_\B$.
 
 \end{prop}
 
 \noindent{\sl Proof:} We prove (a), since (b) is similar. For $a''\in\A^\perp$, the value
 $\Phi(a'')$ is, by Proposition \ref{prop:cartright}, found by including $a''$ into an M-square
 \eqref{eq:M-sq} and projecting to $a\in \A$. The value $\rho_\A (a'')$ is found by including
 $a''$ into (the same) M-square and projecting to $\B={^\perp\A}$, getting
 some $a'\in\B$.  Finally,
 for any  $a'\in \B$ (in particular, for the one just obtained)  the value of $G(a')$
 is found by including $a''$ into an M-square and projecting to $a$.
 Since we can take the same M-square to serve all three cases, the claim follows. \qed


 \subsection{$4$-periodic semiorthogonal decompositions and spherical adjunctions}
 \label{subsec:4-period}

In this section, we provide an alternative definition of spherical adjunctions in terms of
$4$-periodic semiorthogonal decompositions. This perspective is due to \cite{halpern-shipman}.

 \paragraph{ BiCartesian decompositions, orthogonals and adjoints.} 
 Let $\C$ be a stable $\infty$-category, and $(\A,\B)$ be a semiorthogonal decomposition (SOD)
of $\C$. We assume that $\A$ and $\B$ are strictly full, so $\B={^\perp\A}$ and $\A=\B^\perp$. 
As in other situations, we say that 
 $(\A,\B)$ of $\C$  is a {\em biCartesian} SOD, if it is both Cartesian and coCartesian. 
 In this case
the $\oo$-functor $p: \chi(\A,\B)\to\Delta^1$, see 
  Definition   \ref{def:Cart-SOD}, is a biCartesian fibration.  Further, we have the
  Cartesian  gluing
  functor $G:\B\to\A$ (coming from the fact that $(\A,\B)$ is Cartesian) and
  the coCartesian gluing functor $F: \A\to\B$ (coming from the fact that $(\A,\B)$ is coCartesian). 
  These functors are adjoint to each other: $(F,G)$ is an adjoint pair and $\chi(\A,\B)\to \Delta^1$
  provides the adjunction datum. 
 Moreover, 
every consecutive pair in the
sequence 
\[
	\A^{\perp}\quad,\quad \A\quad,\quad \B={^\perp\A}\quad,\quad {^\perp}\B={^{\perp\perp}}\A
\]
forms an SOD, more precisely, $(\A^\perp, \A)$ is coCartesian, $(\A,\B)$ is biCartesian
and $(\B, {^\perp\B})$ is Cartesian. 
These categories are related by the gluing functors:
\[
\xymatrix{
\A^\perp \ar[r]^\Phi& 
\A \ar@<.4ex>[r]^{F } &\B  \ar@<.4ex>[l]^{G }& ^{\perp}\B,\ar[l]_\Psi
}
\]
with $\Phi, F$ being coCartesian gluing functors of coCartesian SODs and
 $G, \Psi$ being
right Cartesian  functors of Cartesian SODs. 

\begin{prop}\label{prop:bicart-adj}
\begin{itemize}
\item[(a)] Let $(\A,\B)$ being a biCartesian SOD. Then:
\begin{itemize}
\item[(a1)] The coCartesian gluing functor $\Phi$ of $(\A^\perp, \A)$ is identified, by composing
with the mutation functor $\rho_\A$, with the right adjoint of the coCartesian
 gluing functor $F$
for  $(\A,\B)$, i.e., $\Phi=G\circ\rho_\A$. 
 
 \item[(a2)] The Cartesian  gluing functor  $\Psi$ for $(\B, {^\perp\B})$ is identified, by
 composing with the mutation functor $\lambda_\B$, with the left adjoint of the
 Cartesian gluing functor $G$ for $(\A,\B)$, i.e., $\Psi=F\circ\lambda_\B$. 
 
\end{itemize}

\item[(b)] Let $(\A,\B)$ be a coCartesian SOD of $\C$, with the corresponding gluing functors
$\A\buildrel F\over\lra \B \buildrel \Psi\over\lla {^\perp\B}$. Then $(\A,\B)$ is
biCartesian, if and only if $F$ has a right adjoint.

\item[(c)]  Let $(\A,\B)$ be a Cartesian SOD of $\C$, with the corresponding gluing functors
$\A^\perp\buildrel\Phi\over\lra \A \buildrel G\over\lla \B$.  Then $(\A,\B)$ is
biCartesian,  if and only if $G$ has a left adjoint.

\end{itemize}
\end{prop}

\noindent{\sl Proof:} Part (a) follows from Proposition \ref{prop:L-R-glue}.  Let us prove
(b), part (c) being similar. The ``only if'' part follows since  for a biCartesian SOD
$(\A,\B)$, the biCartesian fibration
$\chi(\A,\B)$ provides an adjoint pair of functors. So see the ``if'' part, we identify
$\chi(\A,\B)$ with the covariant Grothendieck construction $\Gamma(F)$. So the existence
of a right adjoint to $F$ means that $\Gamma(F)\to\Delta^1$ is a biCartesian fibration and
so $(\A,\B)$ is a bicartesian SOD. \qed

\vskip .2cm

 One advantage of expressing an adjunction $p: \X \to \Delta^1$ of stable $\infty$-categories in
terms of a semiorthogonal decomposition of the stable $\infty$-category $\C =
\Fun_{\Delta^1}(\Delta^1, \X)$ consists of the accessibility of {\em further} iterated adjoints by
mutating the semiorthogonal decomposition to the left and the right.  
To illustrate this, we use the following concepts. 

\begin{defi}
(a) An {\em admissible chain of orthogonals} in a stable $\oo$-category $\C$ is a sequence
$(\A_i)_{i\in \ZZ}$ formed by strictly full stable subcategories such that for each
$i$ we have $\A_i = {^\perp\A}_{i-1}=\A_{i+1}^\perp$ and each pair $(\A_i, A_{i+1})$
forms a semiorthogonal decomposition of $\C$. 

\vskip .2cm

(b) A semiorthogonal decomposition $(\A,\B)$ of $\C$ is called $\oo$-{\em admissible},
if it can be included into a (necessarily unique)
 admissible chain of orthogonals $(\A_i)_{i\in\ZZ}$
so that $\A=\A_0$ and $\B=\A_1$. 
\end{defi}

 Proposition \ref{prop:bicart-adj} implies at once the following.
 
 \begin{cor}\label{cor:oo-adm-adj}
 Let $(\A,\B)$ be a semiorthogonal decomposition  of $\C$. The following
 are equivalent:
 \begin{itemize}
 \item[(i)] $(\A,\B)$ is $\oo$-admissible.
 
 \item[(ii)] $(\A,\B)$ is coCartesian and the corresonding coCartesian
  gluing functor $F:\A\to\B$ admits
 iterated left and right adjoints of all orders. 
 
  \item[(iii)] $(\A,\B)$ is Cartesian and the corresonding Cartesian
  gluing functor $G:\B\to\A$ admits
 iterated left and right adjoints of all orders. 
 \end{itemize} \qed
 \end{cor}
 
 
 \paragraph{$4$-periodic SODs.}

\begin{defi} Let $\C$ be a stable $\infty$-category. An admissible chain of orthogonals
$(\A_i)_{i\in \ZZ}$ is called {\em $4$-periodic}, if $\A_i=\A_{i+4}$ for all $i$.
An SOD $(\A,\B)$ of $\C$ is called  {\em $4$-periodic}, if 
it includes into a $4$-periodic admissible chain of orthogonals $(\A_i)$
so that $\A=\A_0, \B=\A_1$. 
\end{defi}

\begin{prop}
Let $(\A,\B)$ be an SOD of a stable $\oo$-category $\C$.
The following are equivalent:
\begin{itemize}
\item[(i)] $(\A,\B)$ is $4$-periodic.

\item[(ii)]  $(\A,\B)$ is biCartesian and the pair $({^\perp}\B, \A^{\perp})$ also forms
	an SOD of $\C$. \qed
\end{itemize}
 \end{prop}
 
 In the next paragraph we give a criterion for a biCartesian SOD to be $4$-periodic.


\paragraph{The relative suspension functor of a biCartesian SOD.} 
Let $\C$ be a stable $\infty$-category equipped with a biCartesian SOD
$(\A,\B)$. Denote  $\A' = {^{\perp}\B}$ and $\B' = {\A^{\perp}}$. 
By a {\em bimutation cube} ({\em MM-cube} for short) we will mean a
 cubical diagram, i.e., a functor $Q\in  \Fun((\Delta^1)^3,\C)$ of the form
\begin{equation}\label{eq:cubical}
		\xymatrix@=0.4cm{ & x' \ar[rr]\ar[dd] & & y' \ar[dd]\\ 
			x \ar[ur]\ar[rr]\ar[dd]_{\textstyle{Q\quad = \quad\quad}
			}
			&   &  y \ar[dd]\ar[ur] & \\
			& 0 \ar[rr] & & z' \\ 
		0 \ar[ur] \ar[rr]&  & z \ar[ur] &  }
\end{equation}
satisfying the following conditions:
\begin{enumerate}
	\item[(MM1)] the front and back faces of $Q$ are biCartesian squares
	  in $\C$,
	\item[(MM2)]  the vertices $x$, $x'$ are objects of $\A$, the vertices $y$, $y'$ are objects
		of $\B$,
	\item[(MM3)]  the morphism $x \to y$ is coCartesian (and therefore $z\in \A'$),  
	while the morphism $x' \to y'$ is Cartesian (and therefore $z'\in\B'$). 
\end{enumerate}

Let  $\K \subset \Fun((\Delta^1)^3,\C)$ be  the full subcategory spanned by  MM-cubes.

\begin{lem}\label{lem:cube} Consider the 
 diagram of $\infty$-categories
\[
		\C \overset{\Fib}{\longleftarrow} \arr{\A,\B} \overset{p}{\longleftarrow} \K
		\overset{q}{\longrightarrow} \arr{\A',\B'} \overset{\Fib}{\longrightarrow} \C
		\]
where $p$  projects an MM-cube $Q$ to the diagonal arrow  $x \to y'$, 
while $q$ projects $Q$ to the arrow $z \to z'$. 
As usual,  $\Fib$ denotes the functor   associating to a morphism its fiber. 
	In this diagram, the functors $p$ and $q$ are trivial Kan fibrations.
\end{lem}
\noindent{\sl Proof:} Consider first $p: \K\to \arr{\A,\B}$. Let 
$\E\subset \Fun((\Delta^1)^2,\C)$ be the full subcategory formed by possible
upper squares of MM-cubes, i.e., by diagrams
\[
\xymatrix{
x'\ar[r]&y'
\\
x\ar[u] \ar[r]&y\ar[u]
}
\]
with $x,x'\in\A$, $y,y'\in\B$ the  arrow $x \to y$  coCartesian and 
	 the arrow $x' \to y'$  Cartesian. We factor $p$ as the composition
	 $\K \buildrel p_2\over\to \E \buildrel p_1\over\to   \arr{\A,\B}$.
	 Any arrow $x\to y'$ can be represented, in an essentially unique way,
	 as a composition of a coCartesian arrow $x'\to y'$ and an arrow $x\to x'$, 
	 as well as a composition of an arrow $y\to y'$ and a Cartesian arrow
	 $x\to y$. This fact upgrades, in a standard way, see  \cite[4.3.2.15]{lurie:htt} 
	 to the conclusion that $p_1$ is a trivial Kan fibration. Further,  
	 $p_2$ is a trivial Kan fibration by a similar upgrade of the ``uniqueness of triangle'' property. 
	 
	 \vskip .2cm
	 
	 Look now at $q: \K\to \arr{\A', \B'}$. Using the fact the $\Cof: \arr{\A,\B}\to\C$
	 is an equivalence, we get that $\Cof$ induces an equivalence
	 $\Cof_1: \Fun(\Delta^1,  \arr{\A,\B})\to \Fun(\Delta^1, \C)$. We consider
	 the preimage, under $\Cof_1$, of $\arr{\A', \B'}\subset \Fun(\Delta^1, \C)$.
	 This preimage is easily identified with the category $\E$ of possible
	 upper parts of MM-cubes. Therefore $q$ is a trivial Kan fibration as well.\qed
	 
	  \vskip .5cm

	  Note that  $\Fib: \arr{\A,\B} \to \C$ is an equivalence by definition of an SOD and
	  Proposition \ref{prop:<AB>}(a), but  $\Fib: \arr{\A',\B'} \to \C$ is not necessarily an
	  equivalence, since $(\A', \B')$ may not be an SOD. 
	 
\begin{defi}\label{def:rel-susp}
The {\em relative suspension functor} of a biCartesian SOD is the
(essentially unique)  functor $\tau: \C\to\C$ obtained as the composite 
of an inverse of $\Fib \circ p$ with $\Fib \circ q$ in the diagram of
 Lemma \ref{lem:cube}. 
\end{defi}

\begin{exa}
	 Consider the trivial semiorthogonal
	decomposition with $\A = \C$ and $\B \subset \C$ the full subcategory spanned by the zero
	objects. Then  
 $\tau: \C \to \C$ is the usual suspension functor of the stable $\oo$-category $\C$. 
\end{exa}

\begin{cor}\label{cor:sphericaltau} Let $\C$ be a stable $\infty$-category equipped with a
	biCartesian SOD  $(\A,\B)$. Then the following conditions are
	equivalent:
	\begin{enumerate}
		\item[(i)]  $(\A,\B)$ is $4$-periodic.
		\item[(ii)] The relative suspension $\tau: \C \to \C$ is an equivalence.
	\end{enumerate}
\end{cor}
\begin{proof}
Indeed, (i) is equivalent to the property that $(\A', \B')$ is an SOD. At the same
time, (ii) is equivalent to the property that $\Fib: \arr{\A', \B'}\to\C$ is an equivalence,
by the two out of three property of equivalences, since all the other functors
in the diagram of  Lemma \ref{lem:cube} are equivalences. Now, the 
property that $\Fib: \arr{\A', \B'}\to\C$ is an equivalence is equivalent to 
$(\A', \B')$ being an SOD by Definition \ref{def:oo-SOD}. 
  \end{proof}
  
  
  \paragraph{$4$-periodicity and spherical adjunctions.} 
 Let $\C$ be a stable $\infty$-category equipped with a biCartesian
	semiorthogonal decomposition $(\A,\B)$.  

\begin{prop}\label{prop:taures}The relative suspension functor 
	$\tau$
	restricts to the functors of $\infty$-categories given by the following formulas:
		\[
 \rho_{\B}: \A \to \A',\; x \mapsto \fib(x \to f(x))[1],
 \]
 
			 where $x \to f(x)$ is a coCartesian morphism in $\C$,
			\[
			\lambda_{\A}[1]: \B \to \B',\; y \mapsto \fib(g(y) \to y)[1],
			\]
		 where $g(y) \to y$ is a Cartesian morphism in $\C$,
			 \[
			 \A' \to \A,\; \fib(x \to f(x)) \mapsto \fib(x \to g(f(x)))[1],
			 \]
		 where $x \to g(f(x))$ is a unit morphism, 
			\[
			B' \to \B,\; \fib(g(y) \to y) \mapsto \fib(f(g(y) \to y),
			\]
	 where $g(f(y)) \to y$ is a counit morphism. 
\end{prop}
\begin{proof} This can be directly read off from \eqref{eq:cubical}.
\end{proof}

 \vskip .2cm
 
 Further, let 
$
	F: \A \llra \B: G
$
be the left and right gluing functors for $(\A,\B)$,
so $(F,G)$ form an adjoint pair, with $p: \chi(\A,\B)\to\Delta^1$ providing the adjunction. 

\begin{prop}\label{prop:sph=4per}
	Let $\C$ be a stable $\infty$-category with a biCartesian semiorthogonal decomposition
	$(\A,\B)$. Then the following are equivalent:
	\begin{enumerate}
		\item[(i)]  $(\A,\B)$ is $4$-periodic.
		\item[(ii)]  The adjunction $p: \chi(\A,\B) \to \Delta^1$ is spherical.
	\end{enumerate}
\end{prop}

\noindent{\sl Proof:} 
By Proposition \ref{prop:taures}, the functor $\tau^2[-1]$ stabilizes the subcategories
	$\A$, $\B$, as well as $\A'$ and $\B'$. 
By Remark \ref{rem:sections}, this implies that
	$\tau^2[-1]$ induces an endofunctor of the biCartesian fibration $\pi: \chi(\A,\B) \to
	\Delta^1$ that preserves Cartesian and coCartesian edges. By \cite[3.3.1.5]{lurie:htt}, the
	functor $\tau^2[-1]$ is an equivalence if and only if its restrictions to $\A$ and $\B$ are 	equivalences.
But the restriction of $\tau^2[-1]$ to $\A$ and $\B$
recovers the twist and cotwist functors from \eqref{eq:twist} and  \eqref{eq:cotwist}, respectively,
defined as
\[
	T_{\A}: \A \lra \A, \; x \mapsto \Cof\{ x \to G(F(x))\},
\]
where $x \to G(F(x))$ is a unit morphism, and
\[
	T_{\B}: \B \lra \B, \; x \mapsto \Fib\{ F(G(y)) \to y\} 
\]
where $F(G(y)) \to y$ is a counit morphism. 
So  $\tau^2[-1]$ is an equivalence if and only if i	
  the spherical twist and cotwist are equivalences. But $\tau^2[-1]$ is
	an equivalence if and only if $\tau$ is an equivalence so that the claimed statement follows
	from Corollary \ref{cor:sphericaltau}. \qed
	
	
	\paragraph{The spherical functor package.}

  \begin{cor}\label{cor:left-sph-right}
  The concepts of left and right spherical $\oo$-functors  (Definition \ref{def:sphad})
 are equivalent. 
 \end{cor}
 
 \noindent{\sl Proof:} 
  Indeed, let $(F: \A \leftrightarrow B: G)$ be an adjoint pair of spherical $\oo$-functors so $F$ is left spherical and $G$ is right spherical. 
 By Proposition \ref{prop:sph=4per}, 
 it gives rise to a stable $\oo$-category $\C$ with a $4$-periodic
 SOD $(\A,\B)$ for which $F$ is a coCartesian gluing functor and $G$ is
 a Cartesian gluing functor. Now $(\B, \A'={^\perp\B})$ is also a $4$-periodic
 SOD of $\C$, a part of the same admissible chain of orthogonals, and 
 Proposition \ref{prop:bicart-adj} implies that its {\em Cartesian} gluing
 functor is identified with $F$, so $F$ is right spherical. Similarly,
 $(\Bc'=\A^\perp, \A)$ is also a $4$-periodic SOD of $\C$ and its {\em coCartesian}
 gluing functor is identified with $G$, so $G$ is left spherical. \qed
 
 \vskip .2cm
 
  We will therefore use the term {\em spherical $\oo$-functor}.

 \begin{cor}\label{cor:sph-oo-adjoints}
 Any spherical  $\oo$-functor $F: \A\to\B$ has iterated 
 left and right adjoints of all orders, each of which
 is a spherical $\oo$-functor itself. 
 \end{cor}
 
 \noindent{\sl Proof:} This follows from Proposition \ref{cor:oo-adm-adj},
   as a $4$-periodic SOD is $\oo$-admissible. \qed
   
   \vskip .2cm
   
   We will denote the iterated adjoints by
   \[
   \cdots {^{** }F}, {^*F}, F, F^* = G, F^{**} = G^*, \cdots
   \]

\begin{lem}\label{lem:grid} Let $(\A,\B)$ be a $4$-periodic semiorthogonal decomposition of a stable
	$\infty$-category $\C$. Consider the $\infty$-category $\D$ defined 
in \eqref{eq:m}, so an object of $\D$ is a diagram as in 
\eqref{eq:mdiag}. For such a diagram, 
	the following are equivalent:
	\begin{enumerate}
		\item[(i)]  The edge $b \to a$ is $(\A^{\perp},\A)$-Cartesian.
		\item[(ii)]  The edge $a' \to b'$ is $(\A, \B)$-coCartesian.
	\end{enumerate}
	In other words, the base change equivalence 
	\[
		\psi: \arr{\A^{\perp},\A} \simeq \arr{\A,\B}
	\]
	identifies Cartesian with coCartesian edges. 
\end{lem}
\begin{proof}
	This is an immediate consequence of Lemma \ref{lem:flip}: The edge $b \to a$ being
	$(\A^{\perp},\A)$-Cartesian is equivalent to $x$ being an object of $\A^{\perp\perp}$. The
	edge $a' \to b$ being $(\A, \B)$-coCartesian is equivalent to $x$ being an object of
	$^{\perp}\B$. But due to $4$-periodicity, we have $\A^{\perp\perp} = {^{\perp}\B}$.
\end{proof}

\vskip .2cm

Usually, the definition of the spherical functor considered in the
literature includes various additional conditions, see   \cite[Th. 1.1] {AL:spherical}. 
The following corollary shows that these
additional conditions are actually consequences 
 of  the requirement that the twist and cotwist are
equivalences.

\begin{cor}[(The spherical functor package)]
Let $(F: \A \leftrightarrow B: G)$ be an adjoint pair of spherical $\oo$-functors, so $G=F^*$ and $F={^*G}$. Then:  
 
	\begin{enumerate}
	 
		\item[(1)]  The cotwist with respect to  the adjoint pair
		of spherical functors $(G, G^* = F^{**})$ is an inverse to the twist with
			respect to $(F,G)$.
		\item[(2)] The twist with respect to $({^*F} = {^{**}G} , F)$ is an inverse to the cotwist with
			respect to $ (F,G)$.
		\item[(3)]  There is a canonical equivalence $F^{**} \simeq F T^{-1}_{\A}$.
		\item[(4)]  There is a canonical equivalence ${^{**}G} \simeq T_{\A}^{-1} G$.
	\end{enumerate}
\end{cor}

\noindent{\sl Proof:} 
 To show (1), consider a diagram of the form \eqref{eq:mdiag} satisfying the equivalent conditions of Lemma
\ref{lem:grid}. The condition on $b \to a$ to be Cartesian implies that the morphism $y \to a$ is a
counit morphism for the adjoint pair $(G,G^*)$ modelled by the biCartesian semiorthogonal
decomposition $(\A^{\perp},\A)$. In particular, we have $y \simeq G(G^*(a))$, and the exact triangle
\begin{equation}\label{ex:triangle}
		\begin{tikzcd}
			a' \ar{r} \ar{d} & y\ar{d} \\ 
			0 \ar{r}  & a 
		\end{tikzcd}
\end{equation}
exhibits $a'$ as a cotwist of $a$, with respect to $(G, G^*)$. The condition that $a' \to b'$ be
coCartesian implies that the morphism $a' \to y$ is a unit for the adjoint pair
 $(F,G)$ modelled
by $(\A,\B)$. In particular, the exact triangle \eqref{ex:triangle} exhibits $a$ as a twist of $a'$
with respect to the   $(F,G)$. We deduce that the cotwist with respect to 
$(G,G^*)$ and the
twist with respect to $(F,G)$ are inverse to one another. This shows (1),
part (2) is similar. 

Similarly, the statement $F \simeq G^*  T_{\A}$ follows from direct inspection of  \eqref{eq:mdiag}
showing (3). Claim (4) then follows from (3) by passing to left adjoints.
\qed

We conclude by noting that, building on the above results, a version of the $2$-out-of-$4$ property,
describing equivalent characterizations of spherical adjunctions of stable $\infty$-categories, has
been established in \cite{christ}, generalizing the results of \cite{AL:spherical} for dg-categories.


\section{The spherical S-construction}\label{sec:spher-S} 

In this chapter, we show that Waldhausen's relative S-construction admits a canonical
paracyclic structure when applied to a spherical functor of stable $\infty$-categories. The idea
that lies behind our proof is that the spherical adjunctions from Proposition \ref{prop:spherical}
organize into a coparacyclic object in the category $\Sph$ of spherical adjunctions which
corepresents the relative S-construction. For technical convenience, we provide a more direct
combinatorial proof. We start with some reminders. 

\subsection{The relative S-construction}\label{subsec:rel-S}
\paragraph{Cyclic and paracyclic objects.} 
Along with the simplex category $\Delta$, we will use Connes' {\em cyclic category}
$\Lambda$ and its universal cover, the {\em paracyclic category} $\Lambda_\oo$
introduced by Elmendorf \cite{elmendorf} under a different name. 
Here is a brief reminder, see
\cite[\S 2.1]{dkss:1} for a detailed survey. 

Objects of $\Lambda$ and $\Lambda_\oo$ are formal symbols $\cn$, $n\in\ZZ_+$
which can be thought of as finite cyclic ordinals $\ZZ/(n+1)\ZZ$, see 
\cite{dk:crossed}. The set $\Hom_{\Lambda_\oo}(\cm,\cn)$ can be identified  \cite{elmendorf} 
with the set of monotone maps $f: \ZZ\to\ZZ$ such that 
\be\label{eq:Lambda-oo-monot}
f(i+m+1)=f(i)+n+1.
\ee
The composition of morphisms is the composition of monotone maps. 
The map $\tau_n: i\mapsto i+1$
 is an invertible element of $\Hom_{\Lambda_\oo}(\cn, \cn)$ known as the
 {\em paracyclic rotation}.  The category $\Lambda$ can be obtained from
 $\Lambda_\oo$ by imposing additional relations $\tau_n^{n+1}=\Id$. 
 Thus we have a functor $p: \Lambda_\oo\to\Lambda$ bijective on objects
 and surjective on morphisms.
  
 We have an embedding  (functor bijective on objects and injective on morphsims)
 $\eps: \Delta\hra \Lambda_\oo$, $[n]\mapsto \cn$
 obtained by extending monotone maps $[m]\to[n]$ using \eqref{eq:Lambda-oo-monot}. 
 The set $\Mor(\Lambda_\oo)$ of morphisms of $\Lambda_\oo$ is generated
 by $\Mor(\Delta)$ and the automorphisms $\tau_n$. 
 The composition $p\eps: \Delta\to\Lambda$ is also an embedding. 
 
 By a {\em paracyclic}   (resp. {\em cyclic}) object in an ordinary category $\Cc$ we will mean
 a contravariant functor $X$  from $\Lambda_\oo$ (resp. $\Lambda$) to $\Cc$ and write $X_n=X(\cn)$. 
 Thus a (para)cyclic object can be seen as a simplicial object $X_\bullet$ together
 with extra {\em (para)cyclic symmetry} given by the action of $\tau_n$ on $X_n$ for each $n$ and commuting,
  in a certain definite way, cf.  \cite[\S 2.1]{dkss:1} with the  simplicial structure. 
  
  By a {\em (para)cyclic $\oo$-category} we will mean a (para)cyclic object $\C_\bullet$ in the category
  of simplicial sets such that each $\C_n$ is an $\oo$-category. 
  
  \paragraph{The S-construction of a stable $\oo$-category.} The Waldhausen S-construction of a
  stable (more generally, exact) $\oo$-category was originally defined in   \cite[\S 7.3]{DK:HSS}.
  In this paper we  modify the approach of \cite{DK:HSS} to get this and related constructions take values in 
  $\oo$-categories,  not just in  Kan complexes  
   (which appear as maximal $\oo$-subgroupoids in those
  $\oo$-categories).

  For $n\geq 0$ let $J(n)$ be the poset formed by pairs $0\leq i<j\leq n$ with order $(i,j)\leq (k,l)$
  iff $i\leq k$ and $j\leq l$.  As usual, we consider $J(n)$ as a category. 
  Thus $J(n)=\Mor [n]$  is the category of morphisms of the poset $[n]$ considered as a category. The correspondence $[n]\mapsto J(n)$ 
  defines a covariant functor $J: \Delta\to\Cat$. 
  
  \begin{defi} Let $\C$ be a stable $\oo$-category. We denote $S_n\C\subset\Fun(J(n),\C)$
  the full $\oo$-subcategory spanned by functors $Y: J(n)\to \C$ satisfying the conditions:
  \begin{itemize}
  \item[(1)] For each $0\leq i\leq n $  we have $Y(i,i)=0$ (is a zero object).
  \item[(2)] For each $0\leq i\leq k\leq j\leq l\leq n$ the square
  \[
  \xymatrix{
  Y(i,j)\ar[d] \ar[r]&Y(i,l)\ar[d]
  \\
  Y(k,j) \ar[r]& Y(k,l)
  }
  \]
  \end{itemize}
  is a biCartesian square in $\C$. When $n$ varies, the $\oo$-categories $S_n\C$ unite into a simplicial
  $\oo$-category $S_\bullet \C$ called the {\em Waldhausen S-construction} of $\C$. 
  \end{defi} 
  
  The following is well known \cite[Rem.4.1.2]{lurie:rotation}:
  
  \begin{prop}\label{prop:A_n-quiv}
  The projection $S_n(\C)\to\Fun(\Delta^{n-1}, \C)$ given by restriction to
  the subposet
  \[
  [n-1] \simeq \bigl\{ (0,1) <(0,2) <\cdots < (0,n)\bigr\} \,\subset \, J_n,
  \]
  is an equivalence of $\oo$-categories. In particular, $S_n(\C)$ is stable. \qed
  \end{prop}
  
  Note that $\Fun(\Delta^{n-1}, \C)$ can be seen as the category of representations
  of the $A_n$-quiver in $\C$.

If $\C=\N^\dg(\Vc)$ is the dg-nerve of a $2$-periodic pre-triangulated dg-category $\Vc$, then 
by  \cite{DK:triangulated}   $S_\bullet\C$ has a natural structure of a cyclic, not just a simplicial 
$\oo$-category. The cyclic symmetry $\tau_n$ is given by the 
(derived) Coxeter functor of the $A_n$-quiver. 
This was  extended to ``$2$-periodic'' stable $\oo$-categories by Lurie \cite {lurie:rotation}. Lurie 
further showed that for $\C$ an arbitrary stable $\oo$-category, $S_\bullet\C$ always has a
paracyclic structure, with $\tau_n^{n+1}=\Sigma^2$ being the functor of shift by $2$ (also see \cite{dyck:a1}). 
We will obtain this latter statement as a consequence of our more general results below. 

\paragraph{The relative $S$-construction of an exact functor.} Let $F: \A\to\B$ be an exact
functor of stable $\oo$-categories. For every $n\geq 0$ we have then the induced
$\oo$-functor $F_*: S_n\A\to S_n\B$. 
 Following \cite[Def. 1.5.4] {waldhausen},  
we define the $\oo$-category $S_n(F)$ as the pullback (in the category of simplicial sets)
\be\label{eq:S-n-F}
\xymatrix{
S_n(F)\ar[d] \ar[r]& S_{n+1}\B\ar[d]^{\del_{n+1}}
\\
S_n\A \ar[r]_{F_*}  &S_n\B. 
}
\ee
 When $n$ varies, the $S_n(F)$ unite into a simplicial $\oo$-category $S_\bullet(F)$
 which we call the {\em relative Waldhausen S-construction} of $F$.
 
 \begin{exas}\label{exas:S_n-F}
 (a) We have $S_0(F)=\B$, while $S_1(F) = \Fun_{\Delta^1}(\Delta^1, \Gamma(F))$
 is the category of sections of the covariant Grothendieck construction of $F$.
 In other words, an object of $S_1(F)$ consists of objects $a\in \A$, $b\in\B$ and
 a morphism $F(a)\to b$ in $\B$. Thus $S_1(F)$ is the stable $\oo$-category
 glued out of $\A$ and $\B$ using $F$ as a coCartesian gluing functor. 
 
 \vskip .2cm
 
 (b)  More generally, for any $n$ by applying Proposition
 \ref{prop:A_n-quiv} 
  we identify $S_n(F)$ with the category of sections of the covariant
  Grothendieck construction of the composite functor
  \[
 \Fun(\Delta^{n-1}, \A) \buildrel \ev_{n-1}\over\lra
  \A \buildrel F\over \lra\B,
  \]
  where $\ev_{n-1}$ is the functor of evaluation at the $(n-1)$st (last)
  vertex of $\Delta^{n-1}$. 
   In particular, $S_n(F)$
 is a stable $\oo$-category. 
 
 This  also means that  ``size-wise'',
$S_n(F)$  looks roughly  like the sum  of $n$ copies of $\A$ and one copy of $\B$. 
 For this reason $S_\bullet (F)$ can be seen as a categorical analog
 of the nerve of the Picard groupoid $[\Psi\buildrel b\over \to\Phi]$ associated to a $2$-term
 complex, i.e., a morphism of abelian groups $b: \Psi\to\Phi$. Indeed,
 as an abelian group,  $\N_n[\Psi\to\Phi] = \Psi^{\oplus n} \oplus\Phi$, cf. 
 \cite[Ex. 4.3.3]{dkss:1}. In fact, the relative S-construction is a particular case of the
 {\em categorified Dold-Kan nerve} of \cite{dyckerhoff:DK} applied to the $2$-term complex
 $\A\buildrel F\over \to\B$ of stable $\oo$-categories. 
 \vskip .2cm
 
 (c) Note that we could use any $\del_i$, $i=0,\cdots, n+1$ instead of $\del_0$
 to form the pullback diagram above, and it would lead to an equivalent construction
 of $S_n(F)$. This is because the paracyclic symmetry (Coxeter functor)
 $\tau_{n+1}: S_{n+1}\B\to S_{n+1}\B$ rotates the $\del_i$.  This phenomenon will
 be used later in the proof of Theorem \ref {thm:2-Seg}. 
  \end{exas}

   \paragraph{Relative S-construction in terms of SOD.}  Given $F: \A\to\B$ as above,
   we can construct a stable $\oo$-category $\C$
   with a coCartesian SOD $(\A, \B)$ 
   by using $F$ as a coCartesian gluing functor. Explicitly, 
   $\C = \Fun_{\Delta^1}(\Delta^1, \Gamma(F))$.  Let us describe the simplicial
   $\oo$-category $S_\bullet(F)$ in terms of  $\C$. This can be  compared with
    the ``simplified'' description of the relative S-construction given by Waldhausen
   in \cite[p. 344] {waldhausen}.  
   
   \begin{defi}\label{def:rel-pushf}
  (a)  Let $\C$ be a stable $\infty$-category equipped with a coCartesian  SOD
$(\A,\B)$. We say that a square 
\[
	\begin{tikzcd}
		a \ar{r}\ar{d}{u} & b\ar{d}\\
		a' \ar{r} & b'
	\end{tikzcd}
\]
with $a,a'$ objects of $\A$ and $b,b'$ objects of $\B$, is a {\em relative pushout square} for $(\A, \B)$,  if the
right-hand square of the induced rectangle
\[
	\begin{tikzcd}
		a \ar{r}{!} \ar{d}{u} & t\ar{d}{v} \ar{r}&  b\ar{d}\\
		a' \ar{r}{!} & t' \ar{r} & b'
	\end{tikzcd}
\]
is a biCartesian square in $\B$. 

(b) Dually, we define a {\em relative pullback square} with respect to a
Cartesian  SOD as a relative pushout square in the opposite category.
   \end{defi}

   Let now $[n]^\rhd = \{0,1,\cdots, n, \oo\}$ be the poset obtained from $[n]$
   by adjoining a maximal element $\oo$. As before,
   let $W(n)= \Mor ([n]^\rhd)$ be
the poset   formed by all pairs $(i\leq j)$ with $i,j\in [n]^\rhd$.
   It is clear that the various posets $W(n)$ organize into a cosimplicial object
\begin{equation}\label{eq:simplicial}
	\Delta \lra \Cat,\; [n] \mapsto W(n).
\end{equation}

\begin{defi}\label{def:S-n-C-A-B} 
Let $\C$ be a stable $\infty$-category equipped with a coCartesian semiorthogonal
	decomposition $(\A,\B)$   We define 
	$
		S_n(\C; (\A, \B)) \subset \Fun(W(n),\C)
	$
	as the full subcategory spanned by those diagrams
	$
		X: W(n) \lra \C
	$
	satisfying:
	\begin{enumerate}
		\item The objects $X(i,\oo)$, $0 \le i \le n$, lie in $\B$.
		\item The objects $X(i,j)$, $0 \le i \le j \le n$, lie in $\A$.
		\item For every $0 \le i \le \oo$, the object $X(i,i)$ is a zero object in
		 $\A$.
		\item For every $0 \le i \leq k \leq j \leq l \le n$, the square
			\[
				\begin{tikzcd}
					X(i, j) \ar{r}\ar{d} & X(i,l) \ar{d}\\
					X(k,j) \ar{r} & X(k,l)
				\end{tikzcd}
			\]
			is a biCartesian square in $\A$. 
		\item For every $0 \le i\le k\le j  \le n$, the square 
			\[
				\begin{tikzcd}
					X(i,j) \ar{r}\ar{d} &X(i,\oo)  \ar{d}\\
					X(k,j) \ar{r} & X(k,\oo)
				\end{tikzcd}
			\]
			is a relative pushforward square with respect to $(\A,\B)$.
	\end{enumerate}
	Via the functoriality from \eqref{eq:simplicial}, we obtain a simplicial $\infty$-category
	\[
		S_{\bullet}(\C;(\A,\B)): \Delta^{\op} \lra \Cat_{\infty}
	\]
	which we call the {\em relative S-construction} of $\C$ with respect to
	 $(\A,\B)$. 
\end{defi}

\begin{prop}
The simplicial $\infty$-category
$S_{\bullet}(\C;(\A,\B))$ is equivalent to $S_\bullet(F)$, where $F: \A\to\B$
is a coCartesian gluing functor for $(\A,\B)$. 
\end{prop}

 \noindent{\sl Proof:} To compare with the definition   \eqref{eq:S-n-F}
 of $S_\bullet(F)$, we
 recall the interpretation of
 the induced  rectangle in part (a) of Definition \ref{def:rel-pushf} in terms of
 $F$. Indeed, from the definition of $(\A,\B)$-coCartesian arrows,
 we have that the induced rectangle has
  $t=F(a)$, $t'=F(a')$ and $v=F(u)$. After this identification, 
   Definition \ref{def:S-n-C-A-B} becomes identical to
  the fiber product in
   \eqref{eq:S-n-F}: introducing
   the extra label $\oo$  corresponds to introducing the label $n+1$
   in the definition of $S_n(F)$ by putting $\del_{n+1}:S_{n+1}\B\to 
   S_n\B$
   as the right vertical arrow of the fiber product diagram. \qed
   
   
   \subsection{The relative S-construction of a spherical functor}
   \label{subsec:rel-S-sph}

   \paragraph{ The main result: paracyclic structure.} 
   We now formulate the main result of this paper whose proof will be finished  in \S \ref{par:spher-vs-rel-S}
   
   \begin{thm}\label{thm:main}
  (i)  Suppose that $F: \A\to\B$ is a spherical functor of stable $\oo$-categories.
   The  simplicial $\oo$-category $S_\bullet(F)$ has a natural lift
   to a paracyclic $\oo$-category.
   
   \vskip .2cm
   
   (ii) Equivalently, suppose that $(\A,\B)$ is a $4$-periodic SOD of a stable
   $\oo$-category $\C$. The simplicial $\oo$-category $S_{\bullet}(\C;(\A,\B))$
   has a natural lift
   to a paracyclic $\oo$-category. 
   \end{thm}
  
 It suffices to  prove the form (ii) of the theorem.  To do this, we give another definition of 
 $S_{\bullet}(\C;(\A,\B))$ whose paracyclic nature will be clear from the outset. 
 This is analogous to  the approach to the absolute S-construction in  \cite[\S4.3] {lurie:rotation}. 
 
 \paragraph{The spherical S-construction.} 
 We start with   extending each poset $W(n)$ as follows. 
   
\begin{defi}\label{def:R(n)}
Let $\ZZ^{\rhd}$ denote the poset obtained from $\ZZ$ by adjoining a maximal element  $\infty$. Similarly, 
let $\ZZ^{\lhd}$ be $\ZZ$ with a minimal element adjoined,
denoted by $-\infty$. Let 
$
	R(n) \subset \ZZ^{\lhd} \times \ZZ^{\rhd}
$
 be the poset consisting of: 
\begin{itemize}
	\item all pairs $(i,j)$ of integers such that $0 \le j-i \le n+1$,
	\item all elements of the form $(-\infty,j)$ where $j \in \ZZ^{\rhd}$,
	\item all elements of the form $(i,\infty)$ where $i \in \ZZ^{\lhd}$.
\end{itemize}
\end{defi}

\noindent It is evident that the various posets $R(n)$ organize into a coparacyclic object
\begin{equation}\label{eq:paracyclic}
		R: \Lambda_{\infty} \lra \Cat,\; \cn \mapsto R(n). 
\end{equation}

\begin{defi}\label{defi:sphericalS} Let $\C$ be a stable $\infty$-category equipped with a $4$-periodic  SOD  $(\A,\B)$ and set $\B' = \A^{\perp}$, $\A'={^\perp}\B$  so that every consecutive pair of the
	sequence
	\[
\cdots \quad 	, \quad 	\B'\quad,\quad \A\quad,\quad \B\quad , \quad \A' \quad , \quad \B' \quad , \quad \cdots
	\]
	forms an  SOD. We define 
	\[
	S_n^\sph(\C; (\A, \B))\, \subset\, \Fun(R(n), \C)
	\]
	to be the full subcategory spanned by those diagrams
	$
		X: R(n) \to \C
	$
	which satisfy the following conditions:
	\begin{enumerate}
		\item The objects $X(-\infty,j)$, $j \in \ZZ$ lie in $\B'$.
		\item The objects $X(i,j)$ for $-\infty < i \le j < \infty$ lie in $\A$.
		\item The objects $X(i,\infty)$, $i \in \ZZ$ lie in $\B$.
		\item For every $i \in \ZZ$, the object $X(i,i)$ is a zero object in $\C$. 
		\item The object $X(-\oo, \oo)$ is also a zero object in $\C$.
		\item For every $j \in \ZZ$, the morphism $X(-\infty,j) \to X(j-n-1,j)$ is $(\B',\A)$-coCartesian (so its fiber lies  in $\B$).
		\item For every $i \in \ZZ$, the morphism $X(i, i+n+1) \to X(i,\infty)$ is $(\A,\B)$-Cartesian (so its fiber lies in $\B'$). 
		\item For every $i \in \ZZ$, the square
			\[
				\begin{tikzcd}
					X(-\infty, i+n+1) \ar{r}\ar{d}{!} & X(-\infty,\infty)=0 \ar{d}\\
					X(i,i+n+1) \ar{r}{*} & X(i,\infty)
				\end{tikzcd}
			\]
			is a biCartesian square in $\C$.
		\item For every $-\infty < i\le k\le j \le l < \infty$ with $l-i \le n+1$, the square
			\[
				\begin{tikzcd}
					X(i, j) \ar{r}\ar{d} & X(i,l) \ar{d}\\
					X(k,j) \ar{r} & X(k,l)
				\end{tikzcd}
			\]
			is a biCartesian square in $\A$. 
		\item For every $-\infty < i \le k\le  j < \infty$ with $j-i \le n+1$, the square 
			\[
				\begin{tikzcd}
					X(i, j) \ar{r}\ar{d} & X(i,\infty) \ar{d}\\
					X(k,j) \ar{r} & X(k,\infty)
				\end{tikzcd}
			\]
			is a relative pushout square with respect to $(\A,\B)$.
		\item For every $-\infty < i \leq j\leq l < \infty$ with $l-i \le n+1$, the square 
			\[
				\begin{tikzcd}
					X(-\infty, j) \ar{r}\ar{d} &X(i,j)  \ar{d}\\
					X(-\infty,l) \ar{r} & X(i,l)
				\end{tikzcd}
			\]
			is a relative pullback square with respect to $(\B',\A)$.
	\end{enumerate}
	Via the functoriality from \eqref{eq:paracyclic}, we obtain a paracyclic $\infty$-category
	\[
		S^{\on{sph}}_{\bullet}(\C;(\A,\B)): \Lambda_{\infty}^{\op} \lra \Cat_{\infty}
	\]
	which we call the {\em spherical $S$-construction} of $\C$ with respect to $(\A,\B)$.

\end{defi}

	\paragraph{Spherical vs. relative S-construction.}\label{par:spher-vs-rel-S}
Theorem \ref{thm:main} in the form (ii) would follow from the next result.

\begin{thm}\label{thm:sphericalS} Let $\C$ be a stable $\infty$-category equipped with a $4$-periodic  SOD 
$(\A,\B)$. Then the functor given by restriction along the inclusions $W(n)
	\subset R(n)$ induces an equivalence
	\[
		S^{\on{sph}}_{\bullet}(\C;(\A,\B))|\Delta^{\op} \,\simeq \, 
		S_{\bullet}(\C;(\A,\B))
	\]
	between the simplicial object underlying the spherical $S$-construction and the relative
	$S$-construction of $\C$ with respect to $(\A,\B)$.
\end{thm}
 \noindent{\sl Proof:}
To show this, we use an elaboration of the argument in the proof of
\cite[Prop.4.3.3]{lurie:rotation}, given by iterated application of \cite[4.3.2.15]{lurie:htt} using
the various universal properties and Proposition \ref{prop:cartesian} to control those properties
via relative Kan extensions. Most crucial will be Lemma \ref{lem:grid} which implies
that the mutation of an $(\A,\B)$-relative pushout square is an $(\B',\A)$-relative pullback square
and vice versa.

We prove that for each $n\geq 0$ the functor 
\[
S^{\on{sph}}_{n}(\C;(\A,\B)) \lra S_{n}(\C;(\A,\B))
\]
given by restriction from $R(n)$ to
$W(n)$, is an equivalence. 

To this end, we represent $R(n)$
as the union of an infinite chain of inclusions of subposets
\be\label{eq:exhaustion}
W(n)=W_0(n) \subset W_1(n)\subset \cdots \subset R(n) = \bigcup_i W_i(n)
\ee
and, for each $i\in\ZZ_{\geq 0}$, define a full subcategory $S_n^{(i)}\subset
\Fun(W_i(n), \C)$ so that:
\begin{itemize}
\item[(1)] $S_n^{(0)}= S_n(\C; (\A,\B))$, and $S_n^\sph(\C; (\A, \B))$ is
the limit of the $S_n^{(i)}$. 

\item[(2)] Objects  of $S_n^{(i)}$, $i\geq 1$,  are precisely  left or right
(depending on $i$) Kan extensions of  those of  $S_n^{(i-1)}$. 
\end{itemize}
\noindent The statement about equivalence will then follow by iterated
application of \cite[4.3.2.15]{lurie:htt}, as all the
restrictions 
\[
S_n(\C; (\A,\B)) = S_n^{(0)} \lla S_n^{(1)} \lla \cdots\lla
S_n^\sph(\C; (\A, \B))
\]
will be shown to be equivalences. 

\vskip .2cm

We will define each $W_i(n)$ by adding a single element $(p_i,q_i)$  to $W_{i-1}(n)$,
and define $S_n^{(i)}$ as the category of diagrams $W_i(n)\to \C$ which
satisfy all those conditions of Definition \ref{defi:sphericalS} which make
sense within $W_i(n)$. The adding of the object $X(p_i,q_i)$ corresponding to
the new element  will be referred to as {\em filling}.  This object will
be determined from the previous part of the diagram either as some
pullback, in which case diagrams from $S_n^{(i)}$ will be right Kan extensions
of those from $S_n^{(i-1)}$, or as some pushout, in which case diagrams 
 from $S_n^{(i)}$ will be left Kan extensions. 

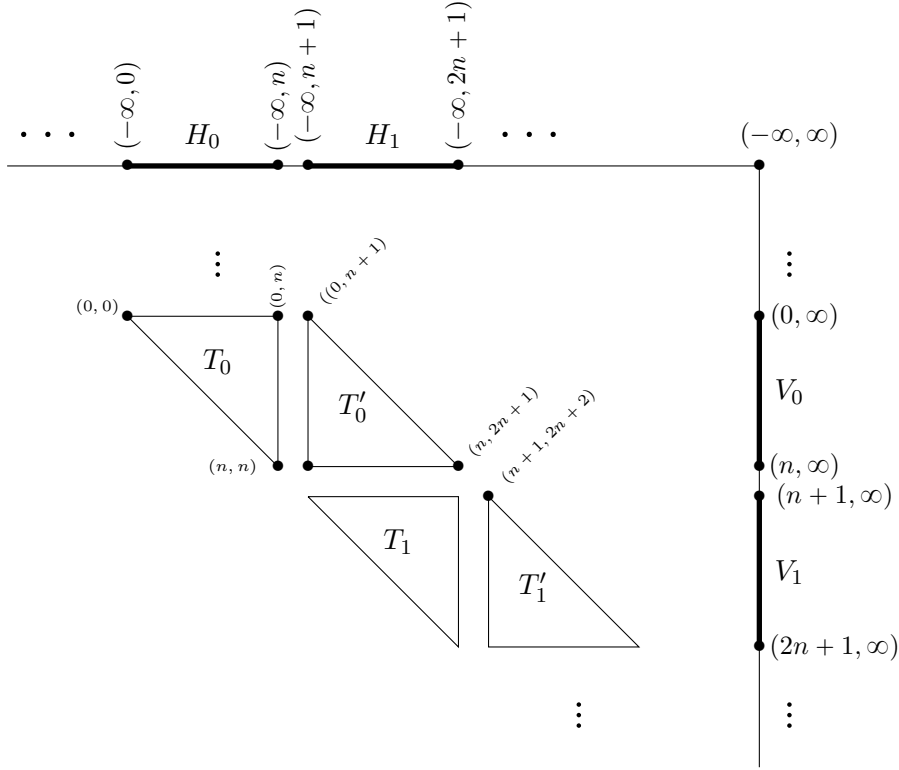
\begin{figure}[h]
\centering
\begin{tikzpicture}[scale=0.4]

\draw (-1,0) -- (-6,5) --(-1,5) --(-1,0); 
\draw (0,0) -- (0,5) -- (5,0) -- (0,0); 
\draw (0,-1) --(5,-1) -- (5, -6) -- (0,-1); 
\draw (6,-1) -- (6, -6) -- (11,-6) -- (6, -1); 

\draw (15, 10) -- (15, -10); 
\draw (-10, 10) -- (15, 10); 

\draw[line width=2] (15,0) -- (15,5); 
\draw[line width=2] (15,-1) -- (15,-6); 
\draw[line width=2] (-6,10) -- (-1,10); 
\draw[line width=2] (0,10) -- (5, 10); 

\node at (9,-8) {\huge$\vdots$}; 
\node at (-3, 7) {\huge$\vdots$}; 

\node at (-3,3.5){$T_0$}; 
\node at (1.5,2){$T'_0$}; 
\node at (3, -2.5){$T_1$}; 
\node at (7.5, -4){$T'_1$}; 

\node at (16, 2.5){$V_0$}; 
\node at (16, -3.5){$V_1$}; 

\node at (-3.5, 11){$H_0$}; 
\node at (2.5, 11){$H_1$}; 

\node at (-8.5, 11){\huge$\cdots$}; 
\node at (7.5, 11){\huge$\cdots$}; 

\node at (16, -8){\huge$\vdots$}; 
\node at (16,7){\huge$\vdots$}; 

\node at (15,5){$\bullet$}; 
\node at (15,0){$\bullet$}; 
\node at (15,-1){$\bullet$}; 
\node at (15,-6){$\bullet$}; 

\node at (-6,10){$\bullet$}; 
\node at (-1,10){$\bullet$}; 
\node at (0,10){$\bullet$}; 
\node at (5,10){$\bullet$}; 

\node at (16.5,5){\small$(0,\oo)$}; 
\node at (16.5,0){\small$(n,\oo)$}; 
\node at (17.5,-1){\small$(n+1,\oo)$}; 
\node at (17.5,-6){\small$(2n+1,\oo)$}; 

\node[rotate=90] at (-6,12){\small$ (-\oo, 0)$}; 
\node[rotate=90] at (-1,12){\small$ (-\oo, n)$}; 
\node[rotate=90] at (0,13){\small$ (-\oo, n+1)$}; 
\node[rotate=90] at (5,13){\small$ (-\oo, 2n+1)$}; 

\node at (-6,5){$\bullet$}; 
\node at (-7, 5.3){\tiny$(0,0)$}; 

\node at (15,10){$\bullet$}; 
\node at (16, 11){\small$(-\oo, \oo)$}; 

\node at (-1,5){$\bullet$};
\node[rotate=90] at (-1,6){\tiny$(0,n)$}; 

\node at (0,5){$\bullet$}; 
\node[rotate=45] at (1.5, 6.5) {\tiny$((0,n+1)$}; 

\node at (-1,0){$\bullet$}; 
\node  at (-2.5,0){\tiny$(n,n)$}; 

\node at (5,0){$\bullet$}; 
\node[rotate=45] at (6.5,1.5){\tiny$(n, 2n+1)$}; 

\node at (0,0){$\bullet$}; 

\node at (6,-1){$\bullet$}; 
\node[rotate=45] at (8, 1){\tiny$(n+1, 2n+2)$}; 

\end{tikzpicture}

\caption{Subdivision of $R(n)$ into triangles and intervals.}\label{fig:domains}
\end{figure}

To explain the order of adding elements, we subdivide the poset $R(n)$
into the following subsets  which we visualize as geometric figures 
 in Fig. \ref{fig:domains}:
\begin{itemize}

\item ``Triangles'' $T_m$, $m\in\ZZ$ consisting of $(i,j)$ with $m(n+1) \leq i\leq j < (m+1)(n+1)$. 

\item ``Triangles'' $T'_m$, $m\in \ZZ$, consisting of $(i\leq j)$ with $(m-1)(n+1) \leq i  < m(n+1)$ and $m(n+1) \leq j < (m+1)(n+1)$. 

\item ``Vertical intervals'' $V_m$, $m\in \ZZ$, consisting of $(i,\oo)$ with $m(n+1) \leq i <
(m+1)(n+1)$.

\item ``Horizontal intervals'' $H_m$, $m\in\ZZ$, consisting of $(-\oo,j)$
with $m(n+1) \leq j<
(m+1)(n+1)$.
\end{itemize}

\noindent Our initial set $W_0(n) = W(n)$ is the union of $T_0$ and $V_0$. 
We now start the filling procedure by filling the triangle $T'_0$ from top to bottom.

 \vskip .2cm

So we take $(p_1, q_1)=(0, n+1)$, the top vertex of $T'_0$. The corresponding object
$X(0, n+1)$ is found by factoring $X(0,n)\to X(0,\oo)$ as
 $X(0,n)\to X(0, n+1)\buildrel *\over\to
X(0,\oo)$ with the $*$-arrow being $(\A, \B)$-Cartesian. This amounts to a
right Kan extension. Alternatively, in terms of the gluing functors $F: \A\leftrightarrow \B:G$
classifying the SOD $(\A,
B)$,  we 
have $X(0, n+1)=G(X(0,\oo)$ and the arrow $X(0,n+1)\to X(0,\oo)$
 involving is  given by the counit of the
adjunction. 

Next, we fill the object straight down, at $(p_2, q_2)=(1, n+1)$. This object $X(1, n+1)$
is found to fit into  a biCartesian square
\[
\xymatrix{
X(0,n)\ar[d] \ar[r]& X(0, n+1)
\ar@{-->}[d]
\\
X(1,n) \ar@{-->}[r]& X(1, n+1),
}
\]
i.e., as a pushout. This amounts to a left Kan extension.  

Next, we fill the object at $(p_3, q_3)=(1, n+2)$ by factoring $X(1, n+1)\to X(1,\oo)$
as $X(1, n+1)\to X(1, n+2)\buildrel *\over\to X(1,\oo)$. 

Continuing like this, we will the entire triangle $T'_0$. We further fill $T_1$ in the same
way, row by row from the top by forming pushots and starting  each row from
$X(i,i)=0$.

\vskip .2cm

After this, we fill $X(n+1,\oo)$ in the vertical  interval $V_1$ from the condition that 
the left square below is an $(\A,\B)$-relative pushout square in $\C$ or, equivalently
(in terms of the gluing functor $F$)
the right square below is biCartesian in $\B$:
\[
\xymatrix{
X(n, 2n+1)\ar[d] \ar[r]& X(n,\oo)\ar@{-->}[d]
\\
X(n+1, 2n+1) \ar@{-->}[r]& X(n+1, \oo)
}
\quad\quad
\quad
\xymatrix{
F(X(n, 2n+1))\ar[d] \ar[r]& X(n,\oo)\ar@{-->}[d]
\\
F(X(n+1, 2n+1)) \ar@{-->}[r]& X(n+1, \oo). 
}
\]
This amounts to a left Kan extension. 

Continuing like this we can fill,  going down row by row, all the $T_m, T'_m, V_m$ with
$m\geq 0$.  For this, we do not yet need the fact that our adjunction was spherical
or, equivalently, that $(\A,\B)$ was $4$-periodic. 

\vskip .2cm

To fill the entire diagram, we need to alternate the above ``downward moves'' with
the ``upward moves'' in the opposite direction which are as follows.

\vskip .2cm

Note that condition (8) of Definition \ref{defi:sphericalS} recovers the objects in $H_1$
as mutations of the corresponding objects in $V_0$, the objects in $H_2$ as
mutations of those in $V_1$ etc. 
After doing this (and here we use $4$-periodicity), Lemma \ref{lem:grid}  implies that the mutation of an $(\A,\B)$-relative pushout square
is an $(\B',\A)$-relative pullback square and vice versa. Thus the diagram with
the objects in $H_{\geq 1}$ filled in,  will satisfy all  the relevant conditions of 
Definition \ref{defi:sphericalS}. 

After this, already having $X(-\oo, n+1)$ (from $H_1$), we fill $X(-\oo, n)$
(from $H_0$ which is not yet filled) from the condition that it fits into $(\B', \A)$-relative
pullback square 
\[
\xymatrix{
X(-\oo, n) \ar@{-->}[r] \ar@{-->}[d]& X(-\oo, n+1)\ar[d]
\\
X(0,n) \ar[r]&X(0, n+1).
}
\]
 This amounts to a right Kan extension. 
 
 \vskip .2cm
 
 After this, we fill $H_0$ and the bottom row of $T'_{-1}$ from right to left.
 
 That is, we find $X(-1,n)$, the object at the right end of the bottom row of
 $T'_{-1}$ by factoring $X(-\oo,n)\to X(0,n)$ as $X(\oo,n) \buildrel !\over\to X(-1,n)
 \to X(0,n)$, which amounts to a left Kan extension. 
 
 Then we find $X(-1, n-1)$ to fit into a biCartesian square
 \[
 \xymatrix{
 X(-1, n-1)\ar@{-->}[d] \ar@{-->}[r]& X(-1,n)\ar[d]
 \\
 X(0,n-1)\ar[r]& X(0,n), 
 }
 \]
 a right Kan extension, and so on. 
 
 \vskip .2cm
 
 By alternating the two types of moves, downward and upward,  we exhaust the set $R(n)$
 by a single sequence to get a chain of inclusions as in  
 \eqref{eq:exhaustion}. 
 \qed

\paragraph{Examples: the spherical $S_0$ and $S_1$.} 
 Let us describe the low-dimensional cells of the spherical $S$-construction  more explicitly and illustrate
the action of the paracyclic shift. 

\begin{exa}
The $\infty$-category $S_0^\sph(\C, (\A, \B))$ consists of diagrams of the form
\begin{equation}\label{eq:r0}
	\adjustbox{scale=.8,center}{%
		\begin{tikzcd}
		{\mathbf{ \cdots}} \ar{r}& X(-\infty,0) \ar{r}\ar{dd}{!} & X(-\infty,1)\ar{ddd}{!} \ar{r} &
		X(-\infty,2)\ar{dddd}{!}
		\ar{r} & \dots \ar{r} & X(-\infty,\infty) \ar{d} \\
		& & & & & \vdots \ar{d}\\
		& X(-1,0) \ar{d}\ar{rrrr}{*} &   &  & &X(-1,\infty) \ar{d} \\
		& X(0,0) \ar{r} & X(0,1)\ar{d} \ar{rrr}{*}    &  & &X(0,\infty) \ar{d} \\
		&             & X(1,1) \ar{r}  & X(1,2) \ar{rr}{*}   &  &X(1,\infty) \ar{d}\\
		&&&&&{\mathbf{ \vdots}}
	\end{tikzcd}
	}
\end{equation}
satisfying the conditions of Definition \ref{defi:sphericalS}. By the argument of the proof
of Theorem \ref{thm:sphericalS}, the projection map onto the subdiagram given by $X(0,\infty)$ is an equivalence, 
showing that $S^{\on{sph}}_0(\C;(\A,\B)) \simeq \B$. To exhibit the action of
the paracyclic shift, which simply acts by shifting all indices by $+1$, we restrict attention to
the the subdiagram 
\[
	\begin{tikzcd}
		X(0,1) \ar{r}{*}\ar{d} & X(0,\infty) \ar{d}\\
		X(1,1)=0 \ar{r} & X(1,\infty)
	\end{tikzcd}
\]
of \eqref{eq:r0} which, by the conditions of Definition \ref{defi:sphericalS}, is an
$(\A,\B)$-relative pushout square. Refining this square to the rectangle
\[
	\begin{tikzcd}
		X(0,1) \ar{r}{!}\ar{d} & Y \ar{r}\ar{d}&  X(0,\infty) \ar{d}\\
		X(1,1)=0 \ar{r}{!} & Y'=0 \ar{r} & X(1,\infty)
	\end{tikzcd}
\]
this means that the right-most square is biCartesian in $\B$ (with $Y'=0$ since $X(1,1)=0$).
 But this is simply saying that the
diagram exhibits $X(1,\infty)$ as $T_{\B}(X(0,\infty))[1]$, since the map $Y \to X(0,\infty)$ is a
counit map for the adjunction classified by the biCartesian  SOD  $(\A,\B)$.
Therefore, the paracyclic shift acts via the suspended spherical cotwist $T_{\B}[1]$.
\end{exa}

\begin{exa}
	The $\infty$-category  $S_1^\sph(\C; (\A,\B))$  consists of diagrams of the form
	
	\[
	\adjustbox{scale=.8,center}
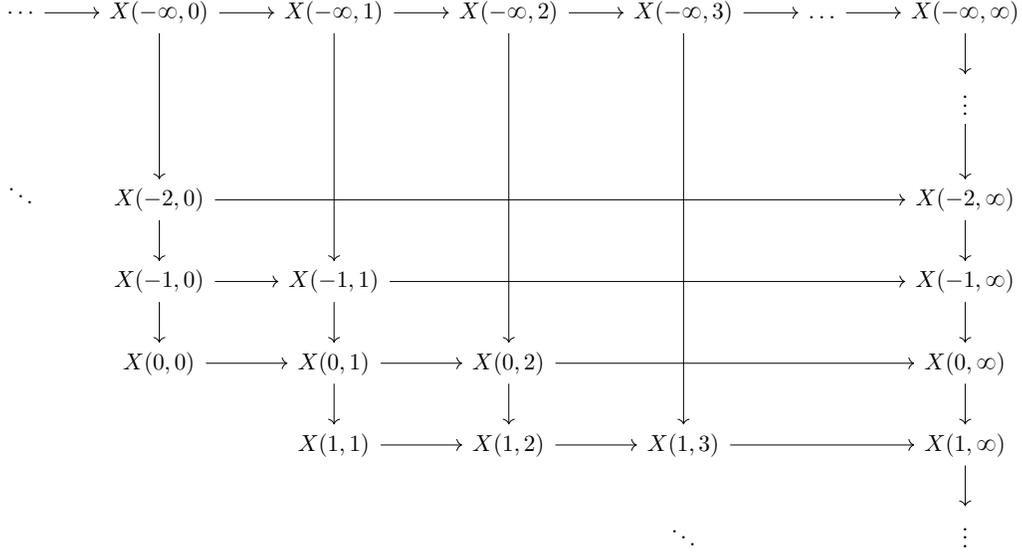
{%
		\begin{tikzcd}
			\cdots\ar{r} & X(-\infty,0) \ar{r}\ar{dd} & X(-\infty,1)\ar{ddd} \ar{r} & X(-\infty,2)\ar{dddd}
			\ar{r} & X(-\infty,3) \ar{ddddd}
			\ar{r} & \dots \ar{r} & X(-\infty,\infty) \ar{d} \\
			& & & & & & \vdots \ar{d}\\
			\ddots & X(-2,0) \ar{rrrrr}\ar{d} &  &  &  & &X(-2,\infty) \ar{d} \\
			& X(-1,0) \ar{d}\ar{r} & X(-1,1)\ar{d} \ar{rrrr}  &  &  & &X(-1,\infty) \ar{d} \\
			& X(0,0) \ar{r} & X(0,1)\ar{d} \ar{r}  & X(0,2)\ar{d} \ar{rrr}  &  & &X(0,\infty) \ar{d} \\
			&    & X(1,1) \ar{r}  & X(1,2) \ar{r}  & X(1,3) \ar{rr} &  &X(1,\infty)\ar{d} \\
			&&&&{\mathbf{\ddots}}&&{\mathbf{ \vdots}}
		\end{tikzcd}
		}
		\]
	satisfying the conditions of Definition \ref{defi:sphericalS}. The projection functor to
	the subdiagram 
	\[
	 		X(0,1) \lra X(0,\infty)
	\]
	defines a trivial Kan fibration thus providing an equivalence $S^{\on{sph}}_1(\C;(\A,\B)) \simeq
	\{\A,\B\} \overset{\fib}{\simeq} \C$. To unravel the action of the paracyclic twist, consider the subdiagram 
	\[
		\begin{tikzcd}
			X(0,1) \ar{r}\ar{d} & X(0,2) \ar{r}{*}\ar{d} & X(0,\infty) \ar{d}\\
			X(1,1) \ar{r} & X(1,2) \ar{r} & X(1,\infty)
		\end{tikzcd}
	\]
	which we may refine to a cubical diagram 
	\[
	\adjustbox{scale=.8,center}{%
		\begin{tikzcd}
			 	& Y \ar[rr]\ar[dd] & & X(0,\infty) \ar[dd]\\ 
				X(0,1) \ar{ur}{!} \ar[rr]\ar[dd] &   &  X(0,2) \ar[dd]\ar{ur}{*} & \\
				& Y' \ar[rr] & &  X(1,\infty)\\ 
				X(1,1) \ar{ur}{!} \ar[rr]&   & X(1,2) \ar[ur] &  
		\end{tikzcd}
		}
	\]
	Since the front and back faces of this cube are biCartesian squares, the cofiber of the
	morphism $X(1,2) \to X(1,\infty)$ in $\C$ is equivalent to the totalization of top
	commutative square
	\[
		\begin{tikzcd}
			X(0,1) \ar{r}{!} \ar{d} &  Y \ar{d}\\
			X(0,2) \ar{r}{*} & X(0,\infty)
		\end{tikzcd}
	\]
	in $\C$. But therefore, by inspection of \eqref{eq:cubical} and the subsequent definitions,
	the fiber of $X(1,2) \to X(1,\infty)$ is precisely the image of $\fib(X(0,1) \to
	X(0,\infty))$ under the relative suspension functor $\tau$ from Definition \ref{def:rel-susp}. In
	conclusion, the action of the paracyclic shift on $S^{\on{sph}}_1(\C;(\A,\B)) \simeq \C$ is
	given  by the relative suspension functor $\tau$. 
	In the case when $\B = 0$, so that the relative $S$-construction simply becomes the absolute
	$S$-construction of $\A$, we have $\tau = [1]$ in agreement with the known paracyclic
	structure in this case.
\end{exa}


\subsection{Spherical functors and admissible paracyclic stable $\infty$-categories}\label{subsec:sph-admis}

In this section we note some important additional features of the relative S-construction $S_\bullet(F)$ (for any exact $F$, spherical or not). 
We conjecture that these features, together with the paracyclic structure
which is specific to the spherical case, characterize  the   spherical S-constructions,

\paragraph{Admissibility: faces adjoint to degenerations.} \label{par:admis}
 Suppose that $\E_{\bullet}: \Delta^{\op} \to \Cat_\oo$ is a simplicial  $\infty$-category, i.e., a simplicial object in $\Set_\Delta$ such that each $\E_n$ is
an $\oo$-category. Then the
simplicial identities imply the following relations between the face and degeneracy maps  relating
$\E_{n-1}$ and $\E_n$: 
\be\label{eq:face-deg-rels}
\del_0s_0\buildrel =\over\to  \Id_{\E_{n-1}}, \quad 
\Id_{\E_{n-1}}\buildrel =\over\to  \del_1 s_0,\quad
\del_1 s_1\buildrel =\over\to \Id_{\E_{n-1}},
\quad \cdots , \quad
\Id_{\E_{n-1}} \buildrel =\over\to   \del_n s_{n-1}. 
\ee

\begin{defi}\label{def:admis}
 We say that $\E_{\bullet}$ is {\em admissible}, if the above identity maps form counit and
	units, respectively, of adjunctions
	\[
		\del_0 \dashv s_0 \dashv \del_1 \dashv \dots \dashv s_{n-1} \dashv \del_n.
	\]
	Analogously, we call a paracyclic  $\infty$-category admissible if the chain of adjunctions continues
	infintely on both sides. We denote $(\Cat_\oo)_\Delta^\adm$,
	resp. $(\Cat_\oo)_{\Lambda_\oo}^\adm$ the full subcategory in
	$(\Cat_\oo)_\Delta$, resp. $(\Cat_\oo)_{\Lambda_\oo}$
	spanned by admissible simplicial, resp. paracyclic $\oo$-categories. 
\end{defi}

\begin{prop}
(a) Let $F: \A\to\B$ be an exact functor of stable $\oo$-categories.
The relative S-construction $S_\bullet(F)$ is an admissible
simplicial stable $\oo$-category. 

\vskip .2cm

(b) Let $\C$ be a stable $\infty$-category equipped with a $4$-periodic semiorthogonal
	decomposition $(\A,\B)$. The spherical S-construction
 $S^\sph_\bullet(\C; (\A, \B))$ is an admissible paracyclic stable $\oo$-category. 
\end{prop}

\noindent {\sl Proof:} It is enough to show (a), since (b) then follows, in virtue of Theorem \ref{thm:main},
by applying an appropriate paracyclic rotation to bring any desired pair $(\del_i, s_i)$ or $(s_i, \del_{i+1})$
of paracyclic face and degeneracy operators into the simplicial range. 

So let us focus on (a) and consider the statement about the adjunction $(\del_i, s_i)$.
 By Proposition \ref{prop:oo-unit},  it suffices 
to prove that for any
two objects $X\in S_n(F)$, $Y\in S_{n-1}(F)$, 
the map of the mapping spaces
\[
 \eta: \Map_{S_n(F)}(X, s_i Y) 
 \buildrel \simeq\over\lra \Map_{S_{n-1}(F)}(\del_i X, Y). 
\]
induced by the   identity map $\del_is_i\buildrel =\over \to\Id$  from \eqref{eq:face-deg-rels},  
 is an equivalence.  But  $s_i$ consists in adding zeroes to the diagram $Y$,
 and we can think of this added $0$ as  an initial and final object in the strict sense. 
 With this understanding,
 both the source and the target of $\eta$ involve  exactly the same data, so $\eta$ is an identification
  ``on the nose''.  The statement about the adjunction $(s_i, \del_{i+1})$ is similar. \qed

\vskip .2cm

The patterns of adjunctions of Definition \ref{def:admis}
 can be implemented universally by
upgrading the categories $\Delta$ and $\Lambda_\oo$ to
(classical, strict) $2$-categories $\DDelta$, $\LLambda_{\infty}$
in which the corresponding coface and codegeneration maps become
$1$-morphisms equipped with canonical adjunction data (unit and counit
$2$-morphisms). 

\vskip .2cm

Following \cite[\S3.1.3]{dyckerhoff:DK}, we define 
  $\DDelta$ to have the same objects (finite ordinals $[m]$) and $1$-morphisms
  (monotone maps) as $\Delta$.  For two
  monotone maps $f, f': [m]\to [n]$, there is a unique $2$-morphism
  $f\Rightarrow f'$ if and only if $f\leq f'$ pointwise. 
  In other words, $\DDelta$ is obtained by regarding each poset $[m]$
  as a category and considering $\Hom_\Delta([m], [n])$, the set of
  monotone maps, not just as a set but as a poset, therefore  a category, 
  and this category is $\Hom_\DDelta([m], [n])$. In $\DDelta$, 
  we have  ``internal''  adjunctions corresponding to  those of 
  Definition \ref{def:admis}.    In each unit-counit pair of these adjunctions, 
  one member, corresponding to 
   \eqref{eq:face-deg-rels},  is  an identity $2$-morphism but the
  other one is a $2$-morphism represented by a nontrivial inequality $f \leq f'$. 
 Passing to nerves, we may 
interpret $\DDelta$ as enriched in $\Cat_{\infty}$. 
  
 \begin{rem}
 A more ``relaxed'' (and augmented) version of the $2$-category $\DDelta$ can be obtained
 as $\Hom_{\Adj^{\leq 3}}(a,a)$, where
$\Adj^{\leq 3}$  is  the Gray $3$-category generated by the universal
 $2$-dimensional adjunction, see  \cite { macdonald-scull, macdonald-stone}
 and 
 Remark \ref{rem:adj-sphere} above. This $2$-category can be seen
 as generated by the {\em universal $2$-comonad}, just as the augmented
 simplex category $\Delta_+$ is generated by the universal comonad
 in the classical sense, see \cite{{auderset}, macdonald-scull}. 
  \end{rem}
 
 Similarly, to define $\LLambda_\oo$, we interpret $\Hom_{\Lambda_\oo}(\cm, \cn)$
 as the set of quasi-periodic monotone maps $\ZZ\to\ZZ$, see 
 \eqref {eq:Lambda-oo-monot}, make it into a poset by pointwise
 inequality and interpret it as a category $\Hom_{\LLambda_\oo}(\cm, \cn)$. 
We regard  $\LLambda_\oo$ as enriched in $\Cat_{\infty}$ as well. 

Note that the category $\Cat_\oo$   can  itself be considered as enriched in 
$\Cat_{\infty}$. 

 \begin{defi}
 A {\em 2-simplicial} (resp. {\em $2$-paracyclic}) $\oo$-category
 is a $\Cat_{\infty}$-enriched
functor from $\DDelta^{(\op,-)}$ (resp. from $\LLambda_\oo^{(\op, -)}$) to $\Cat_\oo$. We denote by $(\Cat_\oo)_{\DDelta}$ (resp. 
$(\Cat_\oo)_{\LLambda_\oo}$)
the $\oo$-category formed by $2$-simplicial, resp. $2$-paracyclic
$\oo$-categories. 
 \end{defi}
 
 Thus a $2$-simplicial $\oo$-category is equipped with the
 ``missing'' adjunction  natural transformations completing each identity transformation in 
 \eqref{eq:face-deg-rels} to a unit-counit pair. Similarly for $2$-paracyclic $\oo$-categories.

\begin{prop} The forgetful functors
	$
		(\Cat_\oo)_{\DDelta} \to (\Cat_\oo)_{\Delta}
	$
	and 
	$
		(\Cat_\oo)_{\LLambda_{\infty}} \to (\Cat_\oo)_{\Lambda_{\infty}}
	$
	take values in the full subcategories of admissible objects. 
\end{prop}

\noindent {\sl Proof:} This follows from the property, discussed above,
that in $\DDelta$ as well as in $\LLambda_\oo$ the corresponding
members of the unit-counit pairs are  identity $2$-morphisms. \qed

\begin{conj} The forgetful functors
	$
		(\Cat_\oo)_{\DDelta} \to  (\Cat_\oo)_{\Delta}^{\adm}
	$
	and 
	$
		(\Cat_\oo)_{\LLambda_{\infty}} \to (\Cat_\oo)_{\Lambda_{\infty}}^{\adm}
	$
	induce equivalences of $\infty$-categories. 
\end{conj}

In other words, the datum of a $2$-simplicial (resp. $2$-paracyclic) lift of a given simplicial
(resp. paracyclic)   $\infty$-category is a property and not an extra datum.  

 For example, as already pointed out in
 Example \ref{exas:S_n-F}(b), the relative S-construction is a particular case of the categorified
 Dold-Kan nerve of \cite{dyckerhoff:DK}. It was established in  \cite{dyckerhoff:DK}
 that the relative Dold-Kan nerve always comes with a natural $2$-simplicial structure.


\paragraph{$2$-Segal property.} 
Recall from \cite  {DK:HSS} the concept of a $2$-Segal simplicial object.
For $n\geq 2$ let
\[
P_n\,=\,\Conv\bigl\{ v_n = e^{2\pi i k/{(n+1)}}, \quad k=0,1,\cdots, n\bigr\} \,\subset
\CC
\]
be the regular convex $(n+1)$-gon with vertices labelled by $0,1,\cdots, n$
counterclockwise. Each triangulation $\Tc$ of $P_n$ into straight triangles
lifts to a $2$-dimensional simplicial subset $\Delta^\Tc\subset\Delta^n$
by lifting each triangle $[v_i, v_j, v_k]$ into the $2$-simplex of $\Delta^n$
with vertices $i,j,k$. 

Let $\Cb$ be a combninatorial model category (see \cite[\S4.1]{DK:HSS}
and references therein) and $Y_\bullet$ be a simplicial object in $\Cb$. 
The simplicial subset $\Delta^\Tc\subset\Delta^n$ above gives rise to the
{\em derived membrane space} $Y_\Tc = (\Delta^\Tc, Y_\bullet)_R \in\Ob(\Cb)$
which is the homotopy limit (fiber product) in $\Cb$ of the diagram formed
by copies of $Y_i$ associated to each non-degenerate $i$-simplex of
$\Delta^\Tc$, $i=1,2$. It comes with a natural morphism $f_\Tc: Y_n\to Y_\Tc$.
The simplicial object $Y_\bullet$ is called {\em $2$-Segal}, if $f_\Tc$
is an equivalence for each $n\geq 2$ and each triangulation $\Tc$ of $P_n$. 

We apply this to the {\em Joyal model structure} on the category $\Set_\Delta$
of simplicial sets \cite{cisinski}. Fibrant objects in this structure are 
$\oo$-categories and, more generally, fibrations are so-called inner
fibrations of simplicial sets, see also  \cite[\S2.3]{lurie:htt}. Equivalences between fibrant
objects are equivalences of $\oo$-categories. 

Let $F: \A\to\B$ be an exact functor of stable $\oo$-categories. The relative
S-construction $S_\bullet(F)$, being a simplicial $\oo$-category is, in
particular, a simplicial object in $\Set_\Delta$. 

\begin{thm}\label{thm:2-Seg}
$S_\bullet(F)$ is $2$-Segal with respect to the Joyal model
structure. 
 \end{thm}
 
 \begin{rems}
 (a) This result generalizes  and refines \cite[Th.7.3.3]{DK:HSS}
 (specialized to the case of stable $\oo$-categories). That is, we consider the 
 relative and not just the absolute S-construction as in \cite{DK:HSS}.
 Further, we consider each $S_n(F)$ as an $\oo$-category rather than
 taking the maximal $\oo$-subgroupoid (i.e., the maximal Kan subcomplex)
 there as in \cite{DK:HSS}.
 
 \vskip .2cm
 
 (b) In the interpretation of $S_n(F)$ as the Fukaya category of the disk
 with coefficients in a schober given by $F$, see the Introduction, 
 the $2$-Segal property corresponds to the descent condition
 such as holds, e. g., for cohomology with coefficients in a perverse sheaf,
 cf. \cite{dkss:1}. 
 \end{rems}
 
 \noindent{\sl Proof of Theorem \ref{thm:2-Seg}:}
 We apply the path space criterion \cite[Th.6.3.2]{DK:HSS} reducing the
 $2$-Segal property of a simplicial object $Y_\bullet$ to the much
 simpler $1$-Segal (or simply Segal) property of the two {\em path space
 objects}
 \[
 P^\lhd Y_\bullet \,=\, i_\lhd^*\, Y_\bullet, \quad P^\rhd Y_\bullet \,=\,
 i_{\rhd}^*\, Y_\bullet. 
 \]
 Here $i_\lhd, i_\rhd: \Delta\hra\Delta$ are two self-embeddings of $\Delta$
(considered as the category of nonempty finite ordinals) given by adding an extra minimal, resp. maximal element to each ordinal. 
Thus the $n$-simplices and the face maps of $ P^\lhd\,  Y_\bullet,
 P^\rhd Y_\bullet$ are given by
 \[
 (P^\lhd Y_\bullet)_n\,=\,  (P^\rhd Y_\bullet)_n\,=\, Y_{n+1},\quad 
 \del_i^{P^\lhd Y_\bullet}=\del_{i+1}^Y, \,\,\, 
  \del_i^{P^\rhd Y_\bullet}=\del_{i}^Y,\,\,\, i=0,\cdots, n. 
 \]
 Recall that a simplicial object $Z_\bullet$ in $\Cb$ is called  $1$-Segal,
 if for any $n\geq 1$ the natural morphism
 \[
 f_n: Z_n\lra Z_1\times^R_{Z_0} Z_1\times^R_{Z_0} \cdots 
 \times^R_{Z_0}  Z_1
 \]
 to the $n$-fold homotopy fiber product is an equivalence. 
 
 We apply the criterion to $Y_\bullet=S_\bullet(F)$ and first look at $P^\rhd Y_\bullet$.  The object $(P^\rhd S_\bullet(F))_n$ being, by definition,
 nothing but $S_{n+1}(F)$, 
  Example \ref{exas:S_n-F}(b) identifies it with the category of sections
  of the covariant Grothendieck construction of the
  functor
  \[
  F(\ev_n): \Fun(\Delta^n, \A)\lra \B. 
  \]
 Further, in this identification  the simplicial structure on $P^\rhd S_\bullet(F)$ 
 is  given by the functoriality of $\Delta^n$ in $[n]$. This means
 that the target of the morphism $f_n$ for $Z_\bullet =P^\rhd S_\bullet(F)$
 is identifed with the category of sections
  of the covariant Grothendieck construction of the
  functor
\[
F(\ev'_n): \Fun\bigl(\Delta^1 \amalg_{\Delta^0} \Delta^1 \amalg_{\Delta^0} \dots \amalg_{\Delta^0}
		\Delta^1, \,\A \bigr) \lra \B,
\]
where $\ev'_n$ is the functor of evaluation at the last vertex of the chain
of  $(\Delta^1)$'s. The morphism $f_n$ itself is induced by restriction
with respect to the embedding
\[
\Delta^1 \amalg_{\Delta^0} \Delta^1 \amalg_{\Delta^0} \dots \amalg_{\Delta^0}
		\Delta^1 \, \subset  \, \Delta^n. 
\] 
 This embedding being an inner-anodyne morphism \cite[Def. 2.0.0.3]{lurie:htt},
 it follows that $f_n$ is equivalence. This shows that
 $P^\rhd S_\bullet(F)$ is $1$-Segal.
 
 \vskip .2cm
 
 To analyze the other path space $P^\lhd$, we first use a seemingly different model
 of the relative S-construction.  For $n\geq 0$ define the $\oo$-category 
 $S'_n(F)$ by the fiber product as in \eqref{eq:S-n-F} but with $\del_0$, not
 $\del_{n+1}$ as the right vertical arrow. We unite the various $S'_n(F)$
 into a simplicial $\oo$-category $S'_\bullet(F)$ by using the simplicial
 structure on the $S_{\bullet+1}(\B)$, the top right corners of the diagrams
  \eqref{eq:S-n-F}, via the self-embedding $i_\lhd: \Delta\hra\Delta$.
  This is, in fact, the original definition of the relative S-construction in
  \cite{waldhausen}. Now, Example  \ref{exas:S_n-F}(b) is immediately
  modified to this case. It identifies $(P^\lhd S'_\bullet(F))_n= S'_{n+1}(F)$
 with the category of sections of the {\em contravariant} Grothendieck 
 construction of the functor
 \[
 F(\ev_0): \Fun(\Delta^n, \A) \lra B,
 \]
where $\ev_0$ is the evaluation at the $0$th vertiex of $\Delta^n$.
The simplicial structure on $P^\lhd S'_\bullet(F)$ is again given by
the functoriality of $\Delta^n$ in $[n]$. So the same argument as before
shows that $P^\lhd S'_\bullet(F)$  is $1$-Segal.

\vskip .2cm

It remains to notice that the simplicial $\oo$-categories $S_\bullet(F)$
and $S'_\bullet(F)$ are in fact equivalent. The equivalence is obtained
by using the paracyclic structure on $S_\bullet(\B)$
(given by   \cite[Prop.4.3.3]{lurie:rotation} or applying our
Theorem \ref{thm:sphericalS} to the spherical functor $\B\to 0$). 
This structure gives, for each $n\geq 0$, the paracyclic rotation
$\tau_{n+1}: S_{n+1}(\B)\to S_{n+1}(\B)$ which sends $\del_{n+1}$
to $\del_0$ and shifts the other $\del_i$.

This identification $S_\bullet(F)\simeq S'_\bullet(F)$ shows that both
$P^\rhd$ and $P^\lhd$ of this simplicial $\oo$-category are $1$-Segal, so
it is $2$-Segal. \qed

\paragraph{Combining the properties.}  
Combining all the properties established in this and the previous sections,
we arrive at:

\begin{cor}
Let $F:  \A\to\B$
be a spherical functor of stable $\oo$-categories. 
Then the relative $S$-construction
	$S_\bullet(F)$ is an admissible $2$-Segal paracyclic stable
	$\infty$-category. \qed
\end{cor}

 Let $\Sph$ be the $\oo$-category formed by spherical adjunctions, and
$\St_{\Lambda_{\infty}}^{\adm, \text{$2$-Segal}}$ be
the $\oo$-category formed by stable $\infty$-categories and admissible $2$-Segal paracyclic 
	stable $\infty$-categories.

\begin{conj}\label{conj:sphad=adm2Segpar}
 The relative $S$-construction induces an equivalence
	$
		\Sph \lra \St_{\Lambda_{\infty}}^{\adm, \text{$2$-Segal}}.
	$
	 
\end{conj}


\vskip 1cm

\small{
 T.D.:  Universit\"at Hamburg,
Fachbereich Mathematik,
Bundesstrasse 55,
20146 Hamburg, Germany. Email: 
{\tt tobias.dyckerhoff@uni-hamburg.de}

\smallskip

M.K.: Kavli IPMU, 5-1-5 Kashiwanoha, Kashiwa, Chiba, 277-8583 Japan. Email: 
\hfil\break
{\tt mikhail.kapranov@protonmail.com}

\smallskip

 V.S.: Institut de Math\'ematiques de Toulouse, Universit\'e Paul Sabatier, 118 route de Narbonne, 
31062 Toulouse, France. Email: 
 {\tt schechtman@math.ups-tlse.fr }

 \smallskip 
 
 Y.S.:  Dept.  Math.,  Kansas State University, Manhattan, KS 66506 USA.
 Email: {\tt soibel@math.ksu.edu}  
 
   }          

\end{document}